\documentclass[12pt,oneside]{imsart}
\usepackage[utf8]{inputenc}
\RequirePackage{amsthm,amsmath,amsfonts,amssymb}
\RequirePackage[numbers,sort&compress]{natbib}
\RequirePackage[colorlinks,
citecolor=blue,
urlcolor=blue,
breaklinks]{hyperref}
\RequirePackage{graphicx}

\numberwithin{equation}{section}
\theoremstyle{plain}

\newtheorem{theorem}{Theorem}[section]
\newtheorem{lemma}[theorem]{Lemma}
\newtheorem{corollary}[theorem]{Corollary}
\newtheorem{proposition}[theorem]{Proposition}
\newtheorem{definition}[theorem]{Definition}
\theoremstyle{remark}

\newtheorem{remark}[theorem]{Remark}
\usepackage{indentfirst}
\usepackage{txfonts}

\newcommand{\red}{\color{red}}

\usepackage{geometry}
\usepackage{float}
\usepackage{arydshln}
\usepackage{booktabs}
\usepackage{multirow,makecell}
\usepackage{enumerate}

\def\bt{\begin{theorem}}
\def\et{\end{theorem}}
\def\bl{\begin{lemma}}
\def\el{\end{lemma}}
\def\bd{\begin{definition}}
\def\ed{\end{definition}}
\def\bc{\begin{corollary}}
\def\ec{\end{corollary}}
\def\bp{\begin{proposition}}
\def\ep{\end{proposition}}
\def\br{\begin{remark}}
\def\er{\end{remark}}

\def\bB{{\mathbf B}}
\def\bC{{\mathbf C}}
\def\bD{{\mathbf D}}

\def\bH{{\mathbf H}}
\def\bP{{\mathbf P}}
\def\bQ{{\mathbf Q}}

\def\1{{\mathbf{1}}}

\def\cP{{\mathcal P}}
\def\cR{{\mathcal R}}

\def\mB{{\mathbb B}}
\def\mE{{\mathbb E}}
\def\mI{{\mathbb I}}
\def\mM{{\mathbb M}}
\def\mN{{\mathbb N}}
\def\mP{{\mathbb P}}
\def\mR{{\mathbb R}}
\def\mS{{\mathbb S}}
\def\mL{{\mathbb L}}

\def\sA{{\mathfrak A}}
\def\sB{{\mathfrak B}}
\def\sC{{\mathfrak C}}
\def\sD{{\mathfrak D}}
\def\sI{{\mathfrak I}}
\def\sL{{\mathfrak L}}
\def\sF{{\mathfrak F}}
\def\sS{{\mathfrak S}}
\def\sK{{\mathfrak K}}

\def\p{\partial}
\def\geq{\geqslant}
\def\leq{\leqslant}
\def\ge{\geqslant}
\def\le{\leqslant}
\def\eps{\varepsilon}

\def\[{{\Big[}}
\def\]{{\Big]}}
\def\<{{\langle}}
\def\>{{\rangle}}
\def\({{\Big(}}
\def\){{\Big)}}

\def\e{{\rm e}}

\def\dif{{\mathord{{\rm d}}}}

\def\bpf{\begin{proof}}
\def\epf{\end{proof}}
\def\div{{\rm div}}
\def\var{{\rm var}}

\endlocaldefs

\begin{document}

\begin{frontmatter}
\title{SDE driven by multiplicative cylindrical $\alpha$-stable noise with distributional drift}
\runtitle{Stable-driven SDE with distributional drift}

\begin{aug}
\author[A]{\fnms{Zimo}~\snm{Hao}\ead[label=e1]{zhao@math.uni-bielefeld.de}}
\and
\author[B]{\fnms{Mingyan}~\snm{Wu}\ead[label=e2]{mingyanwu.math@gmail.com}\orcid{0000-0002-7356-2130}}
\address[A]{Fakult\"at f\"ur Mathematik, Universit\"at Bielefeld\printead[presep={,\ }]{e1}}

\address[B]{
School of Mathematical Sciences, Xiamen University\printead[presep={,\ }]{e2}}
\end{aug}

\begin{abstract}
For $\alpha \in (1,2)$, we study the following stochastic differential equation driven by a non-degenerate symmetric $\alpha$-stable process in $\mathbb{R}^d$:
 \begin{align*}
\dif X_t=b(t,X_t){\mathord{{\rm d}}} t+\sigma(t,X_{t-}){\mathord{{\rm d}}} L_t^{(\alpha)},\ \ X_0 =x \in \mathbb{R}^d,
\end{align*}
where $b$ belongs to $ L^\infty(\mathbb{R}_+;\mathbf{C}^{-\beta}(\mathbb{R}^d))$ with some $\beta\in(0,\alpha-1)$, and $\bC^\beta$ denotes a Besov space (see Definition \ref{iBesov} below). The coefficient $\sigma:\mathbb{R}_+\times \mathbb{R}^d \to \mathbb{R}^d \otimes \mathbb{R}^d$ is a measurable matrix-valued function. The noise $L_t^{(\alpha)}=(L_t^{(\alpha),1},...,L_t^{(\alpha),d})$ consists of independent $1$-dimensional symmetric $\alpha$-stable processes, and is referred to as a cylindrical $\alpha$-stable process. We establish the well-posedness of weak solutions to the SDE, and provide quantitative stability estimates with respect to the drift coefficients. 
\end{abstract}



\end{frontmatter}
\tableofcontents

\section{Introduction}

First of all, consider the following stochastic differential equation (abbreviated as SDE):
 \begin{align}\label{SDE:B}
\dif X_t=b(t,X_t)\dif t+\sigma(t,X_{t})\dif W_t,\ \ X_0 =x \in \mR^d,
\end{align}
where $W$ is a $d$-dimensional standard Brownian motion, $\sigma: \mR_+\times\mR^d\to \mR^d\otimes \mR^d$ is a $d\times d$-matrix valued measurable function, and $b:\mR_+\times\mR^d\to\mR^d$ is a measurable
vector field. Suppose that $\sigma$ satisfies the following uniformly elliptic and Lipschitz condition:

\medskip
\noindent
{\bf (H$^\sigma$)} 
{ There is a constant $c_0>1$ such that for all $t\ge 0$ and $x,\xi\in \mR^d$,
\begin{align}\label{eq:sigma}
c_0^{-1}|\xi|\le|\sigma(t,x)\xi|\leq {c}_0|\xi|,\ \ \hbox{and} \ \  |\nabla_x \sigma(t,x)|\leq {c}_0.
\end{align}}

\noindent
Under the condition {\bf (H$^\sigma$)}, assuming $b\in L^\infty(\mR_+;\bC^{-\beta})$ (here, $\bC^{-\beta}$ is a Besov space, see Definition \ref{iBesov} below), the weak and strong well-posedness are well-known for $\beta<0$ (see \cite{XXZZ20} and the references therein). For $\beta\in[0,\frac{1}{2}]$, the weak well-posedness was established in \cite{ZZ18}. As for $\beta\in(\frac12,\frac23)$ and $\sigma=\mI_{d\times d}$, the weak well-posedness was given in \cite{DD16, CC18} via using rough path and paracontrolled calculus. Recently, when  $\beta\in(\frac12,1]$ and $\sigma=\mI_{d\times d}$, under some additional conditions on $\div b$, the weak well-posedness was investigated in \cite{HZ23, GP23}. 

\medskip

Since it is well known that a $d$-dimensional standard Brownian motion is a special case of an $\alpha$-stable process with $\alpha = 2$, a natural question arises: can the well-posedness results established for Brownian motion be extended to $\alpha$-stable processes with $\alpha < 2$? Precisely, under what conditions on $\beta$ does the weak well-posedness hold for SDEs \eqref{SDE:B} driven by $\alpha$-stable processes with $\alpha<2$? In particular, if we replace the independent components of the $d$-dimensional standard Brownian motion $W_t = (W_t^1, \dots, W_t^d)$ by  independent $1$-dimensional standard $\alpha$-stable processes, we shall get a $d$-dimensional cylindrical $\alpha$-stable process $L_t^{(\alpha)} = (L_t^{(\alpha),1}, \dots, L_t^{(\alpha),d})$. 
It is also worth noting that cylindrical $\alpha$-stable processes arise naturally in the study of $N$-particle systems as well as propagation of chaos (see \cite{Ca22} for example), in which we can see that the joint independence of $\{L^{(\alpha),i}\}_{i=1}^N$ plays a vital role in the following $N$-particle system:
\begin{align*}
\dif X^{N,i}_t=\frac{1}{N}\sum_{i\ne i}K(X^{N,i}_t-X^{N,j}_t)\dif t+\dif L^{(\alpha),i}_t,
\end{align*}
where $K:\mR^d\to\mR^d$ is the interaction kernel. 

\medskip

In this paper, we fix $\alpha\in (1,2)$ and $T>0$, and consider the following SDE, driven by a large class of $\alpha$-stable processes $L^{(\alpha)}$, on the time interval $[0,T]$ with distribution drifts:

\begin{align}\label{in:SDE}
\dif X_t=b(t,X_t)\dif t+\sigma(t,X_{t-})\dif L_t^{(\alpha)},\ \ X_0 =x \in \mR^d,
\end{align}
where $b\in L^\infty_T\bC^{-\beta}:=L^\infty([0,T];\bC^{-\beta})$ with $\beta \in (0,1)$, and $L^{(\alpha)}$ is a $d$-dimensional symmetric $\alpha$-stable process on some probability space $(\Omega,\sF,\mP)$ with the L\'evy measure given by
\begin{align}\label{CC02}
\nu^{(\alpha)}(A)=\int^\infty_0\left(\int_{\mS^{d-1}}\frac{1_A (r\theta)\Sigma(\dif\theta)}{r^{1+\alpha}}\right)\dif r,\quad A\in\sB(\mR^d),
\end{align} 
where $\Sigma$ is a finite measure over the unit sphere $\mS^{d-1}$. 
Throughout this paper, we always assume that the following condition holds:

\medskip\noindent
$\bf (ND)$ The L\'evy measure is {\it non-degenerate}, that is, for each $ \theta_0\in \mS^{d-1}$,  $$
\int_{\mS^{d-1}}|\theta\cdot\theta_0|\Sigma(\dif \theta)>0.
$$

For the cylindrical $\alpha$-stable process, we have
$\Sigma(\dif \theta) = c_\alpha \sum_{i=1}^d \delta_{\theta_i}(\dif \theta)$,  
where $c_\alpha > 0$ is a constant and $\delta_{\theta_i}$ denotes the Dirac measure at the point $\theta_i := (0, \dots, 1 \text{ (in the $i$-th position)}, \dots, 0)$. This measure satisfies the non-degeneracy condition \textbf{(ND)}. Compared to the standard $\alpha$-stable process, whose L\'evy measure is  
$\nu^{(\alpha)}(\dif z) = {c_{d,\alpha} \dif z}/{|z|^{d+\alpha}}$,  
with $\Sigma$ being the uniform (Lebesgue) measure on the unit sphere $\mS^{d-1}$, the L\'evy measure of the cylindrical $\alpha$-stable process is significantly more singular, which will present additional analytical challenges when we study the corresponding equations. We will see this in Subsection \ref{Com23}.

\subsection{Main results}

Since $b$ is a distribution, the drift term is not meaningful in the classical sense, as it is not possible to assign a value to a distribution at the point $X_t$. So to define solutions, a natural approach is to use mollifying approximations. Let $\phi_m(x) := m^d \phi(mx)$, $m \in\mN$, be a family of mollifiers, where $\phi \in C_c^\infty(\mathbb{R}^d)$ is a smooth probability density function with compact support. The smooth approximation of $b$ is then defined by convolution as follows:  
\begin{align}\label{BN1}
b_m (t, x) := b(t, \cdot) * \phi_m (x).
\end{align}
Denote the space of probability measures on $\mR^d$ by $\cP(\mR^d)$. We state the following definition of weak solutions.

\begin{definition}[Weak solutions]\label{Def1}  
Let $ (\Omega, \sF, (\sF_t)_{t \geq 0}, \mathbb{P})$ be a stochastic basis, and let $(X, L)$ be a pair of $\mathbb{R}^d$-valued, c$\rm\grave{a}$dl$\rm\grave{a}$g, $(\sF_t)$-adapted processes on $ (\Omega, \sF, (\sF_t)_{t \geq 0}, \mathbb{P})$. We call $( X, L)$ with $ (\Omega, \sF, (\sF_t)_{t \geq 0}, \mathbb{P})$ a weak solution of the SDE \eqref{in:SDE} with initial distribution $\mu \in \mathcal{P}(\mathbb{R}^d)$ if $L$ is an $(\sF_t)$-$\alpha$-stable process with the L\'evy measure $\nu$ given by \eqref{CC02} which satisfies the condition $\bf (ND)$, and $\mathbb{P} \circ X_0^{-1} = \mu$, and  
$$
X_t = X_0 + A^b_t + \int_0^t \sigma(s,X_{s-})\dif L_t,\quad \text{for all $t \in [0, T]$,} \quad \text{a.s.},
$$  
where $A^b_t := \lim_{n \to \infty} \int_0^t b_n(s, X_s) \, \mathrm{d}s$ exists in the $L^2$-sense.  
\end{definition}
\br
To the best of our knowledge, this notion of weak solutions was first introduced in \cite{BC01} (see also \cite{ZZ17}). It is important to note that the definition of a weak solution in Definition \ref{Def1} depends on the choice of mollifiers $\phi_m$.
\er

Before presenting our results, we first discuss the regime of the order $\beta$. For simplicity, we assume $\sigma = \mI_{d \times d}$ and $b = b(x)$, and begin with a basic scaling analysis of the drift term. Let $X$ be a weak solution to \eqref{in:SDE}. For any $\lambda > 0$, define the rescaled processes and drifts by  $$
X_t^\lambda := \lambda^{-1} X_{\lambda^\alpha t},  \quad 
L_t^\lambda := \lambda^{-1} L^{(\alpha)}_{\lambda^\alpha t},  \quad
b_\lambda(x) := \lambda^{\alpha - 1} b(\lambda x).
$$  
Formally, it follows that  
$$
\mathrm{d} X_t^\lambda = b_\lambda(X_t^\lambda) \, \mathrm{d}t + \mathrm{d}L_t^\lambda.
$$  
Moreover, the $\bC^{-\beta}$ norm of the rescaled drift satisfies  
$$
\|b_\lambda\|_{\bC^{-\beta}} \approx \lambda^{\alpha - 1 - \beta} \|b\|_{\bC^{-\beta}}.
$$ 
As $\lambda \to 0$, we distinguish the following three regimes based on the behavior of this norm:  
$$
\text{\bf Subcritical: } \beta < \alpha - 1; \quad  
\text{\bf Critical: } \beta = \alpha - 1; \quad  
\text{\bf Supercritical: } \beta > \alpha - 1.
$$
This paper focuses on the subcritical case. Notably, there exists a well-known counterexample due to \cite{TTW74}, which shows that in the one-dimensional setting, there exists a drift $b \in \bC^{-\beta}$ with $\beta > \alpha - 1$ (supercritical case) for which both pathwise uniqueness and uniqueness in law fail. Furthermore, to solve SDE \eqref{in:SDE}, a key step is to obtain regularity estimates for the associated Kolmogorov equation:  
\begin{align}\label{in:PDE}
    \partial_t u = \sL^{(\alpha)}_\sigma u + b \cdot \nabla u + f,
\end{align} 
where the nonlocal operator $
\sL^{(\alpha)}_\sigma$ is defined as  
\begin{align}\label{eq:nolocal}
\sL^{(\alpha)}_\sigma f(t,x) := \int_{\mR^d} \left( f(t,x+\sigma(t,x)z) - f(t,x) - \sigma(t,x) z \cdot \nabla f(t,x) \right) \nu^{(\alpha)}(\dif z).
\end{align}
Analogous to the classical Schauder theory for nonlocal equations, when the drift coefficient $b \in \bC^{-\beta}$ with $\beta > 0$, the solution $u$ attains at most $\bC^{\alpha - \beta}$ regularity. Thus, in this case, the product $b \cdot \nabla u$ is  well-defined when $(\alpha - 1 - \beta) -\beta > 0 $, that is, $\beta < \frac{\alpha - 1}{2}$ (see Subsection \ref{sec:3.1} for more details).

\medskip

Now, we summarize several known results concerning the well-posedness of SDE \eqref{in:SDE} in the subcritical regime. 
The results are systematically categorized according to:
\begin{enumerate}
\item The structure L\'evy measure: standard and cylindrical;
    \item The noise structure: additive noise case, and multiplicative noise case;
    \item The singularity strength characterized by the exponent $\beta \in (-\infty,0)\cup (0,\frac{\alpha-1}{2}) \cup [\frac{\alpha-1}{2},\alpha-1)  $.
\end{enumerate}
For clarity, we adopt the following abbreviations:
\begin{center}
{\bf SW}: \emph{Strong Well-posedness of SDEs}\\
 (existence of a strong solution and pathwise uniqueness);\\
{\bf WW}: \emph{Weak Well-posedness of SDEs}\\
(existence of a weak solution and uniqueness in law);\\
{\bf GW}: \emph{Well-posedness of Generalized Martingale Problems}
(see \cite{KP20} for the definition);
\end{center}
\begin{table}[H]\label{table}
\centering
\renewcommand{\arraystretch}{1.5}
\setlength{\tabcolsep}{1.25em}
\caption{Well-posedness results}
\begin{tabular}{c|c|c|c|c}
\toprule
\multicolumn{3}{c|}{$\alpha\in (1,2)$, $b \in \bC^{-\beta}$}   & Standard & Cylindrical  \\ \midrule
\multirow{5}{*}{Additive}  &\multicolumn{2}{c|}{$\beta<0 $} 
&  {\bf SW}:~Pri$^{\mbox{\tiny 12}}_{\mbox{\tiny\cite{Pr12}}}$,~Zh$^{\mbox{\tiny 13}}_{\mbox{\tiny\cite{Zh13b}}}$ 
& {\bf WW\&SW}:~CZZ$^{\mbox{\tiny 21}}_{\mbox{\tiny\cite{CZZ21}}}$ \\ \cline{2-5}
& \multirow{4}{*}{$\beta>0$}
&$\beta \in(0,\frac{\alpha-1}{2})$  
& \thead{ { \bf SW}:~ABM$^{\mbox{\tiny 18}}_{\mbox{\tiny\cite{ABM18}}}$
\\ 
{\bf WW}:~CM$^{\mbox{\tiny 19}}_{\mbox{\tiny\cite{CM19}}}$
}
& {\bf GW}:~KP$^{\mbox{\tiny 20}}_{\mbox{\tiny\cite{KP20}}}$
\\ \cline{3-5}
&
&$\beta \in [\frac{\alpha-1}{2},\frac{2(\alpha-1)}{3})$   
& \multicolumn{2}{c}{{\bf GW}:~KP$^{\mbox{\tiny 20}}_{\mbox{\tiny\cite{KP20}}}$}
\\ \cdashline{3-5}
&
&$\beta \in [\frac{2(\alpha-1)}{3},\alpha-1)$   
& --
& --
\\

\midrule

\multirow{3}{*}{Multiplicative}  &\multicolumn{2}{c|}{$\beta<0 $}  
& \multicolumn{2}{c}{ 
{\bf WW\&SW}:~CZZ$^{\mbox{\tiny 21}}_{\mbox{\tiny\cite{CZZ21}}}$
}
\\ \cline{2-5}
& \multirow{2}{*}{$\beta>0$}
& $\beta\in(0,\frac{\alpha-1}{2})$  
& {\bf WW}: LZ$^{\mbox{\tiny 22}}_{\mbox{\tiny\cite{LZ22}}}$
&   --\\
\cline{3-5}
&
&$\beta \in[\frac{\alpha-1}{2},\alpha-1)$  
& --
&   --\\
\bottomrule
\end{tabular}
\end{table}
\begin{enumerate}[(i)]
\item Let $\beta < 0$. The strong well-posedness was first established for $\beta < \frac{\alpha}{2} - 1$ in \cite{Pr12} (see also \cite{Zh13} for the Sobolev case) in the setting of additive noises. For multiplicative noises, including the cylindrical $\alpha$-stable process, weak and strong well-posedness were obtained in \cite{CZZ21} for both the ranges $\alpha > 1$ and $\alpha \le 1$, under the conditions $\beta < (\alpha - 1) \vee 0$ and $\beta < \frac{\alpha}{2} - 1$ respectively.

\vspace{0.5em}

\item Let $\beta >0$.

\vspace{0.5em}

$\bullet$ Let $\beta \in (0, \frac{\alpha - 1}{2})$. The well-posedness of the weak solution was established in \cite{LZ22} for SDEs driven by a multiplicative standard $\alpha$-stable process. However, the method used in that work does not extend to the cylindrical $\alpha$-stable case (see Section \ref{Com23} for details). In the additive case, we refer to \cite{CM19}, where the authors also claimed the applicability of their method to the multiplicative case in dimension one (see \cite[Remark 15]{CM19}). It is worth noting that in one dimension, the cylindrical and standard $\alpha$-stable processes coincide.

\vspace{0.5em}

$\bullet$ Let $\beta \in \left[\frac{\alpha - 1}{2}, \frac{2(\alpha - 1)}{3}\right)$ and $\sigma = \mI_{d \times d}$. The existence and uniqueness of a generalized martingale solution were obtained in \cite{KP20} using the paracontrolled distributions framework. However, extending their results from additive to multiplicative noise is nontrivial, as paracontrolled calculus relies heavily on Fourier analysis of the operator $\e^{t\Delta^{\alpha/2}}$. We also refer to \cite{KP25} for the rough weak solutions. 
\end{enumerate}

Here is our main result about the well-posedness of the weak solution to SDE \eqref{in:SDE} in the sense of Definition \ref{Def1}. 
 
\bt[Weak well-poseness]\label{thm:G-mart}
Let $T>0$, $\alpha\in(1,2)$ and $\beta\in(0,\alpha-1)$.
Assume that the condition {\bf (H$^\sigma$)} holds with some constant $c_0$, and
\begin{align}\label{con:b}
    (i)~~b\in L^\infty_T\bC^{-\beta},\quad \text{if $\beta\in (0,\tfrac{\alpha-1}{2}$)}; \quad  \quad  (ii)~~b,\div b\in L^\infty_T\bC^{-\beta}, \quad \text{if $\beta\in [\tfrac{\alpha-1}{2}, \alpha-1)$}.
\end{align}
Then for any $\mu\in \cP(\mR^d)$, there is a unique weak solution to SDE \eqref{in:SDE} in the sense of Definition \ref{Def1} such that the following Krylov's estimate holds: for any $ f\in L_T^1 C_b(\mR^d)$ and $0\leq t_0<t_1\leq T$,
\begin{align}\label{sec1:main2:es1}
\Big\|\int_{t_0}^{t_1}f(s,X_s)\dif s\Big\|_{L^p(\Omega)}\lesssim_c (t_1-t_0)^{\frac{\alpha-\beta}{\alpha}}\|f\|_{L_T\bC^{-\beta}},
\end{align}
where the constant $c$ is independent of $t_0$, $t_1$ and $f$.
\et

\br[Our contributions]
Our result fills all the gaps in Table \ref{table}. More precisely, we obtain the weak well-posedness for the multiplicative cylindrical case with $\beta\in (0,\frac{\alpha-1}{2})$. Moreover, if $\div b\in L^\infty_T\bC^{-\beta}$ additionally, we establish the unique weak solution for the additive case with $\beta\in(\frac{2(\alpha-1)}{3},\alpha-1)$, and the multiplicative case with $\beta\in [\frac{\alpha-1}{2},\alpha-1)$.
\er

For any $\sigma$ satisfying {\bf (H$^\sigma$)} and $b_1,b_2\in L^\infty_T\bC^{-\beta}$ satisfying \eqref{con:b} with some $\beta\in(0,{\alpha-1})$, based on Theorem \ref{thm:G-mart}, there are unique weak solutions $X^1$ and $X^2$ to SDE \eqref{in:SDE} with drift $b=b_1$ and $b=b_2$ respectively. We denote the time marginal law of $X^i$ by $\bP_i(t)\in\cP(\mR^d)$, $i=1,2$. 
Our second main result is the following stability result. For convenience of the notation,  we introduce the following parameter set:
\begin{align*}
\Theta:=(T,d,\alpha,\beta,c_0).
\end{align*}
\bt[Stability estimates]\label{thmSt}
Let $T>0$ and $\alpha\in(1,2)$. Assume that $\bP_1(0)=\bP_2(0)$, $\sigma$ satisfies {\bf (H$^\sigma$)} with constant $c_0$, and $b_1,b_2\in L^\infty_T\bC^{-\beta}$ satisfy \eqref{con:b} with some $\beta\in[0,{\alpha-1})$. Then,
\begin{itemize}
    \item[(i)] when $\beta<\frac{\alpha-1}{2}$, for any  $\theta\in[\beta,\alpha-1-\beta)$ and $\eps>0$, there is a constant $c=(\Theta,\theta,\eps,\|b_1\|_{L^\infty_T\bC^{-\beta}})>0$ such that for any $t\in [0,T]$,
\begin{align*}
\|\bP_1(t)-\bP_2(t)\|_{\rm var} \lesssim_c t^{\frac{\alpha-1-2\theta-\eps}{\alpha} }\|b_1-b_2\|_{L^\infty_T\bC^{-\theta}};
\end{align*}
\item[(ii)] when $\beta\ge\frac{\alpha-1}{2}$, for any  $\theta\in[\beta,\alpha-1)$ and $\eps>0$, there is a constant $c>0$ depending on $\Theta,\theta,\eps,\|b_1\|_{L^\infty_T\bC^{-\beta}},\|\div b_1\|_{L^\infty_T\bC^{-\beta}}$, such that for any $t\in [0,T]$,
\begin{align*}
\|\bP_1(t)-\bP_2(t)\|_{\rm var} \lesssim_c t^{\frac{\alpha-1-\theta-\eps}{\alpha} }\left(\|b_1-b_2\|_{L^\infty_T\bC^{-\theta}}+\|\div b_1-\div b_2\|_{L^\infty_T\bC^{-\theta}}\right).
\end{align*}

\end{itemize}

\et 
\br
The stability result clearly indicates that the weak solution obtained in Theorem \ref{thm:G-mart} is independent of the specific choice of mollifier functions $\phi_m$.
\er

\subsection{Comparison with \cite{LZ22} for $\beta\in(0,\frac{\alpha-1}{2})$}\label{Com23}

The primary difference between our work and \cite{LZ22} lies in the type of stochastic processes considered. In \cite{LZ22}, the process studied is the standard (i.e., rotationally invariant) $\alpha$-stable process $L^{\mathrm{sta}}_t = (L^{\mathrm{sta},1}_t, \dots, L^{\mathrm{sta},d}_t)$, where the components $L^{\mathrm{sta},i}_t$ are not independent. The L\'evy measure of $L^{\mathrm{sta}}_t$ is given by $\dif z / |z|^{d+\alpha}$,  which determines the infinitesimal generator of the solution $X$ to SDE \eqref{in:SDE} (with $L^{(\alpha)}$ substituted by $L^{\mathrm{sta}}$) as follows:

\begin{align}\label{in:kappa}
\begin{split}
    \sL^{\mathrm{sta}}_{\sigma} f(t,x) &= \int_{\mathbb{R}^d} \left( f(\sigma(t,x)+z) - f(x) - \sigma(t,x)z \cdot \nabla f(x) \right) \frac{\dif z}{|z|^{d+\alpha}} \\
    &= \int_{\mathbb{R}^d} \left( f(x+z) - f(x) - z \cdot \nabla f(x) \right) \frac{\dif z}{|\sigma^{-1}(t,x)z|^{d+\alpha} |\det \sigma(t,x)|}.
\end{split}
\end{align}
Therefore, the operator considered in \cite{LZ22} takes the form
\begin{align}\label{eq:XM98}
\sA^{(\alpha)}_{\kappa} f(t,x) := \int_{\mathbb{R}^d} \left(f(x+z) - f(x) - \nabla f(x) \cdot z\right) \kappa(t,x,z)\frac{ \dif z}{|z|^{d+\alpha}},
\end{align}
where $\kappa : \mathbb{R}_+ \times \mathbb{R}^{2d} \to \mathbb{R}_+$ is a measurable function. In particular, the choice
$$
\kappa(t,x,z) := \frac{|z|^{d+\alpha}}{\left( |\sigma^{-1}(t,x)z|^{d+\alpha} |\det \sigma(t,x)| \right)},
$$
yeilds the operator identity $\sL^{\mathrm{sta}}_{\sigma} = \sA^{(\alpha)}_{\kappa}$ under the standard case. However, when the L\'evy measure $\nu^{(\alpha)}$ fails to be absolutely continuous 
with respect to the Lebesgue measure (as occurs in the cylindrical case), 
such a representation becomes unattainable. That is, there exists no measurable function $\kappa$ for which the operator equality $\sL^{(\alpha)}_{\sigma} = \sA^{(\alpha)}_{\kappa}$ holds.

\paragraph{Challenges in Analyzing $\sL^{(\alpha)}_{\sigma}$.}

Notice that, when dealing with the PDE \eqref{in:PDE}, direct treatment of $\sL^{(\alpha)}_{\sigma}$ presents significantly greater challenges compared to working with $\sA^{(\alpha)}_{\kappa}$. For instance, when $\kappa = \kappa(t,x)$, the operator $\sA^{(\alpha)}_{\kappa} $ reduces to
$$
\sA^{(\alpha)}_{\kappa} f(t,x) = \kappa(t,x) \Delta^{\frac{\alpha}{2}} f(x),
$$
which allows us to apply classical functional inequalities (e.g., \eqref{eq:pro}) to obtain the estimate:
$$
\|\sA^{(\alpha)}_{\kappa}f(t)\|_{\bC^{-\beta}} = \|\kappa(t)\Delta^{\frac{\alpha}{2}}f\|_{\bC^{-\beta}} \lesssim \|\kappa(t)\|_{\bC^{\beta+\varepsilon}} \|\Delta^{\frac{\alpha}{2}}f\|_{\bC^{-\beta}} \lesssim \|\kappa(t)\|_{\bC^{\beta+\varepsilon}} \|f\|_{\bC^{\alpha - \beta}},
$$
for any $\varepsilon > 0$. In contrast, the operator $\sL^{(\alpha)}_{\sigma}$ exhibits a fundamentally more complex structure, which renders the derivation of norm estimates particularly challenges. A representative difficulty arises when attempting to establish estimates of the form:
$$
\| f(\cdot + \sigma(\cdot)z) - f(\cdot) - \sigma(\cdot)z \cdot \nabla f(\cdot) \|_{\bC^{-\beta}} \lesssim c_\sigma |z|^{\alpha + \varepsilon} \|f\|_{\bC^{\alpha - \beta}},
$$
especially in the regime where $-\beta < 0$ (i.e., for negative regularity exponents). To overcome these technical obstacles, in this paper, we employ integral-type heat kernels estimates developed in \cite{CHZ20} to obtain norm bounds of the form (see Lemma \ref{Lem31} below):
$$
\| \sL^{(\alpha)}_{\sigma} f(t) \|_{\bC^{-\beta}} \lesssim \| \nabla \sigma(t) \|_{\infty} \| f \|_{\bC^{\alpha - \beta}}.
$$
This approach necessitates the Lipschitz continuity of the map $x \mapsto \sigma(t,x)$, which justifies our assumption of the Lipschitz property for $\sigma$ throughout this work.

\subsection{Structure and notations}

\noindent
{\bf Outline of paper.} The rest of this paper is organized as follows. In Section \ref{Pre}, we introduce some basic concepts and estimates about Besov spaces and $\alpha$-stable processes. In Section \ref{sec:PDE}, we establish Schauder's estimate and obtain the well-posedness for the non-local parabolic equation with $\alpha$-stable measures being singular and drifts being distributional. In Section \ref{sec:SDE}, we show the first main result of this paper, Theorem \ref{thm:G-mart}. 
Finally, we prove the stability estimate Theorem  \ref{thmSt} in Section \ref{sec:SDE-2}. 

\medskip
\noindent
{\bf Notation and Conventions.} Throughout this paper, we use the following conventions and notations: As usual, we use $:=$ as a way of definition. Define $\mN_0:= \mN \cup \{0\}$ and $\mR_+:=[0,\infty)$. The letter $c=c(\cdots)$ denotes an unimportant constant, whose value may change in different places. We use $A \asymp B$ and $A\lesssim B$ to denote $c^{-1} B \leq A \leq c B$ and $A \leq cB$, respectively, for some unimportant constant $c \geq 1$. We also use $A   \lesssim_c  B$ to denote $A \leq c B$ when we want to emphasize the constant. 
\begin{itemize}
\item  Let $\mM^d$ be the space of all real $d\times d$-matrices, and $\mM^d_{non}$ the set of all non-singular matrices. Denote the identity $d\times d$-matrix by $\mI_{d \times d}$.
 
\item  For every $p\in [1,\infty)$, we denote by $L^p$ the space of all $p$-order integrable functions on $\mR^d$ with the norm denoted by $\|\cdot\|_p$. 

\item  For a Banach space $\mB$ and $T>0$, $q\in[1,\infty]$, we denote by
$$
L_T^q\mB:= L^q([0,T];\mB),\ \  \mL^q_T:=L^q([0,T]\times \mR^d).
$$
\end{itemize}

\section{Preliminary}\label{Pre}

\subsection{Besov spaces}\label{Sec:Be}

In this subsection, we introduce some basic concepts and properties of Besov spaces. Let $\sS(\mR^d)$ be the Schwartz space of all rapidly decreasing functions on $\mR^d$, and $\sS'(\mR^d)$ 
the dual space of $\sS(\mR^d)$ called Schwartz generalized function (or tempered distribution) space. Given $f\in\sS(\mR^d)$, 
the Fourier transform $\hat f$ and the inverse Fourier transform  $\check f$ are defined by
$$
\hat f(\xi) :=(2 \pi)^{-d/2}\int_{\mR^d} \e^{-i\xi\cdot x}f(x)\dif x, \quad\xi\in\mR^d,
$$
$$
\check f(x) :=(2 \pi)^{-d/2}\int_{\mR^d} \e^{i\xi\cdot x}f(\xi)\dif\xi, \quad x\in\mR^d.
$$
For every $f\in\sS'(\mR^d)$, the Fourier and the inverse  transforms are defined by
\begin{align*}
\<\hat{f},\varphi\>:=\<f,\hat{\varphi}\>,\qquad \<\check{f},\varphi\>:=\<f,\check{\varphi}\>, \ \  \forall\varphi\in\sS(\mR^d).
\end{align*}
Let $\chi:\mR^{d}\to[0,1]$ be a radial smooth function with
\begin{align*}
\chi(\xi)=
\begin{cases}
1, & \ \  |\xi|\leq 1,\\
0, &\ \ |\xi|>3/2.
\end{cases}
\end{align*}
For $\xi \in \mR^d$, define $\psi(\xi):=\chi(\xi)-\chi(2\xi)$ and for $j\in\mN_0$,
\begin{align*}
\psi_j(\xi){:=}\psi(2^{-j}\xi).
\end{align*}
Let $  B_r := \{\xi\in \mR^d :  |\xi|\leq r\}$ for $r>0$. It is easy to see that $\psi\geq 0$,  supp$\psi\subset B_{3/2}/B_{1/2}$, and
\begin{align}\label{eq:SA00}
\chi(2\xi)+\sum_{j=0}^{k}\psi_j(\xi)=\chi(2^{-k}\xi)\to 1,\ \ \hbox{as}\ \ k\to\infty.
\end{align}
Since $\check \psi_j (y) = 2^{jd}\check \psi (2^j y), j \geq 0$, we have
\begin{align*}
\int_{\mR^d}|x|^\theta |\nabla^k\check \psi_j|(x) \dif x\leq {c} 2^{(k-\theta)j} , \ \ \theta>0,\ \ k\in\mN_0,
\end{align*}
where the constant $c$ is equal to $\int_{\mR^d}|x|^\theta |\nabla^k\check \psi|(x) \dif x$ and $\nabla^k$ stands for the $k$-order gradient. The block operators $\cR_j, j \geq 0$ are defined on $\sS'(\mR^d)$ by
\begin{align}\label{eq:Block}
\cR_j f (x):=(\psi_j\hat{f})^{\check\,}(x)=\check\psi_j* f (x) = 2^{jd} \int_{\mR^d} \check\psi(2^{j}y)f(x-y) \dif y,
\end{align}
and 
$
\cR_{-1}  f(x):= (\chi(2\cdot)\hat{f})^{\check\,}(x)=(\chi(2\cdot))\check{ }* f(x).
$

\br\label{rem:AA01}
For $j\ge-1$,  by definitions, it is easy to see that
\begin{align*}
\cR_j=\cR_j\widetilde\cR_j,\quad \text{where $\widetilde\cR_j{:=}\sum_{\ell=-1}^1\cR_{j+\ell}$ with $\cR_{-2}:=0$},
\end{align*}
and $\cR_j$ is symmetric in the sense of 
\begin{align*}
\int_{\mR^d}\cR_jf(x)g(x)\dif x=\int_{\mR^d}f(x)\cR_jg(x)\dif x, \ \ f \in \sS'(\mR^d),\, g \in \sS(\mR^d).
\end{align*}
\er

Now we state the definitions of Besov spaces.
\bd[Besov spaces]\label{iBesov}
For every $s\in\mR$ and $p,q\in[1,\infty]$, the Besov space $\bB_{p,q}^s(\mR^d)$ is defined by
$$
\bB_{p,q}^s(\mR^d):=\Big\{f\in\sS'(\mR^d)\, \big| \, \|f\|_{\bB^s_{p,q}}:= \[ \sum_{j \geq -1}\left( 2^{s j} \|\cR_j f\|_{p} \right)^q \]^{1/q}  <\infty\Big\}.
$$
In particular, we set
$$
\bC^s(\mR^d):=\bB^s_{\infty,\infty}.
$$
Notice that $\bB_{p,q}^s(\mR^d)$ is a Banach space for all $s,p,q$.
\ed 
\br
It is worth discussing here the equivalence between the Besov and H\"older  spaces, which will be used in various contexts in this paper without much explanation. For $s>0$, let $\sC^s(\mR^d)$ be the classical $s$-order H\"older space consisting of all measurable functions $f:\mR^d\to\mR$ with
\begin{align*}
\|f\|_{\sC^s}:=\sum_{j=0}^{[s]}\|\nabla^jf\|_\infty+[\nabla^{[s]}f]_{\sC^{s-[s]}}<\infty,
\end{align*}
where $[s]$ denotes the largest integer less than or equal to $s$, and
\begin{align*}
\|f\|_\infty:=\sup_{x\in\mR^d}|f(x)|,\quad [f]_{\sC^\gamma}:=\sup_{h\in\mR^d}\frac{\|f(\cdot+h)-f(\cdot)\|_\infty}{|h|^\gamma},~\gamma\in(0,1).
\end{align*}
If $s>0$ and $s\notin\mN$, we have the following equivalence between $\bC^s (\mR^d)$ and $ \sC^s  (\mR^d)$: (cf. \cite{Tr92})
\begin{align*}
\|f\|_{\bC^s}\asymp\|f\|_{\sC^s}.
\end{align*}
However, for any $n\in\mN_0$, we only have one side control that is
$
\|f\|_{\bC^n}\lesssim\|f\|_{\sC^n}.
$ 
\er

Let $S_k$ be the cut-off low frequency operator defined by
\begin{align*}
S_kf:=\sum_{j=-1}^{k-1}\cR_j f\to f,\ \ k\to\infty.
\end{align*}
For $f,g\in\sS'(\mR^d)$, define
$$
f\prec g:=\sum_{  j\geq 1} \left( S_{j-1}f\right) \cR_j g,\ \ f\circ g:=\sum_{|i-j|\leq1} (\cR_i f)\cR_jg.
$$
Then the Bony decomposition of $fg$ is formally given by (cf. \cite{BCD11}, Definition 2.81)
\begin{align}\label{Bony}
fg=\left( f\prec g+ f\circ g\right) +g\prec f=:f\preceq g+f\succ g.
\end{align}

Recall the following Bernstein's inequality (cf. \cite{BCD11}, Lemma 2.1).

\bl[Bernstein's inequality]\label{Bern}
For every $k\in\mN_0$, there is a constant $c=c(d,k)>0$ such that for all $j\ge-1$ and $1\leq p_1\leq p_2 \leq \infty$,
\begin{align*} 
\|\nabla^k\cR_j f\|_{p_2}  \lesssim_c  2^{(k+ d (\frac{1}{p_1}-\frac{1}{p_2}))j}\|\cR_j f\|_{p_1}.
\end{align*}
In particular, for any $s \in \mR$ and $1\leq p, q \leq \infty$,
\begin{align}\label{S2:Bern}
\|\nabla^k f\|_{\bB^{s}_{p,q}} \lesssim_c  \|f\|_{\bB^{s+k}_{p,q}}.
\end{align}
\el

\br[Mollification in Besov spaces]
Let $\phi_n:=n^{-d}\phi(nx)$, $(n > 0)$, be the mollifier for fixed $\phi \in C^{\infty}_c(\mR^d)$ being a smooth function with compact support and unit integral. unit integral. Let $\beta\in\mR$ with $  \eps \in [0,1]$. It is easy to check that there is a constant $c>0$ such that for all $f\in\bC^{\beta+\eps}$ and $n\in\mN$,
\begin{align}\label{0725:new00}
    \|f-f_n\|_{\bC^{\beta}}\le c n^{-\eps}\|f\|_{\bC^{\beta+\eps}}.
\end{align}
\er

We introduce the following paraproduct estimates, which can be found in \cite{BCD11}, Theorems 2.82 \& 2.85, or \cite{HZ23}, Lemma 2.5. We also refer to \cite{GIP15} and \cite{Bo81} for readers who want to know more about Bony’s paraproduct.

\bl[Paraproduct estimates]\label{productlaw}
Let $p,p_1,p_2,q, q_1,q_2\in[1,\infty]$ with  $\frac1p=\frac1{p_1}+\frac1{p_2}$ and $\frac1q=\frac1{q_1}+\frac1{q_2}$ and $\alpha,\beta\in\mR$.
\begin{enumerate} 
\item If $\beta<0$, then there is a constant $c=c(\alpha,\beta,d,p,q,p_1,q_1,p_2,q_2)>0$ such that
$$
\|f\prec g\|_{\bB^{\alpha+\beta}_{p,q}}\lesssim_ c \|f\|_{\bB^{\beta}_{p_1,q_1}} \|g\|_{\bB^{\alpha}_{p_2,q_2}}.
$$
Moreover, for $\beta=0$, we have
$$
\|f\prec g\|_{\bB^{\alpha}_{p,q}}\lesssim_c
\|f\|_{p_1} \|g\|_{\bB^{\alpha}_{p_2,q}}.
$$

\item If $\alpha+\beta>0$,  then there is a constant $c=c(\alpha,\beta,d,p,q,p_1,q_1,p_2,q_2)>0$ such that
$$
\|f\circ g\|_{\bB^{\alpha+\beta}_{p,q}}\lesssim_c \|f\|_{\bB^{\beta}_{p_1,q_1}} \|g\|_{\bB^{\alpha}_{p_2,q_2}}.
$$
Moreover, when $\alpha+\beta=0$ and $q=1$, we have
\begin{align*} 
\|f\circ g\|_{\bB^{0}_{p,\infty}}\lesssim_ c
\|f\circ g\|_p\lesssim_c \|f\|_{\bB^{-\alpha}_{p_1,q_1}} \|g\|_{\bB^{\alpha}_{p_1,q_2}}.
\end{align*}
\end{enumerate}
\el

\bc
For any $s>0$ and $\eps>0$, there is a constant $c=c(s,\eps)>0$ such that
\begin{align}\label{eq:pro}
\| f g\|_{\bC^{-s}}  \lesssim_c \|f\|_{\bC^{s+\eps}}\|g\|_{\bC^{-s}}.
\end{align}
\ec

We also need the following interpolation inequality (cf. \cite{BCD11}, Theorem 2.80).

\bl[Interpolation inequality]
Let $s_1,s_2\in \mR$ with $s_2>s_1$. For any $p\in [1,\infty]$ and $\theta \in (0,1)$, there is a constant $c=c(s_1,s_2,p)>0$ such that 
\begin{align}\label{InIn:00}
\|f\|_{\bB_{p,1}^{\theta s_1+ (1-\theta) s_2} } \lesssim_c \|f\|_{\bB^{s_1}_{p,\infty}}^\theta\|f\|_{\bB^{s_2}_{p,\infty}}^{1-\theta}.
\end{align}
Furthermore, for any $ s_2> 0 > s_1$,
\begin{align}\label{InIn}
\|f\|_{\infty} \lesssim_c \|f\|_{\bC^{s_1}}^\theta\|f\|_{\bC^{s_2}}^{1-\theta},
\end{align}
where $\theta=s_2/(s_2-s_1)$.
\el
\br
We note that the interpolation and Young's inequality ensure that for any $s_i \in \mR$, $i=0,1,2$, with $ s_0<s_1<s_2$ and $\varepsilon>0$, there is a constant $c_\varepsilon>0$ such that
\begin{align}\label{eq:inter}
\|f\|_{\bC^{s_1}}\le \varepsilon \|f\|_{\bC^{s_2}}+c_\varepsilon \|f\|_{\bC^{s_0}},
\end{align}
\er

\subsection{$\alpha$-stable processes}

Fix $\alpha\in(0,2)$. Let $L^{(\alpha)}_t$ be a $d$-dimensional $\alpha$-stable process with L\'evy measure (or $\alpha$-stable measure) $\nu^{(\alpha)}$ defined as \eqref{CC02}. We say an $\alpha$-stable measure $\nu^{(\alpha)}$ is non-degenerate, if the assumption {\bf (ND)} holds. Note that $\alpha$-stable process $L_t^{(\alpha)}$ has the scaling property,
\begin{align}\label{eq:scal}
 (L_t^{(\alpha)})_{t \geq 0} \overset{(d)}{=}  (\lambda^{-1/\alpha} L^{(\alpha)}_{\lambda t})_{t \geq 0},\ \ \forall\, \lambda >0,
\end{align}
and for any $\gamma_2>\alpha>\gamma_1\geq 0$, 
\begin{align}\label{eq:AE02}
\int_{|z|\leq 1} |z|^{\gamma_2} \nu^{(\alpha)}(\dif z) + \int_{|z|>1} |z|^{\gamma_1}\nu^{(\alpha)}(\dif z)<\infty.
\end{align}
Moreover, it is easy to see that for any $\lambda>0$ and $p\ge2$,
\begin{align}\label{20230418}
\int_{\mR^d}(1\wedge |\lambda z|^p)\nu^{(\alpha)}(\dif z)=\lambda^{\alpha}\int_{\mR^d}(1\wedge|z|^p)\nu^{(\alpha)}(\dif z).
\end{align}
Let $N(\dif r, \dif z)$ be the associated Poisson random measure defined by
$$
N((0,t]\times A)  :=  \sum_{s\in(0,t]} \1_{A}(L_s^{(\alpha)} -  L^{(\alpha)}_{s-}),\ \ A\in \sB(\mR^d\setminus \{0\}) , t >0.
$$
By L\'evy-It\^o's decomposition (cf. \cite{Sa99}, Theorem 19.2), one sees that
\begin{align*}
L^{(\alpha)}_t =  \lim_{\eps \downarrow 0} \int_0^{t} \int_{\eps<|z|\leq 1} z \widetilde N(\dif r,\dif z) + \int_0^{t} \int_{|z|> 1} z   N(\dif r,\dif z),
\end{align*}
where $ \widetilde{N}(\dif r,\dif z):= N(\dif r,\dif z)-\nu^{(\alpha)}(\dif z)\dif r $ is the compensated Poisson random measure. 

In the sequel, we always assume that $\nu^{(\alpha)}$ is symmetric. Hence, we can write
\begin{align}\label{eq:Levy}
L^{(\alpha)}_t = \int_{0}^t \int_{|z|\leq c}z \widetilde N(\dif r , \dif z) + \int_{0}^t \int_{|z|>c} z N(\dif r , \dif z), \ \ \forall c>0.
\end{align}
The following moment estimate is taken from  \cite{CZ18b}, Lemma 2.4,  with some slight modification. 

\bl\label{lem:BH}
Let $T,\delta>0$. Assume that $0\leq \tau_1< \tau_2\leq \tau_1+\delta \leq T$  are two bounded stopping times and $p\in(0,\alpha)$. Let  $g: \mR_+  \times \Omega \to \mM_{non}^d$ be a bounded predictable process, where $\mM^d_{non}$ is the set of all non-singular $d \times d$ matrices. Then, there is a constant $c=c(d,\alpha, p,T)>0$ such that
\begin{align*} 
\mE \left| \int_{\tau_1}^{\tau_2} g(r) \dif L^{(\alpha)}_r \right|^p  \lesssim_c  \delta^{p/\alpha}\left(\left[\mE \|g\|_{L^2([0,T])}^2\right]^{p/2}+ \mE \| g\|_{L^p([0,T])}^p\right).
\end{align*}
\el

\begin{proof}
Noticing that Poisson measures are counting measures, by \eqref{eq:Levy}, we have
\begin{align*} 
\int_{\tau_1}^{\tau_2} g(r) \dif L^{(\alpha)}_r 
= & \int_{0}^{T} \int_{|z|\leq \delta^{1/\alpha}} \widetilde g(r,z)  \widetilde N(\dif r, \dif z)  +  \int_{0}^{T} \int_{|z|>  \delta^{1/\alpha}} \widetilde g(r,z)  N(\dif r, \dif z),
\end{align*}
where $\widetilde g (r,z):=   g(r)z\1_{(\tau_1, \tau_2]}(r) $ is left continuous. On the one hand, by Jensen's inequality and the isometry of stochastic integral,
\begin{align*}
\mE \left|  \int_{0}^{T} \int_{|z|\leq \delta^{1/\alpha}} \widetilde g (r,z)  \widetilde N(\dif r, \dif z) \right|^p 
& \leq  \left[ \mE \left|  \int_{0}^{T} \int_{|z|\leq \delta^{1/\alpha}} \widetilde g (r,z)  \widetilde N(\dif r, \dif z) \right|^2 \right]^{p/2}\\
& =  \left[  \mE \int_{0}^{T} \int_{|z|\leq \delta^{1/\alpha}} |\widetilde g (r,z) |^2 \nu^{(\alpha)}(\dif z) \dif r \right]^{p/2} \\
& \leq  \left[\mE \|g\|_{L^2([0,T])}^2\right]^{p/2}  \left[ \int_{|z|\leq \delta^{1/\alpha}} |z|^2 \nu^{(\alpha)}(\dif z)  \right]^{p/2} \\
& \lesssim \delta^{p/\alpha} \left[\mE \|g\|_{L^2([0,T])}^2\right]^{p/2},
\end{align*}
where we used \eqref{CC02} in the last inequality. On the other hand, using Burkholder’s inequality (cf. \cite[Lemma 2.3]{SZ15}) and \eqref{CC02} for $p\in (1,\alpha)$, we have
\begin{align*}
 \mE \left|  \int_{0}^{T} \int_{|z|>\delta^{1/\alpha}} \widetilde g (r,z)  N(\dif r, \dif z) \right|^p
& \lesssim \mE \left( \int_{0}^{T} \int_{|z|> \delta^{1/\alpha}} |\widetilde g (r,z) | \nu^{(\alpha)}(\dif z) \dif r \right)^{p} \\
& \quad + \mE \int_{0}^{T} \int_{|z|> \delta^{1/\alpha}} |\widetilde g (r,z) |^p \nu^{(\alpha)}(\dif z) \dif r\\
& \lesssim    \mE \|g\|_{L^1([0,T])}^p   \left ( \int_{|z|>\delta^{1/\alpha}} |z| \nu^{(\alpha)}(\dif z) \right)^{p}\\
&\quad+ \mE  \|g\|_{L^p([0,T])}^p \int_{|z|> \delta^{1/\alpha}} |z|^p \nu^{(\alpha)}(\dif z)\\
& \lesssim \delta^{p/\alpha}  \mE \|g\|_{L^{p}([0,T])}^p .
\end{align*}
Observing that $|\sum_{i =1}^n a_i|^p \leq (n^{p-1} \vee 1) \sum_{i =1}^n |a_i|^p =\sum_{i =1}^n |a_i|^p  $ when $p\in (0,1]$, we get
\begin{align*}
\mE \left|  \int_{0}^{T} \int_{|z|>\delta^{1/\alpha}} \widetilde g (r,z)  N(\dif r, \dif z) \right|^p
& \leq  \mE    \int_{0}^{T} \int_{|z|>\delta^{1/\alpha}}  | \widetilde g (r,z)|^p  N(\dif r, \dif z)\\
& =   \mE  \int_{0}^{T} \int_{|z|>\delta^{1/\alpha}}  | \widetilde g (r,z)|^p \nu^{(\alpha)} (\dif z) \dif r \\
& \lesssim \delta^{p/\alpha} \mE \| g\|_{L^p([0,T])}^p,
\end{align*}
where the equality is yielded by the property of martingales (cf. \cite{IW89}). The above calculations derive the desired estimates. 
\end{proof}

We introduce the notation of $q$-variation for stochastic processes.

\bd[$q$-variation]\label{def:SA00} 
Fix $T>0$. Let $0\leq s<t \leq T$, and $q>2$. We say that a process $X:=(X_t)_{t\in[0,T]}$ has a.s. $q$-variation on $[s,t]$ if
\begin{align} 
[X]_{q,\var;[s,t]}:=\left ( \sup_{\cP}\sum_{[t_1,t_2]\in \cP}|X_{t_1}-X_{t_2}|^q \right)^{1/q}<\infty,\quad a.s.,
\end{align}
where the supremum is taken over all partitions $\cP$ of $[s,t]$. In particular, if $[s,t]=[0,T]$, we simply say $X$ has a.s. $q$-variation, and denote by $[X]_{q,\var}:=[X]_{q,\var;[0,T]}$.
\ed

\bl\label{lem:B2}
Let $(L^{(\alpha)}_t)_{t\in[0,T]}$ be an $\alpha$-stable process whose L\'evy measure $\nu^{(\alpha)}$ is given as \eqref{CC02}. Then for any $q>2$, and predictable bounded process $\Sigma =(\Sigma_s)_{s \in [0,T]} $,
\begin{align*}
    \left[\int_0^\cdot \Sigma_{s} \dif L^{(\alpha)}_s\right]_{q,\var}<\infty,\quad a.s.
\end{align*}
\el
\begin{proof}
    By \eqref{eq:Levy}, one sees that 
\begin{align*}
    \int_0^t \Sigma_{s} \dif L^{(\alpha)}_s&=\int_0^t \int_{|z|\le1}\Sigma_{s}\cdot z \,\widetilde{N}(\dif z, \dif s)+\int_0^t \int_{|z|>1}\Sigma_{s}\cdot z  \, N(\dif z, \dif s) 
    =: I_1(t)+I_2(t).
\end{align*}
Observe that $I_2(t)$ is well defined a.s. as the Lebesgue-Stieltjes integral, and equals the absolutely convergent sum $\sum_{s\leq t} \Sigma_s \cdot (L_s^{(\alpha)} - L_{s-}^{(\alpha)})1_{|\Delta L_s|>1}$, which is actually a finite sum since each path of $L^{(\alpha)}$ has at most finitely many jumps of size greater than $1$ before time $t$. Thus, one sees that $t\to I_2(t)$ is a bounded variation function. So we only need to show $[I_1]_{q,\var}<\infty $. Noting that $t\to I_1(t)$ is a c$\rm\grave{a}$dl$\rm\grave{a}$g martingale (cf. \cite{IW89}), by L\'epingle's inequality (cf. \cite{Le76}), we have for any $q
>2$,
\begin{align*}
     \mE [I_1]_{q,\var}^2&\lesssim \mE \left[\sup_{0\le t\le T}\left|I_1(t)\right|^2\right]=\mE \sup_{0\le t\le T}\left|\int_0^t \int_{|z|\le1}\Sigma_{r}\cdot z\widetilde{N}(\dif r,\dif z)\right|^2\\
     &\lesssim \mE  \int_0^T \int_{|z|\le1}|\Sigma_{s}\cdot z|^2\nu^{(\alpha)}(\dif z) \dif s<\infty,
    \end{align*}
    where we used the by Burkholder-Davis-Gundy's inequality (see \cite{SZ15}, Lemma 2.3, or \cite{No75}, Theorem 1) in the last inequality. This completes the proof.
\end{proof}

\subsection{Time-dependent L\'evy-type operators}\label{sec:op}

Let $\alpha \in (1.2)$ and $L^{(\alpha)}$ be an $\alpha$-stable process having symmetric non-degenerate L\'evy measure $\nu^{(\alpha)}$. In this subsection, we start with the following time-inhomogeneous L\'evy process: for $0\leq t< \infty$,
\begin{align}\label{eq:L}
L^\sigma_t :=\int_0^t \sigma_r \dif L_r^{(\alpha)} = \int_0^{t} \int_{\mR^d} \sigma_r z \widetilde N(\dif r,\dif z),
\end{align}
where $\sigma : \mR_+ \to  \mM_{non}^{d}$ is  a bounded measurable function. Define
\begin{align}\label{eq:XM100}
P^\sigma_{s,t}f(x) := \mE \left( x+\int_s^t \sigma_r \dif L_r^{(\alpha)} \right)
\end{align}
for all $f \in C_b^2(\mR^d)$. By It\^o's formula (cf. \cite{IW89}, Theorem 5.1 of Chapter II), one sees that
\begin{align*}
\p_t P^\sigma_{s,t}f(x) = \sL^{(\alpha)}_{\sigma_t}  P^\sigma_{s,t}f(x),
\end{align*}
where 
\begin{align}\label{eq:a} 
\sL^{(\alpha)}_{\sigma_t}  f(x):=\int_{\mR^d}\(f(x+\sigma_tz)-f(x)- \sigma_tz \cdot\nabla f(x)\)\nu^{(\alpha)}(\dif z).
\end{align} 

Below, we always make the following assumption in this subsection:

\medskip
\noindent
$\bf (\bH0)$ There is a constant $a_0 >1$ such that
$$
a_0^{-1}|\xi|\le|\sigma_t\xi|\le a_0|\xi|,\ \ \forall(t,\xi)\in\mR_+\times \mR^d.
$$

\medskip
\noindent
Under the assumptions $\bf (\bH0)$ and $\bf (ND)$, owing to  L\'evy-Khintchine's formula  (cf. \cite{Sa99}, Theorem 8.1) and \eqref{CC02}, for all $ |\xi|\geq1$, we have 
\begin{align*}
|\mE \e^{i\xi\cdot L^{\sigma}_{t} }|\leq&\exp\left(t\int_{\mR^d}(\cos(\xi\cdot \sigma_tz)-1)\nu^{(\alpha)}(\dif z)\right)\\
\leq&\exp\left(-{t}|\xi|^\alpha\int_0^{\infty}\int_{\mS^{d-1}}\frac{1-\cos(\xi/|\xi|\cdot \sigma_t r\theta)}{r^{1+\alpha}}\Sigma(\dif \theta)\dif r\right)
\leq \e^{-ct|\xi|^{\alpha}},
\end{align*}
where the constant $c>0$ depends only on $\alpha$ and $\Sigma(\mS^{d-1})$. Hence, by \cite{Sa99}, Proposition 28.1, the random variable $L_t^\sigma$ defined by \eqref{eq:L} admits a smooth density $p^\sigma(t,x)$ given by Fourier's inverse transform
$$
p^\sigma(t,x)=(2\pi)^{-d/2}\int_{\mR^d}\e^{-i x\cdot \xi}\mE \e^{i\xi\cdot L^\sigma_{t} }\dif \xi,\ \ \forall t>0,
$$
and the partial derivatives of $p^\sigma(t,\cdot)$ at any orders  tend to $0$ as $|x|\to \infty$. 

We need the following heat kernel estimates in integral form with Littlewood-Paley’s decomposition, which is obtained in \cite{CHZ20}, Lemma 3.3 (see also \cite{HWW20}, Lemma 2.12). 

\bl[Heat kernel estimates]
Suppose that $\bf (\bH0)$ holds with constant $a_0\in(0,1)$. Let $p_{s,t}^\sigma$ be the density of the random variable $L_t^\sigma - L_s^\sigma $.
\begin{enumerate}[(i)]
\item {There is a constant $c>0$ such that for all $0\leq s\leq t < \infty$,
\begin{align}\label{eq:XM50}
\| \cR_{-1} p^\sigma_{s,t} \|_1\leq c.
\end{align} }
\item For every $ \gamma\in[0,\alpha)$ and $ \vartheta\ge  \gamma$, there is a constant $c>0$ such that for all $0\leq s < t< \infty$ and $j\in \mN_0$,
\begin{align}\label{CC01}
\int_{\mR^d}|x|^{\gamma}|\cR_jp^\sigma_{s,t}(x)|\dif x \lesssim_c  (t-s)^{-\frac{\vartheta-\gamma}{\alpha}}2^{-j\vartheta},
\end{align}
where the block operators $\cR_j$ are defined by \eqref{eq:Block}. 

\item For each $\gamma \in [0,\alpha)$ and $T>0$, there is a constant $c>0$ such that for any $s,t \in [0,T]$ and $ j \in \mN_0$,
\begin{align}\label{CC01z}
\int_0^t\int_{\mR^d} |x|^{\gamma}|\cR_jp^\sigma_{s,t}(x)|\dif x\dif s \lesssim_c  2^{-j\alpha}.
\end{align}
\end{enumerate}
\el

\section{Nonlocal equations with singular L\'evy measures}\label{sec:PDE}
 
Fix $\alpha \in (1,2)$ and $T>0$. Let $\sigma(t,x): \mR_+ \times \mR^d \to \mM_{non}^d$ be a bounded measurable function and $\lambda\ge0$ be a real number. In this section, we consider the following non-local parabolic equation with time-dependent variable diffusion coefficient $\sigma_t(x) := \sigma(t,x)$:
\begin{align}\label{PDE}
\p_t u=\sL^{(\alpha)}_{\sigma}   u-\lambda u+b\cdot\nabla u+f,\ \ u(0)=0,
\end{align}
where $b,f\in L^\infty_T\bC^{-\beta}$ with some $\beta\in(0,1)$, and $\sL^{(\alpha)}_{\sigma}$ is given by \eqref{eq:nolocal} with $\nu^{(\alpha)}$ being defined by \eqref{CC02} and satisfying condition {\bf (ND)}.

Our aim in this section is to establish Schauder's estimate and obtain the well-posedness for PDE \eqref{PDE}. Before this, we need to give the definitions of solutions (see Definition \ref{AA05} \& \ref{AA05:00} below). We first state the main results of this section.

\bt\label{BB01}
Let $T>0$, and $\alpha\in(1,2)$, and $\beta\in(0,\frac{\alpha-1}{2})$. Assume that {\bf (H$^\sigma$)} holds with some constant $c_0$, and $b,f\in L^\infty_T\bC^{-\beta}$. For any $\lambda\ge0$, there is a unique solution $u$ to PDE \eqref{PDE} in the sense of Definition \ref{AA05} satisfying for any $\gamma>\frac{\alpha}{\alpha-1-2\beta}$ and $\theta\in[0,\alpha)$,
\begin{align}\label{Hzm01}
 (1+\lambda)^{\frac{\theta}{\alpha}} \|u\|_{L^\infty_T\bC^{\alpha-\beta-\theta}} \lesssim_c  \ell_{b,\beta}^\gamma \|f\|_{L^\infty_T\bC^{-\beta}} ,
\end{align}
where $\ell_{b,\beta}:=1+\|b\|_{L^\infty_T\bC^{-\beta}}$, and $c>0$ is a constant depending only on $  \Theta, \gamma,\theta$. 
\et

\bt\label{BB01:00}
Let $T>0$, and $\alpha\in(1,2)$, and $\beta\in[\frac{\alpha-1}{2}, \alpha-1 )$. Assume that {\bf (H$^\sigma$)} holds with some constant $c_0$, and $b,\div b,f\in L^\infty_T\bC^{-\beta}$. For any $\lambda \geq 0$, there is a unique solution $u$ to PDE \eqref{PDE} in the sense of Definition \ref{AA05:00} satisfying for any $\gamma>\frac{\alpha}{\alpha-1-\beta}$ and $\theta\in[0,\alpha)$,
\begin{align}\label{Hzm01:00}
(1+\lambda)^{\frac{\theta}{\alpha}} \|u\|_{L^\infty_T\bC^{\alpha-\beta-\theta}} \lesssim_c  \ell_{b,\beta}^\gamma \|f\|_{L^\infty_T\bC^{-\beta}},
\end{align}
where $\ell_{b,\beta}:=1+\|b\|_{L^\infty_T\bC^{-\beta}}+\|\div b\|_{L^\infty_T\bC^{-\beta}}$,  and $c>0$ is a constant only depending on $  \Theta, \gamma,\theta$. 
\et

\subsection{Well-definedness of solutions}\label{sec:3.1}

First of all, it is well-known that the product of two distributions is not always meaningful. Fortunately, thanks to \eqref{S2:Bern} and \eqref{eq:pro}, when $b(t) \in \bC^{-\beta}$ and $u(t) \in \bC^{1+\beta+\varepsilon}$ with $\varepsilon>0$, it is easy to see that $b\cdot \nabla u (t) :\bC^{-\beta}\times\bC^{\beta+\eps}$ is well-defined. Thus, we have the following definition of strong solutions to PDE \eqref{PDE}.

\bd[Strong solutions]\label{AA05}
Let $\alpha\in(1,2)$, $\beta\in(0,1)$, $T>0$, and $\lambda \geq 0$. Assume that {\bf (H$^\sigma$)} holds with some constant $c_0$. For any $b,f\in L^\infty_T\bC^{-\beta}$, we call a function $u\in\bigcup_{\eps>0} L^\infty_T\bC^{1+\beta+\eps}$ a strong solution to PDE \eqref{PDE} on $[0,T]$, if for any $t\in[0,T]$,
\begin{align}\label{AA03}
\begin{split}
u(t)
= \int_0^t\([\sL_{\sigma}^{(\alpha)}-\lambda] u+b\cdot\nabla u+f\)(s)\dif s,
\end{split}
\end{align}
where $\sL^{(\alpha)}_{\sigma}$ is defined by \eqref{eq:nolocal}.
\ed

On the other hand, according to Schauder's estimates, one sees that the best regularity of the solution  $u(t)$ is in $\bC^{\alpha-\beta}$. However, we also know that the domain of the operator $\sL^{(\alpha)}_{\sigma_t} $ is $\bigcup_{\eps>0}\bC^{\alpha+\eps}$. Therefore, in order to define the solutions of PDE \eqref{PDE}, we need to extend the domain of $\sL^{(\alpha)}_{\sigma_t} $.  
 
\bl[Boundedness of operator $\sL_\sigma^{(\alpha)}$]\label{Lem31}
Let $\alpha\in(1,2)$ and $\beta\in(0,1)$. Assume that the condition {\bf(H$^\sigma$)} holds for some constant $c_0$. Then there is a constant $c=c(d,\alpha,\beta,c_0)>0$ such that for any $u\in C^\infty_b(\mR^d)$ and $t\ge0$,
\begin{align}\label{EE01}
\|\sL_{\sigma}^{(\alpha)} u(t)\|_{\bC^{-\beta}} \lesssim_c  \|u(t)\|_{\bC^{\alpha-\beta}}.
\end{align}
\el

We put the proof of this lemma in the Appendix for the fluency in reading.

\br
Based on Lemma \ref{Lem31}, we extend the domain of the linear operator $\sL_{\sigma_t}^{(\alpha)}$ from $\bC^{\alpha+\eps}$ to $\bC^{\alpha-\beta}$ with $ \beta\in(0,1)$. 
Thus, when we consider Definition \ref{AA05} as a solution, we should have $L_T^\infty\bC^{\alpha - \beta} \subset L_T^\infty\bigcup_{\varepsilon>0} \bC^{1+ \beta+\varepsilon}$, which derives that $\alpha -\beta > 1+\beta$. Thus, under $\beta<\frac{\alpha-1}{2}$, based on \eqref{EE01} and \eqref{eq:pro}, every term in \eqref{AA03} makes sense. Then we can state Theorem \ref{BB01}.
\er

Furthermore, a natural question arises: can the condition $\beta\in (0,\frac{\alpha-1}{2})$ be relaxed to allow a wider range of values for $\beta$? To this end, we introduce the following product (see also \cite{HZ23, Per23}): if ${\rm div} b\in L^\infty_T \bC^{-\beta}$, we define
\begin{align}\label{08:21}
b\cdot\nabla u=\{\div(b\preceq u)+ b \succ \nabla u\}-\div b\preceq u=: b \odot \nabla u-\div b\preceq u.
\end{align}
In particular, if $\div b=0$, then $b \cdot \nabla u = b \odot \nabla u$. Based on the estimates of paraproducts in Lemma \ref{productlaw}, we obtain the following result. As the proof is available in \cite[Lemma 2.6]{HZ23}, we omit it here.
\bl\label{Le26}
For $\beta\in(0,1)$ , there is a constant $c=c(d,\beta)>1$ such that
\begin{align}\label{08:20} 
\|b\odot\nabla u\|_{L^\infty_T\bC^{-\beta}}\lesssim_c \|b\|_{L^\infty_T\bC^{-\beta}}\| u\|_{L_T^\infty\bB^1_{\infty,1}}
\end{align}
and
\begin{align}\label{08:210}
\|\div b\preceq u\|_{L^\infty_T\bC^{-\beta}}\lesssim_c
\|\div b\|_{L^\infty_T\bC^{-\beta}}\|u\|_{L^\infty_T\bB^\beta_{\infty,1}}.
\end{align}
\el

Based on Lemma \ref{Le26}, the product $b\cdot\nabla u$ makes sense as long as $b,\div b \in L^\infty_T\bC^{-\beta}$ and $u\in L^\infty_T \bC^{1+\eps}$ with $\varepsilon >0$. Now we have the following definition of paraproduct solutions to PDE \eqref{PDE}.

\bd[Paraproduct solutions]\label{AA05:00}
Let $\alpha\in(1,2)$, $\beta\in(0,1)$, $T>0$, and $\lambda \geq 0$. Assume that {\bf (H$^\sigma$)} holds with some constant $c_0$. For any $b, f\in L^\infty_T\bC^{-\beta}$ with $\div b \in L^\infty_T\bC^{-\beta}$, we call a function $u\in\bigcup_{\eps>0} L^\infty_T\bC^{1+\eps}$ a paraproduct solution to PDE \eqref{PDE} on $[0,T]$, if for any $t\in[0,T]$,
\eqref{AA03} holds, where $b\cdot\nabla u$ is defined as \eqref{08:21}.
\ed

\br 
In this case, $\sL_\sigma^{(\alpha)} u$ and $ b\cdot\nabla u$ are well-defined if $u\in L^\infty_T\bC^{\alpha-\beta} \subset \bigcup_{\varepsilon>0} L^\infty_T\bC^{1+\eps}$, which is $\alpha-\beta>1$, i.e., $\beta<\alpha-1$. Then we can state Theorem \ref{BB01:00}
\er

\subsection{Zero-drift case}
In this subsection, we assume $b\equiv0$ and investigate the following equation with $\lambda\ge0$:
\begin{align}\label{PDElamz}
\p_t w=\sL^{(\alpha)}_{\sigma}  w-\lambda w+f,\ \ w(0)=0,
\end{align}
where $f,w_0$ are smooth functions. We aim to establish the following result.  

\bp\label{Lem01}
Fix $T>0$. Let $\alpha\in(1,2)$, $\beta\in(0,1)$, and $f$ be a smooth function. Assume {\bf (H$^\sigma$)} holds with constant $c_0$. If $w$ is a solution to PDE \eqref{PDElamz}, then for any $\theta\in[0,\alpha)$, there is a constant $c>0$ depending on $ \Theta,\theta$ such that for all $\lambda\ge0$, 
\begin{align}\label{BB06z}
(\lambda+1)^{\frac{\theta}{\alpha}}\|w\|_{L^\infty_T\bC^{\alpha-\beta-\theta}}
& \le c \|f\|_{L^\infty_T\bC^{-\beta}}.
\end{align}
\ep

To prove this result, 
we only need to establish a priori estimates \eqref{BB06z}. 
Fix $y\in\mR^d$. For any function $h$, let 
$$
h^y(x):=h(x+y).
$$
Then we have
\begin{align*}
\p_t w^y=\sL^{(\alpha)}_y w^y-\lambda w^y+ ( f^y+\sA w^y),
\end{align*}
where $\sA:=\sL^{(\alpha)}_{\sigma^y}-\sL^{(\alpha)}_y$ and 
\begin{align*}
\sL^{(\alpha)}_y g(t,x)& :=\sL^{(\alpha)}_{\sigma^y(0)} g(t,x)\\
& =\int_{\mR^d}\(g(t,x+\sigma_t^y(0) z) -g(t,x) -\sigma_t^y(0) z\cdot\nabla g(t,x)\)\nu^{(\alpha)}(\dif z).
\end{align*}
By  Subsection \ref{sec:op} and \cite{CHZ20}, Section 3, one sees that the operator $\p_t- (\sL^{\alpha}_y-\lambda)$ associates to a semigroup $\e^{-\lambda(t-s)}P_{s,t}^{\sigma^y(0)}$, i.e., 
\begin{align*}
w^y(t,x)=\int_0^t \e^{-\lambda(t-s)}P_{s,t}^{\sigma^y(0)} \(f^y+\sA w^y\)(s,x)\dif s,
\end{align*}{
and $p_{s,t}^{\sigma^y(0)}$ is the heat kernel of $P_{s,t}^{\sigma^y(0)}$. For the sake of simplicity, in the sequel, we write $P_{s,t}^{(y)}=P_{s,t}^{\sigma^y(0)}$ and $p^{(y)}_{s,t}=p_{s,t}^{\sigma^y(0)}$.}
Thus, for $j\ge-1$, acting on both sides of the above equations by $\cR_j$, we get that
\begin{align}\label{BB04}
\begin{split}
\cR_j w(t,y) = \cR_j w^y(t,0) =  \int_0^t \e^{-\lambda(t-s)}\cR_jP_{s,t}^{(y)} \(f^y+\sA w^y  \)(s,0)\dif s.
\end{split}
\end{align}

Let us separately estimate the terms on the right-hand side of \eqref{BB04}. 

\bl
Fix $T>0$. Suppose that $\alpha\in(1,2)$. Assume that {\bf (H$^\sigma$)} holds with constant $c_0$. Then for any $\vartheta\in[0,\alpha]$ and $\beta\in (0,1)$, there is a constant $c>0$ such that for all $ j \ge-1$, $s,t \in [0,T]$, $y \in \mR^d$, and $\lambda \geq 0$,
\begin{align}\label{AA07}
\int_0^t \e^{-\lambda(t-s)}\Big|\cR_jP_{s,t}^{(y)}f^y(s,0)\Big|\dif s \lesssim_c 2^{-(\alpha-\beta-\vartheta)j}( \lambda+1)^{-\frac{\vartheta}{\alpha}} \|f\|_{L^\infty_T\bC^{-\beta}};
\end{align}
and for any $\eta\in [0,\alpha)$, $\varepsilon\in (0,\alpha-\eta)$, there is also a constant $c>0$ such that for all $ j \ge-1$, $s,t \in [0,T]$, $y \in \mR^d$, and $\lambda \geq 0$,
\begin{align}\label{AA06}
\int_0^t \e^{-\lambda(t-s)}\Big|\cR_jP_{s,t}^{(y)}\sA w^y(s,0)\Big|\dif s
&\lesssim_c 2^{-\eta j} (\lambda+1)^{-\frac{\alpha-\eta-\eps}{\alpha}} \|w\|_{L^\infty_T\bC^{\alpha-1+\eps}}.
\end{align}
\el

\begin{proof}
We only prove the case $j \geq 0$, since the case $j=-1$ is similar and easier by \eqref{eq:XM50}. 

\medskip\noindent
(i) For the first inequality, by Remark \ref{rem:AA01} and \eqref{CC01}, we have that for any $\vartheta\in(0,\alpha]$,
\begin{align*}
\int_0^t \e^{-\lambda(t-s)}\Big|\cR_jP_{s,t}^{(y)}f^y(s,0)\Big|\dif s&\le \int_0^t \e^{-\lambda(t-s)}\int_{\mR^d}|\cR_jp_{s,t}^{(y)}(x)\widetilde\cR_jf(s,x+y)|\dif x\dif s\nonumber\\
& \leq  \int_0^t \e^{-\lambda(t-s)} \| \cR_j p^{(y)}_{s,t} \|_1 \| \widetilde \cR_j f(s) \|_\infty \dif s\nonumber\\
&\lesssim 2^{\beta j} \int_0^t \e^{-\lambda (t-s)}(t-s)^{-\frac{\alpha-\vartheta}{\alpha}}2^{-j(\alpha-\vartheta)}\|f(s)\|_{\bC^{-\beta}}\dif s\nonumber\\
&\lesssim 2^{-(\alpha-\beta-\vartheta)j}(\lambda+1)^{-\frac{\vartheta}{\alpha}} \|f\|_{L^\infty_T\bC^{-\beta}} \Gamma(\vartheta/\alpha),
\end{align*}
where the last inequality is provided by the change of variable; and for $\vartheta=0$, similarly,  by Remark \ref{rem:AA01} and \eqref{CC01z},
\begin{align*}
\int_0^t \e^{-\lambda(t-s)}\Big|\cR_jP_{s,t}^{(y)}f^y(s,0)\Big|\dif s
&
\lesssim 2^{\beta j}\|f \|_{L^\infty_T\bC^{-\beta}} \int_0^t \|\cR_jp_{s,t}^{(y)}\|_1  \dif s\nonumber\\
&\lesssim 2^{-(\alpha-\beta)j} \|f\|_{L^\infty_T\bC^{-\beta}}.
\end{align*}

\medskip\noindent
(ii) For the second one, applying \cite{CHZ20}, Lemma 4.4, to $\e^{-\lambda(t-s)}w(s)$, we have that for any $T>0$,  and $\eta\in [0,\alpha)$, $\varepsilon\in (0,\alpha-\eta)$,
\begin{align*}
\int_0^t \e^{-\lambda(t-s)}\Big|\cR_jP_{s,t}^{(y)}\sA w^y(s,0)\Big|\dif s&\lesssim  2^{-\eta j}\int_0^t \e^{-\lambda(t-s)}(t-s)^{-\frac{\eta+\eps}{\alpha}}\|w(s)\|_{\bC^{\alpha-1+\eps}}\dif s\nonumber\\ 
&\lesssim 2^{-\eta j}(\lambda+1)^{-\frac{\alpha-\eta-\eps}{\alpha}} \|w\|_{L^\infty_T\bC^{\alpha-1+\eps}} \Gamma(\tfrac{\alpha-\eps-\eta}{\alpha}),
\end{align*}
where we used the change of variable. Combining the calculations above, we finish the proof.
\end{proof}

Now, we give the

\begin{proof}[Proof of Proposition \ref{Lem01}]
{\bf (Step 1)} First of all, we prove \eqref{BB06z} holds for large enough $\lambda$. Notice that, by \eqref{AA07}, for $\theta\in[0,\alpha]$,
$$
\int_0^t \e^{-\lambda(t-s)}\Big|\cR_jP_{s,t}^{(y)}f^y(s,0)\Big|\dif s \lesssim 2^{-(\alpha-\beta -\theta)j}(  \lambda+1)^{-\frac{\theta}{\alpha}} \|f\|_{L^\infty_T\bC^{-\beta}}.
$$
Moreover, by \eqref{AA06} with  $\eta = \alpha -\beta -\theta\wedge (\alpha-\beta)\in[0,\alpha)$ and $\eps\in (0,\beta+\theta\wedge (\alpha-\beta))$, we have that 
\begin{align*}
    \int_0^t \e^{-\lambda(t-s)}\Big|\cR_jP_{s,t}^{(y)}\sA w^y(s,0)\Big|\dif s
&\lesssim  2^{-(\alpha-\beta-\theta\wedge (\alpha-\beta))j} (\lambda+1)^{-\frac{\beta+\theta\wedge (\alpha-\beta)-\eps}{\alpha}} \|w\|_{L^\infty_T\bC^{\alpha-1+\eps}}\\
&\lesssim  2^{-(\alpha-\beta-\theta)j} (\lambda+1)^{-\frac{\theta}{\alpha}} \|w\|_{L^\infty_T\bC^{\alpha-1+\eps}},
\end{align*}
 provided that $\eps \in (0,\beta \wedge (\alpha-\theta)]$ with $\theta \in [0,\alpha)$, and noting that $
 \beta+\theta\wedge (\alpha-\beta) - \theta = \beta \wedge (\alpha-\theta)$. Thus, taking supremum of $y$ in \eqref{BB04}, we have that there is a constant $c>0$ such that for all $\lambda \geq 0$,
\begin{align}\label{Hzm0z}
\|w(t)\|_{\bC^{\alpha-\beta-\theta}} & = \sup_{j \geq -1}2^{(\alpha-\beta-\theta)j} \|\cR_j w(t)\|_\infty \nonumber \\
& \le c(\lambda+1)^{-\frac{\theta}{\alpha}}
\left( \|f\|_{L^\infty_T\bC^{-\beta}} +  \|w\|_{L^\infty_T\bC^{\alpha-1+\eps}}\right).
\end{align}

Now we need to estimate the term $\|w\|_{L^\infty_T\bC^{\alpha-1+\eps}}$. 
To this end, we take $\eps \in (0,\beta \wedge (\alpha -\theta)\wedge (1-\beta)]$, and $\theta=1-\beta-\eps\in[0,\alpha)$, and $\lambda_0 :=(2c)^{\alpha/\theta}$. Then by \eqref{Hzm0z}, one sees that for all $\lambda \geq \lambda_0$, 
\begin{align}\label{revise:00}
\|w\|_{L^\infty_T\bC^{\alpha-1+\eps}}
 \lesssim  (\lambda+1)^{-\frac{1-\beta-\eps}{\alpha}}\|f\|_{L^\infty_T\bC^{-\beta}}\lesssim  \|f\|_{L^\infty_T\bC^{-\beta}}.
\end{align}
Substituting this into \eqref{Hzm0z} with $\eps \in (0,\beta \wedge (\alpha -\theta)\wedge (1-\beta)]$, we prove the desired claim for all $\lambda \geq \lambda_0$.

\medskip\noindent
{\bf (Step 2)} In order to finish the proof, we next prove that \eqref{revise:00} holds for $\lambda\in[0,\lambda_0]$. If we show this, then as the same proof in {\bf (Step 1)}, combining with \eqref{Hzm0z}, we get the desired estimates. Now we divide the solution $w$ into two terms: $w=  w_1+w_2$ with
\begin{align}
    &\p_t w_1=\sL^{(\alpha)}_\sigma w_1-(\lambda+\lambda_0)w_1+f,\quad w_1(0)=0,\label{eq:WQ00}\\
    &\p_t w_2=\sL^{(\alpha)}_\sigma w_2- \lambda w_2+\lambda_0 w_1,\quad \quad w_2(0)=0.\label{eq:WQ01}
\end{align}
By \eqref{eq:WQ01}, based on the maximum principle, one sees that 
\begin{align*}  \|w_2\|_{\mL^\infty_T}\le T\lambda_0\|w_1\|_{\mL^\infty_T},
\end{align*}
which together with \eqref{eq:WQ00} and \eqref{revise:00}, derives that
\begin{align}\label{revise:01}
    \|w\|_{\mL^\infty_T}\le \|w_1\|_{\mL^\infty_T}+\|w_2\|_{\mL^\infty_T}
    \lesssim \|w_1\|_{\mL^\infty_T}
    \overset{\eqref{revise:00}}{\lesssim}\|f\|_{L^\infty_T\bC^{-\beta}}.
\end{align}
Hence, by \eqref{Hzm0z} with $\theta=0$ and $\eps \in (0,\beta  \wedge (1-\beta))$, in view of the interpolation inequality \eqref{eq:inter}, we have that for all $\lambda \geq 0$,
\begin{align*}
\|w\|_{L^\infty_T\bC^{\alpha-\beta}}
    &\le c\|f\|_{L^\infty_T\bC^{-\beta}}+c\|w\|_{L^\infty_T\bC^{\alpha-1+\eps}}\\
    &\le c\|f\|_{L^\infty_T\bC^{-\beta}}+\left(\frac{1}{2} \|w\|_{L^\infty_T\bC^{\alpha-\beta}}+c\|w\|_{\mL^\infty_T}\right)\\
    & \overset{\eqref{revise:01}}{\leq} c \|f\|_{L^\infty_T\bC^{-\beta}} +\frac{1}{2} \|w\|_{L^\infty_T\bC^{\alpha-\beta}},
\end{align*}
which implies that
$\|w\|_{L^\infty_T\bC^{\alpha-\beta}} \lesssim \|f\|_{L^\infty_T\bC^{-\beta}} $. Thus, for any $\eps \in (0, 1-\beta]$, and $\lambda \in [0,\lambda_0]$, one sees that
\begin{align*}
\|w\|_{L^\infty_T\bC^{\alpha-1+\eps}}\lesssim
\|w\|_{L^\infty_T\bC^{\alpha-\beta}} \lesssim \|f\|_{L^\infty_T\bC^{-\beta}}  
\end{align*}
and completes the proof.
\end{proof}

\subsection{Distributional-drift case}

In this section, we focus on showing Theorems \ref{BB01} and \ref{BB01:00}. We  frequently use the following fact:  by \eqref{eq:pro}, and Bernstein's inequality \eqref{S2:Bern}, and \eqref{08:21}-\eqref{08:210}, one sees that for any $\eps>0$,
\begin{align}\label{eq:YY01}
\|b\cdot\nabla u(s)\|_{\bC^{-\beta}}\lesssim
\begin{cases}   \|b\|_{L_T^\infty\bC^{-\beta}} \|u(s)\|_{\bC^{ 1+\beta+\eps}},\\(\|b\|_{L_T^\infty\bC^{-\beta}} +\|\div b\|_{L_T^\infty\bC^{-\beta}} )\|u(s)\|_{\bC^{ 1+\eps}}.
\end{cases}
\end{align}

The following a priori estimate is crucial. 

\bp\label{Lem01z}
Under the same conditions in Theorem \ref{BB01} (resp. Theorem \ref{BB01:00}), if $u$ is a strong solution (resp. paraproduct solution) to PDE \eqref{PDE} in the sense of Definition \ref{AA05} (resp. Defnition \ref{AA05:00}), then for any $\gamma > \frac{\alpha}{\alpha-1-2\beta}$ (resp. $\gamma > \frac{\alpha}{\alpha-1-\beta}$) and $\theta\in[0,\alpha)$, there is a constant $c=c(\Theta, \gamma,\theta)>0$ such that for all $\lambda\ge\lambda_b:=c\ell_{b,\beta}^{\gamma} $, 
\begin{align}\label{BB06zz}
(\lambda+1)^{\frac{\theta}{\alpha}}\|u\|_{L^\infty_T\bC^{\alpha-\beta-\theta}}
& \le c  \|f\|_{L^\infty_T\bC^{-\beta}},
\end{align}
where $\ell_{b,\beta}:= 1+ \|b\|_{L_T^\infty\bC^{-\beta}}$ (resp. $\ell_{b,\beta}:=1+ \|b\|_{L_T^\infty\bC^{-\beta}} + \| \div b\|_{L_T^\infty\bC^{-\beta}}$) for $\beta\in(0,\frac{\alpha-1}{2})$ (resp. $\beta\in[\frac{\alpha-1}{2},\alpha-1)$).
\ep

\begin{proof}
First of all, for any $\theta\in[0, \alpha)$, by Proposition \ref{Lem01}, one sees that for all $\lambda\ge  0$,
\begin{align}
\|u\|_{L^\infty_T\bC^{\alpha-\beta-\theta}}
& \lesssim   (\lambda+1)^{-\frac{\theta}{\alpha}}\|f+b\cdot\nabla u\|_{L^\infty_T\bC^{-\beta}}\nonumber
\nonumber\\
&\lesssim  (\lambda+1)^{-\frac{\theta}{\alpha}}\left(\|f\|_{L^\infty_T\bC^{-\beta}}+\ell_{b,\beta} \|u\|_{L^\infty_T\bC^{1+\beta+\eps}}\right),\label{Hzm02z}
\end{align}
where we used \eqref{eq:YY01} in the last inequality for $\eps>0$. 

\medskip\noindent
{\bf (Case 1: $\beta<\frac{\alpha-1}{2}$)}   From \eqref{Hzm02z}, in particular, taking $\eps \in (0,\alpha -1-2\beta]$ small enough such that $\frac{\alpha}{\alpha-1-2\beta-\eps}<\gamma$, and letting  $\theta=\alpha -1-2\beta-\eps$, we infer that for all $\lambda\geq0$,
\begin{align*}
\|u\|_{L^\infty_T\bC^{1+\beta+\eps}} 
& \lesssim_c  (\lambda+1)^{-\frac{1}{\gamma}}\|f\|_{L^\infty_T\bC^{-\beta}}+(\lambda+1)^{-\frac1\gamma}\ell_{b,\beta} \|u\|_{L^\infty_T\bC^{1+\beta+\eps}},
\end{align*}
which implies that for any $\lambda\geq \lambda_b:= \left(2c\ell_{b,\beta}\right)^{\gamma}$,
\begin{align*}
\tfrac{1}{2}
\|u\|_{L^\infty_T\bC^{1+\beta+\eps}} 
\lesssim_c (\lambda+1)^{-\frac{1}{\gamma}}\|f\|_{L^\infty_T\bC^{-\beta}} \lesssim \ell_{b,\beta}^{-1} \|f\|_{L^\infty_T\bC^{-\beta}}.
\end{align*}
Substituting this into \eqref{Hzm02z}, we obtain the desired estimates \eqref{BB06zz}. 

\medskip\noindent
{\bf (Case 2: $\beta \in[\frac{\alpha-1}{2}, \alpha-1$))} Similarly, from \eqref{Hzm02z}, taking $\eps>0$ small enough so that $\frac{\alpha}{\alpha-1-\beta-\eps}< \gamma$, and letting $\theta=\alpha-1-\beta-\eps$, we have
 \begin{align*}
\|u\|_{L^\infty_T\bC^{1+\eps}}
\lesssim  (\lambda+1)^{-\frac{1}{\gamma}}\left(\|f\|_{L^\infty_T\bC^{-\beta}}+\ell_{b,\beta}\|u\|_{L^\infty_T\bC^{1+\eps}}\right).
\end{align*}
The remainder is the same as the proof in {\bf (Case 1)}, so we omit it. 
\end{proof}

\subsubsection{Proof of Theorem \ref{BB01}}

We divide the proof into two steps. Firstly, in {\bf (Step 1)}, we show \eqref{Hzm01} as a priori estimate. Then we give the existence and uniqueness of a strong solution in {\bf (Step 2)}. 

\medskip\noindent
{\bf (Step 1)} Letting $u$ be a strong solution to PDE \eqref{PDE}, we show that $u$ satisfies the estimates \eqref{Hzm01} in this step.  

\medskip
First of all, let $\lambda_b$ be the same constant in Proposition \ref{Lem01z}, and $u^{\lambda_b}$ be a strong solution to PDE \eqref{PDE} with coefficient $=\lambda+\lambda_b$ and $u^{\lambda_b}(0)=0$. Denote by $\widetilde u:=u-u^{\lambda_b}$. Noticing that
\begin{align*}
\p_t \widetilde u=\sL^{(\alpha)}_{\sigma}  \widetilde u -  \lambda \widetilde u +b\cdot\nabla \widetilde u+\lambda_b u^{\lambda_0}, \quad \widetilde u(0)=0,
\end{align*}
by maximum principle, we have
\begin{align*}
\|u\|_{\mL_T^\infty} \leq 
\|\widetilde u\|_{\mL_T^\infty} + \| u^{\lambda_b}\|_{\mL_T^\infty} &\le (1+T\lambda_b )\|u^{\lambda_b}\|_{\mL_T^\infty}\lesssim \ell_{b,\beta}^\gamma \|u^{\lambda_b}\|_{\mL_T^\infty},
\end{align*}
which together with \eqref{BB06zz} yields
\begin{align}\label{revise:02}  \|u\|_{\mL_T^\infty} 
\lesssim  (1+\lambda+\lambda_b)^{-\frac{\alpha-\beta}{\alpha}} \ell_{b,\beta}^\gamma\|f\|_{L^\infty_T\bC^{-\beta}} 
 \lesssim
\ell_{b,\beta}^{\frac{\gamma\beta}{\alpha}}\|f\|_{L^\infty_T\bC^{-\beta}}.
\end{align}

\medskip
Next, following the proof of \eqref{Hzm02z}, one sees that for any $\eps>0$ and $\theta\in[0,\alpha)$, and all $\lambda\ge0$,
\begin{align*}
\|u\|_{L^\infty_T\bC^{\alpha-\beta-\theta}}
\lesssim_c  (\lambda+1)^{-\frac{\theta}{\alpha}}\left(\|f\|_{L^\infty_T\bC^{-\beta}}+ \ell_{b,\beta} \|u\|_{L^\infty_T\bC^{1+\beta+\eps}}\right),
\end{align*}
which, by interpolation inequality, and Young's inequality, yields that
\begin{align*}
    (1+\lambda)^{\frac{\theta}{\alpha}}\|u\|_{L^\infty_T\bC^{\alpha-\beta-\theta}}
& \leq c
\|f\|_{L^\infty_T\bC^{-\beta}}+c \left( \ell_{b,\beta}  \|u\|_{\mL^\infty_T}^{\frac{\alpha-1-2\beta-\eps}{\alpha-\beta}}\|u\|_{L^\infty_T\bC^{\alpha-\beta}}^{\frac{1+\beta+\eps}{\alpha-\beta}}\right) \\
&\leq c  \|f\|_{L^\infty_T\bC^{-\beta}}+c \ell_{b,\beta} ^{\frac{\alpha-\beta}{\alpha-1-2\beta-\eps}}\|u\|_{\mL^\infty_T}+\tfrac12\|u\|_{L^\infty_T\bC^{\alpha-\beta}},
\end{align*}
where we used the fact $1+\beta<\alpha -\beta$ and took $\eps>0$ small enough. Letting $\eps$ small enough such that $\frac{\alpha-\beta}{\alpha-1-2\beta-\eps}\le \frac{\alpha-\beta}{\alpha}\gamma$, by \eqref{revise:02}, we have
\begin{align}
    (1+\lambda)^{\frac{\theta}{\alpha}}\|u\|_{L^\infty_T\bC^{\alpha-\beta-\theta}}&\leq c  \|f\|_{L^\infty_T\bC^{-\beta}}
    \ell_{b,\beta}^{\frac{\alpha-\beta}{\alpha}\gamma+\frac{\gamma\beta}{\alpha}}+\tfrac12\|u\|_{L^\infty_T\bC^{\alpha-\beta}}\nonumber\\
&=c\|f\|_{L^\infty_T\bC^{-\beta}}\ell_{b,\beta}^{\gamma}+\tfrac12\|u\|_{L^\infty_T\bC^{\alpha-\beta}},\label{revise:03}
\end{align}
which by taking $\theta=0$ gives that
\begin{align*}  \|u\|_{L^\infty_T\bC^{\alpha-\beta-\theta}}\le c\|f\|_{L^\infty_T\bC^{-\beta}}\ell_{b,\beta}^{\gamma}.
\end{align*}
Substituting this into \eqref{revise:03}, we obtain \eqref{Hzm01}.

\medskip\noindent
{\bf (Step 2)} The uniqueness is directly from \eqref{Hzm01}, and we only need to show the existence.
Let $\rho_m(\cdot):=m^d\rho(m\cdot)$ be the usual mollifier with $\rho\in C^\infty_c(\mR^d)$ and 
\begin{align*}
b_m(t):=b(t)*\rho_m,\ \ f_m(t):=f(t)*\rho_m.
\end{align*}
Then, for $b_m,f_m\in L^\infty_TC^\infty_c (\mR^d)$, and any $\gamma>\beta$,
\begin{align*}
\lim_{m\to\infty}\(\|b_m-b\|_{L_T^\infty\bC^{-\gamma}}+\|f_m-f\|_{L_T^\infty\bC^{-\gamma}}\)=0.
\end{align*}
Moreover,
\begin{align}\label{eq:ZY00}
\sup_{m}\|b_m\|_{L_T^\infty \bC^{-\beta}}\le\|b\|_{L_T^\infty\bC^{-\beta}},\ \ \sup_{m}\|f_m\|_{L_T^\infty\bC^{-\beta}}\le\|f\|_{L_T^\infty\bC^{-\beta}}. 
\end{align}

\medskip
 
Let $u_m$ be the classical solution of PDE \eqref{PDE} with $(b,f)=(b_m,f_m)$. Noting that $u_n(0) - u_m(0)=0$ and
\begin{align*}
\p_t (u_n -u_m) = & \sL^{(\alpha)}_\sigma ( u_n -u_m) - \lambda ( u_n -u_m) + b_n \cdot \nabla ( u_n - u_m) \\
& + \left[  (b_n - b_m)\cdot \nabla u_m+ f_n -f_m\right],
\end{align*}
by the result of {\bf (Step 1)} (i.e., \eqref{Hzm01}) and \eqref{eq:ZY00}, we have that 	
\begin{align}
\| u_n -u_m \|_{L^\infty_T\bC^{\alpha-\beta}} 
& \lesssim \|  (b_n - b_m)\cdot \nabla u_m+ f_n -f_m \|_{L^\infty_T\bC^{-\beta}}\nonumber\\
& \lesssim \| b_n - b_m\|_{L^\infty_T\bC^{-\beta}} \| u_m \|_{L^\infty_T\bC^{\alpha-\beta}} + \| f_n -f_m \|_{L^\infty_T\bC^{-\beta}}\nonumber\\
& \lesssim  \| b_n - b_m\|_{L^\infty_T\bC^{-\beta}} + \| f_n -f_m \|_{L^\infty_T\bC^{-\beta}} \to 0,~~\text{as}~~ m,n \to \infty,\label{0728:new00}
\end{align}
where we used Lemma \ref{productlaw} and Bernstein's inequality in the second inequality since $\beta +1 < \alpha - \beta$. Hence, one sees that the limit of $\{u_{m}\}$, denoted by $u$, in $L_T^\infty \bC^{\alpha-\beta}$, exists.

\medskip

Next, we are going to show that $u$ is a solution of PDE \eqref{PDE}  in the sense of Definition \ref{AA05}. For any $\gamma >\beta$, by Lemma \ref{Lem31}, we have
\begin{align*}
\|\sL^{(\alpha)}_\sigma ( u_{m} - u)(t) \|_{\bC^{-\gamma}} & \lesssim \|\sL^{(\alpha)}_\sigma ( u_{m} - u)(t) \|_{\bC^{-\beta}} \\
& \lesssim \| u_{m}  - u\|_{L^\infty_T\bC^{\alpha-\beta}};
\end{align*}
and, by Lemma \ref{productlaw} and Bernstein's inequality, for $\gamma \in (\beta, \frac{\alpha -1}{2})$,
\begin{align*}
\| b_{m}  \cdot \nabla u_{m} (t) - b \cdot \nabla u(t) \|_{\bC^{-\gamma}}& \lesssim \| b_{m}  \cdot \nabla (u_{m}  -  u )\|_{L^\infty_T\bC^{-\gamma}} \\
&\qquad + \| ( b_{m}    - b ) \cdot \nabla u_{m}\|_{L^\infty_T\bC^{-\gamma}}\\
& \lesssim \| b \|_{L^\infty_T\bC^{-\gamma}}  \| u_{m}  - u\|_{L^\infty_T\bC^{\alpha-\gamma}} \\
&  \qquad +\| b_m -b \|_{L^\infty_T\bC^{-\gamma}} \| u \|_{L^\infty_T\bC^{\alpha-\gamma}}.
\end{align*}
Thus, by \eqref{eq:ZY00}, we get that for $\gamma \in (\beta, \frac{\alpha -1}{2})$, when $m \to \infty$,
$$
\int_0^{\cdot} \left( [\sL^{(\alpha)}_\sigma-\lambda] u_m + b_m \cdot \nabla u_m +f_m \right) (s)\dif s \to \int_0^{\cdot} \left( [\sL^{(\alpha)}_\sigma-\lambda] u + b \cdot \nabla u+f \right)(s) \dif s
$$
in $L^\infty_T\bC^{-\gamma}$. Finally, we see that $u (t)= \int_0^{t} \left( [\sL^{(\alpha)}_\sigma-\lambda] u + b \cdot \nabla u + f\right)(s) \dif s $.

\subsubsection{Proof of Theorem \ref{BB01:00}}

The proof of Theorem \ref{BB01:00} follows the same argument in the one of Theorem \ref{BB01}. So we only to show the estimate \eqref{Hzm01:00} by giving the main difference with  {\bf (Step 1)} in the proof of Theorem \ref{BB01}. 
  Let $\lambda_b$ be the same constant in Proposition \ref{Lem01z} and $u^{\lambda_b}$ be a strong solution to PDE \eqref{PDE} with coefficient $\lambda+\lambda_b$ and $u^{\lambda_b}(0)=0$. Denote by $\widetilde u:=u-u^{\lambda_b}$. 
Thus, by the same proof as \eqref{revise:02}, we have
\begin{align*}
\|u\|_{\mL_T^\infty}\le\|\widetilde u\|_{\mL_T^\infty}+\|u^{\lambda_b}\|_{\mL_T^\infty}\lesssim \ell_{b,\beta}^{\frac{\gamma\beta}{\alpha}}\|f\|_{L^\infty_T\bC^{-\beta}}.
\end{align*}
Following the proof of  \eqref{Hzm02z}, one sees that for any $\eps>0$ and $\theta\in[0,\alpha)$, and all $\lambda\ge0$,
\begin{align*}
\|u\|_{L^\infty_T\bC^{\alpha-\beta-\theta}}
\lesssim  (\lambda+1)^{-\frac{\theta}{\alpha}}\left(\|f\|_{L^\infty_T\bC^{-\beta}}+\ell_{b,\beta} \|u\|_{L^\infty_T\bC^{1+\eps}}\right),
\end{align*}
which, similar to the proof of \eqref{revise:03}, yields that
\begin{align}
    (1+\lambda)^{\frac{\theta}{\alpha}}\|u\|_{L^\infty_T\bC^{\alpha-\beta-\theta}}
&\le c\|f\|_{L^\infty_T\bC^{-\beta}} \ell_{b,\beta}^\gamma+\tfrac12\|u\|_{L^\infty_T\bC^{\alpha-\beta}}.\label{revise:04}
\end{align}
In particular, taking $\theta=0$, we get
\begin{align*}
\|u\|_{L^\infty_T\bC^{\alpha-\beta-\theta}}\le c\ell_{b,\beta}^\gamma\|f\|_{L^\infty_T\bC^{-\beta}}.
\end{align*}
Substituting this into \eqref{revise:04}, we obtain \eqref{Hzm01:00} and complete the proof.

\section{Well-posedness of SDEs with $\bC^{-\beta}$ drifts}\label{sec:SDE}

In this section, based on the PDE estimates established previously, we demonstrate the weak well-posedness of SDEs with distributional drift and present the proof of our main result, Theorem~\ref{thm:G-mart}. Let $\alpha \in (1,2)$ and $T>0$. Recall the SDE with multiplicative noise and distributional drift:
\begin{align}\label{eq:SDE}
\dif X_t = b(t,X_t) \dif t + \sigma(t,X_{t-}) \dif L_t^{(\alpha)}, \quad X_0 = x \in \mathbb{R}^d,
\end{align}
where $b\in L_T^\infty \bC^{-\beta}$ for some $\beta \in (0, \alpha - 1)$, and the diffusion coefficient $\sigma: \mathbb{R}_+ \times \mathbb{R}^d \to \mathbb{R}^d \otimes \mathbb{R}^d$ satisfies the condition {\bf (H$^\sigma$)}. Our goal is to construct a filtered probability space $(\Omega, \sF, (\sF_t)_{t\in[0,T]}, \mathbb{P})$ and a pair of c$\rm\grave{a}$dl$\rm\grave{a}$g processes $(X, L^{(\alpha)})$ such that:
\begin{itemize}
    \item $(L_t^{(\alpha)})_{t\in[0,T]}$ is a symmetric $\alpha$-stable process on this space, with L\'evy measure given by \eqref{CC02}, satisfying the non-degeneracy condition {\bf (ND)};

\item $(X_t)_{t\in[0,T]}$ is a weak solution to the SDE \eqref{eq:SDE} in the sense of Definition \ref{Def1};

\item uniqueness in law holds for such solutions.
\end{itemize}
Notice that, in the Definition \ref{Def1}, the drift term in \eqref{eq:SDE} is interpreted via a mollification procedure. Prescisely, 
$$
\int_0^t b(s,X_s) \dif t:=A_t^{b}:=\lim_{m\to \infty}  \int_0^t b_m(s,X_s)\dif s \quad \text{in $L^2(\Omega)$-sense},
$$
where
\begin{align*}
b_m(t, x) := b(t, \cdot) * \phi_m(x),\quad \phi_m(x) = m^d \phi(mx),\quad m \in \mN,
\end{align*}
with $\phi \in C_c^\infty(\mathbb{R}^d)$ being a smooth function with compact support and unit integral. 

In the following, we fix the mollifier family $\phi,\{\phi_m\}_{m > 0}$, and notice that all the results and arguments in this section do not depend on the choice of $\phi$. We denote the space of all continuous functions  from $[0,T]$ to $\mathbb{R}^d$ by $C([0,T])$, and let $\bD := D([0,T], \mathbb{R}^d)$ be the path space consisting of all c$\rm\grave{a}$dl$\rm\grave{a}$g functions from $[0,T]$ to $\mathbb{R}^d$. Throughout this section, $C([0,T])$ is endowed with the uniform metric, and $\bD$ is equipped with the Skorokhod topology under which it becomes a Polish space.

\subsection{The Krylov class}

In order to define $A_t^{b}$, we introduce the following Krylov's class.
\bd[Krylov's class]
Fix $T>0$. Let $X:=(X_t)_{t\in[0,T]}$ be a $\mR^d$-valued c$\rm\grave{a}$dl$\rm\grave{a}$g stochastic process on a filtered probability space $(\Omega,\sF,(\sF_t)_{t\in[0,T]},\mP)$. For given $ \theta>0$ and $ \beta>0$, we call that the  process $X$ is in the Krylov's class $\sK^\beta_\theta$ (denoted by $X\in \sK^\beta_\theta$) if for any $p\ge 2$, there is a constant $c>0$ such that for all $\delta>0$, $ f\in L_T^\infty C_b  (\mR^d)$, and stopping times $0\le \tau_0\le \tau_1\le T$ with $\tau_1-\tau_0\le \delta$,
\begin{align}\label{krylov}
    \left\|\int_{\tau_0}^{\tau_1}f(s,X_s)\dif s\right\|_{L^p(\Omega)}\lesssim_c \delta^{\frac{1+\theta}{\alpha}}\|f\|_{L^\infty_T \bC^{-\beta}}.
\end{align}
\ed

The following proposition tells us that if a process $X$ is in a Krylov's class, then the process $A^f_t$ defined by a limit is well-defined for all $f\in L^\infty_T\bC^{-\beta}$.  

\bp\label{Pr88}
Fix $T>0$. Suppose that $X:=(X_t)_{t\in[0,T]} \in \sK^\beta_\theta$ with parameters $\beta>0$ and $\theta>0$.
For any $\eps>0$ and  $f\in L^\infty_T\bC^{-\beta+\eps}$, let $f_m(t):=f(t,\cdot)*\phi_m$.
Then the limit 
\begin{align}\label{eq:SD01} 
A_\cdot^{f}:=\lim_{m\to\infty}\int^\cdot_0f_m(s,X_s)\dif s
\end{align}
exists in $L^p(\Omega;C([0,T]))$ for any $p\ge2$, and does not depend on the choice of $\phi$.
Furthermore, for any $\gamma \in (0, \frac{1+\theta}{\alpha} - \frac{1}{p})$, there are two constants $c_1,c_2>0$ such that
\begin{align}\label{S4:AA01}
\big\|[A^{f}_\cdot]_{\sC^{\gamma}}\big\|_{L^p(\Omega)}\lesssim_{c_1} \sup_{t\in[0,T]}\|A^{f}_t\|_{L^p(\Omega)}+\sup_{s\ne t\in[0,T]}\frac{\|A^{f}_t-A^{f}_s\|_{L^p(\Omega)}}{|t-s|^{(1+\theta)/\alpha}}\lesssim_{c_2} \|f\|_{L^\infty_T\bC^{-\beta}}.
\end{align}
\ep

\begin{proof}
In view of \eqref{0725:new00}, we have
$$
\|f-f_m\|_{L^\infty_T\bC^{-\beta}}\lesssim n^{ -1\wedge \eps}\|f\|_{L^\infty_T\bC^{-\beta+\eps}},
$$
which implies that $\{f_m\}_{n\in\mN}$ is a Cauchy sequence in $L_T^\infty\bC^{-\beta}$. Based on \eqref{krylov}, we have for all $0\le s<t\le T$ and $p \geq 2$,
\begin{align*} 
   \mE\left|A^{f_m-f_{m'}}_t-A^{f_m-f_{m'}}_s\right|^p= \mE\left|\int_s^t [f_m-f_{m'}](r,X_r)\dif r\right|^p \lesssim (t-s)^{\frac{(1+\theta)p}{\alpha}}\|f_m-f_{m'}\|_{L^\infty_T\bC^{-\beta}}^p.
\end{align*}
Since $p\geq 2 >\alpha/(1+\theta)$, by the Kolmogorov continuity theorem, we conclude that
\begin{align*}
   \left\|\left(\sup_{t\in[0,T]}|A^{f_m-f_{m'}}_t|\right)\right\|_{L^p(\Omega)}\lesssim \|f_m-f_{m'}\|_{L^\infty_T\bC^{-\beta}}\to 0,
\end{align*}
as $m,m'\to \infty$. Thus, $\{\int_0^.f_m(s,X_s)\dif s\}_{m\in\mN}$ has a limit $A^f_\cdot$ in $L^p(\Omega;C([0,T]))$. Furthermore, \eqref{S4:AA01} follows directly from \eqref{krylov} and Kolmogorov's continuity theorem. 
\end{proof}

Next we show a substitution formula for Young's integrals that will be used to show the uniqueness in law.

\bl\label{Pro44}
Fix $T>0$. Assume that $X$ belongs to $\sK^\beta_\theta$ with $\beta,\theta>0$, and has a.s. $q$-variation with any or some $q>2$. Let $g:\mR_+\times\mR^d\to\mR$ be a bounded function satisfying
\begin{align}\label{HF2}
|g(t,x)-g(s,y)|\lesssim_c {|t-s|^{\frac1{2}}}+|x-y|. 
\end{align}
Then for any function $f\in L^\infty_T\bC^{-\beta+\eps}$ with some $\eps>0$,  the integral $\int^t_0g(s, X_s)\dif A^{f}_s$ is well-defined as the Young integral, and
\begin{align}\label{AA0} 
\int^{\cdot}_0g(s, X_s)\dif A^{f}_s=A^{f\cdot g}_\cdot\ \ a.s.,
\end{align}
where $A^{f},A^{f\cdot g}$ is defined as in Proposition \ref{Pr88}.
\el
\begin{proof}
Choosing $p > \frac{\alpha}{\theta}$, one can choose $\gamma\in(\frac12,\frac{1+\theta}\alpha-\frac1p)$ and  
$q>2$ so that
$$
\gamma+\tfrac1q>1.
$$
By \eqref{S4:AA01}, there is a constant $c>0$ such that for all $f\in L^\infty_T\bC^{-\beta}$,
\begin{align}\label{HF5}
\sup_{t\in[0,T]}\|A^f_t\|_{L^1(\Omega)}+ \big\|[A^f_\cdot]_{\sC^{\gamma}}\big\|_{L^1(\Omega)}\lesssim_c \|f\|_{L^\infty_T\bC^{-\beta}}.
\end{align}
Below, we fix  a sample point $\omega$ such that
\begin{align}\label{HF3}
[A^f_\cdot(\omega)]_{\sC^\gamma}<\infty,\ \ [X_\cdot(\omega)]_{q,\var}<\infty.
\end{align}
For any $0\le s<t\le T$, we define
\begin{align*}
\Gamma_{s,t}(\omega):=g(s,X_s(\omega))(A^f_t(\omega)-A^f_s(\omega));
\end{align*}
and
for any $0\le s<u<t\le T$,
\begin{align*}
\delta \Gamma_{s,u,t} (\omega):=\Gamma_{s,t}(\omega)-\Gamma_{s,u}(\omega)-\Gamma_{u,t}(\omega)=(g(s,X_s(\omega))-g(u,X_u(\omega)))(A^f_t(\omega)-A^f_u(\omega)).
\end{align*}
In the sequel, for the simplicity of notation, we drop the dependence of $\omega$. By \eqref{HF2} and \eqref{HF3}, we have
 \begin{align*}
|\delta \Gamma_{s,u,t}|& \lesssim [A^f_\cdot]_{\sC^\gamma}|t-u|^{\gamma}\left(|u-s|^\frac12 +[X]_{q,\var;[s,u]}\right)\\
&\lesssim  [A^f_\cdot]_{\sC^\gamma}|t-u|^{\gamma}\left(T^{q/2-1}|u-s| +[X]^q_{q,\var;[s,u]}\right)^{1/q},
\end{align*}
which together with the facts that $[X]_{q,\var;[s,u]}^q$ is a control function and $\gamma+\frac1q>1$, by the sewing lemma (cf. \cite[Theorem 2.2]{FZ18}), derives that
 \begin{align*}
\Gamma_t:=\int_0^t g(r,X_r)\dif A^f_r:=\lim_{|\pi|\to0}\sum_{r,s\in\pi}\Gamma_{r,s}\quad \text{exists},
\end{align*}
where the limit is taken over any sequence of partitions $\pi$ of the interval $[0,t]$ with mesh size $|\pi| \to 0$. Consequently, we have
 \begin{align*}
\left|\int_s^t(g(r,X_r)-g(s,X_s))\dif A^f_r\right|=|\Gamma_t-\Gamma_s-\Gamma_{s,t}|\lesssim [A^f_\cdot]_{\sC^\gamma}|t-s|^{\gamma}\left(|t-s|^\frac12 +[X]_{q,\var;[s,t]}\right),
\end{align*}
which, by taking $s=0$, deduces that
 \begin{align*}  
 \left|\int_0^tg(r,X_r)\dif A^f_r\right|   \lesssim_T \|g\|_\infty  |A^f_t|+  [A^f_\cdot]_{\sC^\gamma}\left(1 +  [X]_{q,\var}\right).
\end{align*}
From this, for $f_m(t):=f(t,\cdot)*\phi_m$, one sees that any $R>0$,
\begin{align*}
\mE \left| \int_0^t g(r,X_r)\dif (A^{f_m}_r-A^{f}_r)\right| \1_{\{ [X]_{q,\var} \leq R \}} \overset{\eqref{HF5}}{\lesssim_R } \| f_m - f\|_{L_T^\infty \bC^{-\beta}}\to 0, \quad \text{as $m \to \infty$}.
\end{align*}
Thus, for any $\eps>0$, observing that 
\begin{align*}
\mP \left(\left| \int_0^t g(r,X_r)\dif (A^{f_m}_r-A^{f}_r)\right|  >\eps \right) & \leq \mP \left(\left| \int_0^t g(r,X_r)\dif (A^{f_m}_r-A^{f}_r)\right|  >\eps , [X]_{q,\var} \leq R\right)\\
&\qquad  +  \mP ([X]_{q,\var}> R), 
\end{align*}
we get that
$$
\left| \int_0^t g(r,X_r)\dif (A^{f_m}_r-A^{f}_r)\right| \xrightarrow{~\mP~} 0,\quad \text{as $m \to \infty$}.
$$
On the other hand, we have that
\begin{align*}
 \mE|A^{g\cdot f_m}_t-A^{g\cdot f}_t|
= \mE|A^{g\cdot (f_m-f)}_t| \overset{\eqref{krylov}} {\lesssim} \|g\cdot (f-f_m)\|_{L^\infty_T\bC^{-\beta}} \to 0, \quad\text{as $m\to\infty$}.
\end{align*}
Hence, noting that
\begin{align*}
\int_0^t g(r,X_r)\dif A^{f_m}_r=\int_0^t (g\cdot f_m)(r,X_r) \dif r=A^{g\cdot f_m}_t\quad a.s.,
\end{align*}
by taking $m\to \infty$, we 
deduce \eqref{AA0}.
\end{proof}

Now, we give the following result which will play a crucial role in the proof of uniqueness.

\bc\label{geito}
Let $\alpha\in(1,2)$, $\beta\in (0,\frac{\alpha-1}{2})$ (resp. $\beta\in  [\frac{\alpha-1}{2},\alpha-1)$). Assume that $b \in L_T^\infty \bC^{-\beta}$ (resp. $b,\div b \in L_T^\infty \bC^{-\beta}$), $f \in L_T^\infty \bC^{-\beta}$. Suppose that $(X,L^{(\alpha)})$ is a weak solution of SDE \eqref{eq:SDE} in the sense of Definition \ref{Def1}, and $X\in \sK^{ \beta+\eps}_\theta$ with some $\eps>$ and $\theta>0$. Let $u$ be the unique strong (resp. paraproduct) solution, in the sense of Definition \ref{AA05} (resp. Definition \ref{AA05:00}), for the following backward PDE:
\begin{align}\label{0728:PDE}
    \partial_t u+\sL^{(\alpha)}_\sigma u+b\cdot\nabla u=f,\quad u(T)=0.
\end{align}
Then, for all $t\in[0,T]$,
$$
u(t,X_t)=u(0,X_0)+A_t^{f}+\int_0^t\int_{\mR^d}  \Big(
u(s,X_{s-}+\sigma(s,X_{s-})z)-u(s,X_{s-})
\Big)\widetilde{N}(\dif s,\dif z),
$$
where $ A_t^{f}$ is defined as in Proposition \ref{Pr88}.
\ec
\br
Based on Theorems \ref{BB01} and \ref{BB01:00}, the unique strong solution (paraproduct) $u$ exists. 
\er
 
\begin{proof}[Proof of Corollary \ref{geito}]
Let $b_m(t):=b(t,\cdot)*\phi_m$ and $f_m(t):=f(t,\cdot)*\phi_m$. Since $b_n,f_n \in L^\infty_T C^\infty_b(\mR^d)$, it is well known that there is a unique smooth solution $u_m \in L^\infty_T\bC^2$ to PDE \eqref{0728:PDE} with $(b,f)=(b_m,f_m)$ (see \cite{HWW20} for example). Moreover, by definitions, it is easy to check that
\begin{align*}
|u_m(t,x)-u_m (s,y)| 
& \le |t-s|\left(\|\sL^{(\alpha)}_\sigma u_m\|_{\mL^\infty_T}+\|b_m\cdot\nabla u_m\|_{\mL^\infty_T}+\|f_m\|_{\mL^\infty_T}\right)
+ 2\|\nabla u_m\|_{\mL^\infty_T} |x-y|
\\
& \lesssim_m |t-s|  + |x-y|,
\end{align*}
and then 
\begin{align}
\|\nabla \big(u_m(t)-  u_m (s)\big)\|_{\infty} & \overset{\eqref{InIn}}{\lesssim} \|\nabla \big(u_m(t)-  u_m (s)\big)\|_{\bC^{-1}}^{1/2}\|\nabla (u_m(t)-  u_m (s))\|_{\bC^{1}}^{1/2} \nonumber\\
&\,\,\, \lesssim \| u_m(t)- u_m (s)\|_{\infty}^{1/2}\|u_m\|_{L_T^\infty\bC^{2}}^{1/2}
 \lesssim_m |t-s|^{1/2},
\label{eq:WK00}
\end{align}
where we used Bernstein's inequality in the second inequality. It follows from \eqref{S4:AA01} that 
\begin{align*}
    t\to A^b_t:=\underset{m \to\infty}{\lim-L^p} \left( \int_0^t b_m (s,X_s)\dif s \right) \in \sC^{\frac{1+\theta}{\alpha}-\frac{1}{p}} ,\quad \text{for any $p>2$},
\end{align*}
which implies that
$$
\left[A_\cdot^b \right]_{q,\var}< \infty,\quad a.s., \quad \text{for any $q> \frac{\alpha}{1+\theta}$ }.
$$
From Lemma \ref{lem:B2}, one sees that
\begin{align*}
    \left[\int_0^\cdot \sigma(s,X_{s-})\dif L^{(\alpha)}_s\right]_{q,\var}<\infty,\quad a.s., \quad \text{for any $q>2$}.
\end{align*}
Thus, we get
\begin{align*} 
    \left[X_\cdot \right]_{q,\var}<\infty, \quad a.s., \quad \text{for any $q>2$},
\end{align*}
which together with \eqref{eq:WK00}, in view of Lemma \ref{Pro44}, derives that
\begin{align*}
    \int_0^t \nabla u_m (s,X_s)\dif A^b_s= A^{b\cdot \nabla u_m}_t,
\end{align*}
where the left hand is in the sense of Young's integral, and the right hand is defined by \eqref{eq:SD01}. By \eqref{S4:AA01} again, we also have
$$
\left[A^{b\cdot \nabla u_m}_\cdot  \right]_{q,\var}< \infty,\quad a.s., \quad \text{for any $q>\tfrac\alpha{1+\theta}$}.
$$ 
Applying It\^o's formula to $t\to u_m (t,X_t)$, we have
\begin{align}\label{0728:new02}
\begin{split}
u_m (t,X_t)=& u_m (0,X_0)+\int_0^t (\p_s +\sL^{(\alpha)}_\sigma)u_m(s,X_s)\dif s+A^{b\cdot\nabla u_m}_t\\
&+  \int_0^t\int_{\mR^d}\Big(u_m(s,X_{s-} + \sigma(s,X_{s-})z)-u_m(s,X_{s-}) \Big)\widetilde{N}(\dif s,\dif z),
\end{split}
\end{align}
where
\begin{align*}
 \int_0^t \( \p_s +\sL^{(\alpha)}_\sigma \)u_m(s,X_s)\dif s+A^{b\cdot\nabla u_m}_t=A^{(\p_s +\sL^{(\alpha)}_\sigma+b\cdot\nabla)u_m}_t.   
\end{align*}

Now we estimate every term in \eqref{0728:new02} in turn. Notice that
\begin{align}
& \qquad \|(b-b_m)\cdot\nabla u_m \|_{L^\infty_T\bC^{-\beta}} \nonumber\\
& \overset{\eqref{eq:YY01}}{\lesssim} 
\begin{cases}
    \|b-b_m\|_{L^\infty_T\bC^{-\beta}}\|u_m\|_{L^\infty_T\bC^{\alpha-\beta}},\qquad \qquad\qquad\qquad\quad\quad\beta<\tfrac{\alpha-1}{2},\\
    (\|b-b_m\|_{L^\infty_T\bC^{-\beta}}+\|\div b-\div b_m\|_{L^\infty_T\bC^{-\beta}})\|u_n\|_{L^\infty_T\bC^{\alpha-\beta}},\quad \beta\ge\tfrac{\alpha-1}{2}.
\end{cases}\label{eq:WK90}
\end{align}

\noindent
(i) By \eqref{eq:WK90}, following the proof of \eqref{0728:new00}, it is easy to check that
\begin{align*}
\mE |u(t,X_t) - u_m(t,X_t)|^2 \leq \|u_m-u\|_{\mL^\infty_T}^2
\lesssim \|u_m-u\|_{L^\infty_T\bC^{\alpha-\beta}}^2, \quad \text{as $m \to \infty$}.
\end{align*}

\noindent
(ii) Note that
\begin{align*}
    \|(\p_s +\sL^{(\alpha)}_\sigma+b\cdot\nabla)u_m-f\|_{L^\infty_T\bC^{-\beta}}
\overset{\eqref{0728:PDE}}{\lesssim } 
\|(b-b_m)\cdot\nabla u_m \|_{L^\infty_T\bC^{-\beta}} + \|f-f_m\|_{L^\infty_T\bC^{-\beta}},
\end{align*}
which tends to $0$ as $n\to\infty$ since \eqref{eq:WK90}. Hence, by Proposition \ref{Pr88}, we have 
\begin{align*}
    \lim_{m\to\infty}
    \left( \int_0^t (\p_s +\sL^{(\alpha)}_\sigma)u_m(s,X_s)\dif s+A^{b\cdot\nabla u_m}_t \right)
    =\lim_{m\to\infty}A^{(\p_s +\sL^{(\alpha)}_\sigma+b\cdot\nabla)u_m}_t=A^f_t,\quad 
\end{align*}
in $L^p(\Omega)$ for any $p\ge 2$.

\noindent
(iii) Noting that $\alpha-\beta>\frac{\alpha}{2}$, 
by the isometry of stochastic integral, we have that for small enough $\delta>0$,
\begin{align*}
    &\quad \mE\left|\int_0^t\int_{\mR^d}\Big( [u_m-u](s,\sigma(s-,X_{s-} + \sigma(s,X_{s-})z)-[u_m-u](s, X_{s-})\Big) \widetilde{N}(\dif s,\dif z)\right|^2\\
    &= \mE\int_0^t\int_{\mR^d}\Big|[u_m-u](s,X_{s-} + \sigma(s,X_{s-})z)-[u_m-u](s, X_{s-}) \Big|^2\nu^{(\alpha)}(\dif z)\dif s\\
    &\lesssim \|u_m-u\|_{L^\infty_T\bC^{\frac{\alpha}{2}+\delta}}^2\int_{\mR^d}(|z|^{\alpha+2\delta}\wedge 1)\nu^{(\alpha)}(\dif z)
    \overset{\eqref{eq:AE02}}{\lesssim}\|u_m-u\|_{L^\infty_T\bC^{\alpha-\beta}}^2\to 0,\quad \text{as $m\to \infty$}.
\end{align*}

Combining the calculations above and taking $m\to\infty$ on both sides of \eqref{0728:new02}, we complete the proof.
\end{proof}

\subsection{Proof of Theorem \ref{thm:G-mart}}

With the above preparations in place, we are now ready to establish the existence and uniqueness of the weak solution to the SDE \eqref{eq:SDE}. First, we introduce the following mollified SDE:
\begin{align}\label{eq:SDEap}
\dif X_t^m = b_m(t,X^m_t) \dif t + \sigma (t,X^m_{t-}) \dif L_t^{(\alpha)},\ \ X_0^m\overset{(d)}{=}\mu\in\cP(\mR^d),
\end{align}
where $b_m(t,x):=(\phi_m * b)(t,\cdot)(x)$, $m \in \mN$.
It is well-known (for example, see \cite{CZZ21}, Theorem 1.1) that there is a unique strong solution $X^m$ for SDE \eqref{eq:SDEap} on $(\Omega,\sF, (\sF)_{t \geq 0}, \mP)$. 

\bl\label{lem:b}
Let $T>0$, $\alpha\in(1,2)$ and $\beta \in (0, \frac{\alpha-1}{2})$ (resp. $\beta\in  [\frac{\alpha-1}{2},\alpha-1)$). Assume that $b \in L_T^\infty \bC^{-\beta}$ (resp. $b,\div b \in L_T^\infty \bC^{-\beta}$). Then $X^m\in \sK^\beta_{\alpha-\beta-1}$ with a constant $c>0$ independent of $m$.
\el
\begin{proof}
For simplicity of notation, we drop the superscript $(\alpha)$ in $\nu^{(\alpha)}$. Fix $ f \in L_T^\infty C_b^\infty (\mR^d)$. Recalling PDE \eqref{PDE}, let $u_m^\lambda$ be a classical solution of the backward equation
\begin{align}\label{eq:WW00}
\p_t u_m^\lambda +\(\sL^{(\alpha)}_{\sigma}-\lambda \) u_m^\lambda+b_m\cdot\nabla u_m^\lambda + f = 0,\ \ u_m^\lambda(T)=0.
\end{align}
By It\^o's formula (cf. \cite{IW89}, Theorem 5.1 of Chapter II), we have that for any stopping time $\widetilde \tau \leq T$,
\begin{align*}
u_m^\lambda(\widetilde \tau,X_{\widetilde \tau}^m) &- u_m^\lambda(0,x) = 
 \int_0^{{\widetilde \tau}} (\p_r u_m) (r,X_{r}^m) \dif r 
+ \int_0^{{\widetilde \tau}} b_m(r, X_{r}^m) \cdot  \nabla  u_m^\lambda (r,X^m_{r}) \dif r\\
 &+  \int_0^{\widetilde \tau} \int_{\mR^d} \( u_m (r,X_{r-}^m+\sigma(r,X_{r-}^m)z) - u_m^\lambda (r,X_{r-}^m) \) \widetilde N(\dif r, \dif z)\\
&+ \int_0^{\widetilde \tau} \int_{\mR^d}  \( u_m^\lambda (r,X_r^m+\sigma(r,X_{r-}^m)z) - u_m^\lambda (r,X_{r}^m) \\
&\qquad\qquad\qquad\qquad\qquad -\sigma(r,X_{r-}^m)z \cdot \nabla u_m^\lambda (r,X_{r}^m) \) \nu( \dif z) \dif r,
\end{align*}
which together with \eqref{eq:WW00} derives that 
\begin{align}\label{eq:f}
 \int_{\tau}^{\tau+\delta} f (r,X^m_r) \dif r =&\ \   u_m^\lambda(\tau,X^m_\tau) - u_m^\lambda(\tau+\delta,X^m_{\tau+\delta} )   +  \lambda \int_{\tau}^{\tau+\delta} u_m^\lambda(r,X^m_r ) \dif r +M_{\tau,\tau+\delta}
\end{align}
where we substituted $\tau$ and $\tau+\delta$ for $\widetilde \tau$, and 
$$
M_{\tau,\tau+\delta}:= \int_{\tau}^{\tau+\delta}  \int_{\mR^d} \( u_m^\lambda (r,X^m_{r-}+\sigma(r,X^m_{r-})z) - u_m^\lambda (r,X^m_{r-}) \) \widetilde N(\dif r, \dif z).
$$ 
Notice that, by the condition $\bf (H^\sigma)$, we have that
\begin{align*}
| u_m^\lambda (r,X^m_{r-}+\sigma(r,X^m_{r-})z) - u_m^\lambda (r,X^m_{r-}) |  \leq (2\|u^\lambda_m \|_{\mL_T^\infty}) \wedge (c_0 \|\nabla u^\lambda_m\|_{\mL_T^\infty}|z|) ,
\end{align*}
which together with Kunita's inequality (cf. \cite{Ku04}, Theorem 2.11) yields that for any $p\geq 2$,
\begin{align}\label{eq:WH00}
\mE |M_{\tau,\tau+\delta}|^p & \lesssim  \mE\[\( \int_\tau^{\tau+\delta} \int_{\mR^d} | u_m^\lambda (r,X^m_{r-}+\sigma(r,X^m_{r-})z) - u_m^\lambda (r,X^m_{r-}) |^2 \nu(\dif z) \dif r\)^{p/2}\] \nonumber\\
&\quad +  \mE  \int_\tau^{\tau+\delta} \int_{\mR^d} | u_m^\lambda (r,X^m_{r-}+\sigma(r,X^m_{r-})z) - u_m^\lambda (r,X^m_{r-}) |^p \nu (\dif z) \dif r \nonumber\\
 &\lesssim \delta^{p/2}\left(\int_{\mR^d}\left(\|u_m^\lambda\|_{\mL_T^\infty}\wedge(\|\nabla  u_m^\lambda\|_{\mL_T^\infty})|z|\right)^2\nu (\dif z)\right)^{p/2}\nonumber\\
 &\quad+ \delta\int_{\mR^d}\left(\|u_m^\lambda\|_{\mL_T^\infty}\wedge(\|\nabla  u_m^\lambda\|_{\mL_T^\infty})|z|\right)^p\nu (\dif z).     
\end{align}
Observe that by \eqref{InIn} and Bernstein's inequality, we have that
\begin{align*}
\|\nabla u_m^\lambda\|_{\mL_T^\infty} \lesssim \|\nabla u_m^\lambda \|^{1/2}_{L_T^\infty\bC^{-(\alpha-\beta)+1}} \|\nabla u_m^\lambda\|^{1/2}_{L_T^\infty\bC^{(\alpha-\beta)-1}} 
 \lesssim \|  u_m^\lambda \|^{1/2}_{L_T^\infty\bC^{2-(\alpha-\beta) }} \|  u_m^\lambda\|^{1/2}_{L_T^\infty\bC^{\alpha-\beta }};
\end{align*}
and by \eqref{InIn:00},
$$
\|u_m^\lambda\|_{L_T^\infty\bB_{\infty,1}^0} \lesssim \|u_m^\lambda\|_{L_T^\infty\bC^{-\beta/2}}^{1/2}\|u_m^\lambda\|_{L_T^\infty\bC^{\beta/2}}^{1/2}.
$$
Hence, by Theorems \ref{BB01} and \ref{BB01:00}, one sees that  for any $\lambda \ge0$,
\begin{align}\label{eq:WH01}
\| \nabla u_m^\lambda \|_{\mL_T^\infty } \lesssim \lambda^{-\frac{(\alpha-\beta)-1}{\alpha}} \|f\|_{L_T^\infty\bC^{-\beta}}\ \ \text{and}\ \  \|u_m^\lambda \|_{\mL_T^\infty}\leq \|u_m^\lambda \|_{L_T^\infty\bB_{\infty,1}^0}  \lesssim \lambda^{-\frac{\alpha-\beta}{\alpha}} \|f\|_{L_T^\infty\bC^{-\beta}},
\end{align}
where the implicit constants depend on a polynomial of positive degree in $\ell_{b,\beta}$ (see the definitions in Theorems \ref{BB01} and \ref{BB01:00}). Consequently, from \eqref{eq:f}, by \eqref{eq:WH00}, we obtain that for any $\lambda > \lambda_0$,
\begin{align*}
\mE \left|  \int_{\tau}^{\tau+\delta} f (r,X^m_r) \dif r \right|^p 
& \lesssim (2+\lambda \delta)^p\|u_m^\lambda \|_{\mL_T^\infty}^p + \delta\int_{\mR^d}\left(\|u_m^\lambda\|_{\mL_T^\infty}\wedge(|\nabla  u_m^\lambda\|_{\mL_T^\infty}|z|)\right)^p\nu (\dif z)\\
& \quad+ \delta^{p/2}\left(\int_{\mR^d}\left(\|u_m^\lambda\|_{\mL_T^\infty}\wedge(\|\nabla  u_m^\lambda\|_{\mL_T^\infty}|z|)\right)^2\nu (\dif z)\right)^{p/2}\\
&\overset{\eqref{eq:WH01}}{\lesssim} \Bigg[ \lambda^{-\frac{\alpha-\beta}{\alpha}p}(1+(\lambda\delta)^p)+\delta\lambda^{-\frac{\alpha-\beta}{\alpha}p}\int_{\mR^d}\left(1\wedge(\lambda^{\frac{p}\alpha}|z|^p)\right)\nu (\dif z)\\
& \quad+\delta^{p/2}\lambda^{-\frac{\alpha-\beta}{\alpha}p}\left(\int_{\mR^d}\left(1\wedge(\lambda^{\frac2\alpha}|z|^2)\right)\nu (\dif z)\right)^{p/2} \Bigg] \|f\|_{L_T^\infty\bC^{-\beta}}^p.
\end{align*}
Observe that $\|b_m\|_{L^\infty_T\bC^{-\beta}} \lesssim \|b\|_{L^\infty_T\bC^{-\beta}}$. Therefore, by \eqref{20230418}, we have
\begin{align*}
\mE \left|  \int_{\tau}^{\tau+\delta} f (r,X^m_r) \dif r \right|^p 
&
 \lesssim \lambda^{-\frac{p(\alpha-\beta)}{\alpha}} (1+(\lambda\delta)^p+(\lambda\delta)^{p/2}+ \lambda\delta)  \|f\|_{L_T^\infty\bC^{-\beta}}^p.
\end{align*}
Taking $\lambda  = \lambda_0 \delta^{-1}$, we get the desired result.
\end{proof}

Now we give the tightness of $\{X^m,m\geq 1\}$. 

\bt\label{thm:t}
The sequence $\{X^m \}_{m\in\mN}$  in $\bD$ is tight.
\et

\begin{proof}
Fix $T>0$. Let $\tau$ be a bounded stopping time  satisfying $0\leq \tau \leq \tau  + \delta\leq T$ with  $\delta\in(0,1)$. By SDE \eqref{eq:SDEap}, we have
\begin{align*}
X_{\tau+\delta}^m - X_{\tau }^m
=   \int_{\tau}^{\tau+\delta} \sigma(r,X_{r-}^m) \dif L_r^{(\alpha)}  
 +   \int_{\tau}^{\tau+\delta} b_m(r,X_r^m) \dif r. 
\end{align*}
For $p\in [1,\alpha)$, by Chebyshev's inequality, Lemma \ref{lem:BH}, and Lemma  \ref{lem:b}, we have that  for each $R>0$, 
\begin{align}\label{eq:WW02}
\mP(|X_{\tau+\delta}^m - X_{\tau}^m| \geq R)   \leq  R^{-p} \mE   |X_{\tau+\delta}^m - X_{\tau}^m|^p 
  \lesssim \delta^{ p/\alpha} (\|\sigma\|_{\mL_T^\infty}^p + \|b_m \|_{L_T^\infty \bC^{-\beta}}^p),
\end{align}
where the implicit constant in the last inequality is independent of $m,\tau$, and $\delta$. Furthermore,
$$
\lim_{\delta \downarrow 0}\sup_m \sup_{\tau\leq T}\sup_{h \in [0,\delta]} \mP(|X_{\tau+h}^m - X_{\tau}^m| \geq R) =0,
$$
provided by $\|b_m\|_{L_T^\infty\bC^{-\beta}} \leq \|b\|_{L_T^\infty\bC^{-\beta}} $. On the other hand, since all $X^m\in \sK^\beta_{\alpha-\beta-1}$ with a common constant independent of $m$, by \eqref{krylov} and Chebyshev's inequality again, one sees that for any $t>0$,
\begin{align}\label{eq:WW04}
\lim_{R \to \infty}\sup_{m} \mP\left( \left|\int_0^t b_m(r,X_r^m) \dif r\right|\geq R\right) =0.
\end{align}
Moreover, by Burkholder-Davis-Gundy's inequality for jump processes (see \cite{SZ15}, Lemma 2.3, or \cite{No75}, Theorem 1) and \eqref{eq:AE02}, there is a constant $c>0$ only depending on $T,\|\sigma\|_{\mL_T^\infty},\alpha,\nu^{(\alpha)}$ such that
\begin{align*}
\mE \left(\sup_{t \in [0,T]} \left|\int_0^t \sigma(r,X_{r-}^m) \dif L_r^{(\alpha)}\right| \right)
& \lesssim 
\left[ \mE \left( \int_0^T \int_{|z|\leq 1}  \left| \sigma(r,X_{r-}^m)z\right|^2 \nu^{(\alpha)}(\dif z) \dif r \right)\right]^{1/2} \\
& \qquad
+ \mE \left( \int_0^T \int_{|z|>1}  \left| \sigma(r,X_{r-}^m)z\right| \nu^{(\alpha)}(\dif z) \dif r\right)   \leq c ,
\end{align*}
which together with \eqref{eq:WW04}  yields that for any $t \in T$,
\begin{align}\label{eq:WW03}
\lim_{R\to \infty} \sup_m \mP \(|X_t^m| \geq R\) =0,
\end{align}
that is $\{ X_t^m\}_{m \geq 1}$ is tight for every $t\in [0,T]$. 
Hence, by \eqref{eq:WW02}, \eqref{eq:WW03} and Aldous's criterion for tightness (cf. \cite{Ka3rd}, Theorem 23.8 and Lemma 23.12), we conclude the proof. 
\end{proof}

Now we are in a position to give

\begin{proof}[Proof of Theorem \ref{thm:G-mart}]
{\bf (Existence)} For the simplicity of notation, we drop the superscript $(\alpha)$. By Theorem \ref{thm:t}, one sees that the sequence $\{(X^m,L)\}_{m\in\mN}$  in $\bD\times \bD$ is tight. Hence,  by Prokhorov's theorem (cf. \cite{Ka3rd}, Theorem 23.2), there is a subsequence, still denoted by $\{m\}$ for simplicity, and a probability measure $\bQ$ such that $ \mP\circ(X^m,L)^{-1}$ converges to $\bQ$ weakly in $\cP(\bD)$.
Then by Skorokhod's representation theorem (cf. \cite{Ka3rd}, Theorem 5.31), there is a probability space $  (\widetilde{\Omega},\widetilde{\sF},\widetilde{\mP})$ and $\bD\times\bD$-valued random variables $(\widetilde{X}^m,\widetilde{L}^m)$ and $(\widetilde{X},\widetilde{L})$ thereon such that 
$$
\widetilde{\mP}\circ(\widetilde{X}^m,\widetilde{L}^m)^{-1}=  {\mP}\circ({X}^m,{L})^{-1}, \qquad \widetilde{\mP}\circ(\widetilde{X} ,\widetilde{L} )^{-1}=\bQ,
$$ 
and
\begin{align}\label{0728:new06}
(\widetilde{X}^m,\widetilde{L}^m)\rightarrow(\widetilde{X},\widetilde{L})\quad \text{in $\bD\times \bD$},\ \ \widetilde{\mP}-\text{a.s.}
\end{align}
Defining $\widetilde \sF^m_t := \sigma \{ \widetilde X^m_s, \widetilde L^m_s; s \leq t\}$ and obeserving that 
$$
\mP [L_t - L_s \in \cdot \mid \sF_s] = \mP [L_t - L_s \in \cdot ],
$$
we have
$$
\widetilde \mP [\widetilde L^m_t - \widetilde L^m_s \in \cdot \mid \widetilde \sF^m_s] = \widetilde \mP [\widetilde L^m_t - \widetilde L^m_s \in \cdot ],
$$
which implies that $\widetilde L^m$ is an $(\widetilde \sF^m_t)$-adapted $\alpha$-stable L\'evy process. Hence, it is easy to see that
\begin{align}\label{0728:SDE}
\dif \widetilde{X}_t^m=b_m(t,\widetilde{X}^m_t)\dif t+\sigma(t,\widetilde{X}^m_{t-})\dif \widetilde{L}^m_t,
\end{align}
and $\widetilde L$ is an $\alpha$-stable process with respect to the filtration $\widetilde \sF_t:=\sigma(\widetilde X_s, \widetilde L_s; s\le t)$. In the following, we show that $(\widetilde X,\widetilde L)$ with $  (\widetilde{\Omega},\widetilde{\sF},(\widetilde{\sF}_t),\widetilde{\mP})$ is a weak solution to SDE \eqref{eq:SDE} in the sense of Definition \ref{Def1}. To this end, we investigate the convergence of the two terms on the right-hand side of \eqref{0728:SDE}.

\medskip
First, Let $\rho_n(\cdot):=n \rho(n\cdot)$ be the usual mollifiers with $\rho \in C_c^\infty(\mR) $, and $\sigma_n(t,x):=[\sigma(\cdot,x)1_{[0,T]}(\cdot)*\rho_n](t)$, $n \in \mN$. Note that
\begin{align*}
 \left|\int_0^t \sigma(s,\widetilde{X}^m_{s-})\dif \widetilde{L}^m_s - \int_0^t \sigma(s,\widetilde{X}_{s-})\dif \widetilde{L}_s\right| \leq J^{n,m}_1(t) + J^{n,m}2(t)+J^{n}_3(t),
\end{align*}
where 
$$
J^{n,m}_1(t):= \left|\int_0^t \sigma_n(s,\widetilde{X}^m_{s-})\dif \widetilde{L}^m_s - \int_0^t \sigma_n(s,\widetilde{X}_{s-})\dif \widetilde{L}_s\right| ,
$$
and
$$
J^{n,m}_2(t):= \left|\int_0^t (\sigma-\sigma_n)(s,\widetilde{X}^m_{s-})\dif \widetilde{L}^m_s\right| ,
$$
and
$$
J^{n}_3(t):= \left|\int_0^t (\sigma-\sigma_n)(s,\widetilde{X}_{s-})\dif \widetilde{L}_s\right| .
$$
Next, we estimate $J^{n,m}_1(t),J^{n,m}_2(t),J^{n}_3(t)$ in turn. For $J^{n,m}_1(t)$, notice that for each $n \in \mN$,
\begin{align*}
\sigma_n(\cdot,\widetilde{X}^m_\cdot)\to \sigma_n(\cdot,\widetilde{X}_\cdot)\ \ \text{in $\bD$},\ \ \widetilde{\mP}-{\rm a.s.}, \quad \text{as $m \to \infty$}.
\end{align*}
Since $\widetilde{L}^m$ and $\widetilde{L}$ are $(\widetilde{\sF}^m_t)$- and $(\widetilde{\sF}_t)$- semimartingales, respectively, by \cite[Theorem 2.2]{KP91}, we have that for any $n\in\mN$,
$$
    \int_0^\cdot \sigma_n(s,\widetilde{X}^m_{s-})\dif \widetilde{L}^m_s \to \int_0^\cdot \sigma_n(s,\widetilde{X}_{s-})\dif \widetilde{L}_s\ \ \text{in $\bD$},\ \ \widetilde{\mP}-{\rm a.s.}, \quad \text{as $m\to \infty$}.
$$
Then for each $n\in\mN$,
\begin{align}\label{0805:02}
J_1^{n,m}(\cdot) \to 0 \ \ \text{in $\bD$},\ \ \widetilde{\mP}-{\rm a.s.}, \quad \text{as $m\to \infty$}.
\end{align}
For $J^{n,m}_2(t)$, based on Burkholder-Davis-Gundy's inequality for jump processes (see \cite{SZ15}, Lemma 2.3, or \cite{No75}, Theorem 1), we have
\begin{align}
\widetilde\mE |J^{n,m}_2(t)| &\, \lesssim  \left[ \widetilde\mE \left( \int_0^t |(\sigma_n-\sigma)(s,\widetilde{X}^m_{s})|^2\dif s\int_{|z|\le1}|z|^2\nu^{(\alpha)}(\dif z)\right) \right]^{1/2}\nonumber\\
    &\qquad +\widetilde\mE\left[\int_0^t |(\sigma_n -\sigma)(s,\widetilde{X}^m_{s})|\dif s\int_{|z|>1}|z|\nu^{(\alpha)}(\dif z)\right]\nonumber\\
    & \overset{\eqref{eq:AE02}}{\lesssim} \left[\widetilde\mE\int_0^t |(\sigma_n -\sigma)(s,\widetilde{X}^m_{s})|^2\dif s\right]^{1/2}\nonumber\\
    & \lesssim  \left[\widetilde\mE\int_0^t |(\sigma_n-\sigma) (s,\widetilde{X}_{s})|^2\dif s\right]^{1/2} +\left[\widetilde\mE\int_0^t \left( |\widetilde{X}_{s}-\widetilde{X}^m_{s}|^2 \wedge 1 \right)\dif s\right]^{1/2},\label{0805:01}
\end{align}
where we also used the Jensen inequality in the second inequality, and the condition ($\bH^\sigma$) in the last inequality. For $J^{n}_3(t)$, similar to $J^{n,m}_2(t)$, it is easy to check that
\begin{align}\label{eq:XD00}
\widetilde\mE| J^{n}_3(t) | \lesssim \left[\widetilde\mE\int_0^t |(\sigma_n-\sigma) (s,\widetilde{X}_{s})|^2\dif s\right]^{1/2},
\end{align}
Thus, by \eqref{0805:01} and \eqref{eq:XD00}, one sees that for each $t$,
\begin{align*}
\widetilde{\mE}J_2^{n,m}(t)+\widetilde{\mE}J_3^{n}(t)  \to 0, \quad \text{as $n,m \to 0$},
\end{align*}
provided by the dominated convergence theorem and \eqref{0728:new06}. Finally, from this, by \eqref{0805:02}, we have that for Lebesgue everywhere $t\in[0,T]$,
\begin{align*}
    \int_0^t \sigma(s,\widetilde{X}^m_{s-})\dif \widetilde{L}^m_s \overset{\widetilde{\mP}}{\to} \int_0^t \sigma(s,\widetilde{X}_{s-})\dif \widetilde{L}_s,\quad \text{as $m\to \infty$}.
\end{align*}
and we only need to consider the convergence of the drift terms. 

\medskip
On the other hand, observing that $b\in L^\infty_T\bC^{-\beta}\subset L^\infty_T\bC^{-\widetilde{\beta}}$ for any $\widetilde \beta\in[\beta,\alpha-1)$, by Lemma \ref{lem:b}, it is easy to check that $\widetilde{X}^m\in \bigcap_{\widetilde \beta \in [\beta, \alpha-1)} \sK^{\widetilde \beta}_{\alpha-\widetilde \beta-1}$ with a common constant $c>0$ independent of $m$. Thus, by \eqref{0728:new06} and Fatou's lemma, we have
$ \widetilde{X}  \in \bigcap_{\widetilde \beta \in [\beta, \alpha-1)} \sK^{\widetilde \beta}_{\alpha-\widetilde \beta-1}$. Then in view of $b\in L^\infty_T\bC^{-\widetilde{\beta}+\eps}$ and $\widetilde X,\widetilde  X^m \in \sK^{\widetilde \beta}_{\alpha-\widetilde \beta-1}$ with $\eps=\widetilde{\beta}-\beta\in (0,\alpha-1-\beta)$, by Proposition \ref{Pr88}, for any $m\in\mN$,  
\begin{align}\label{eq:SD00}
A^{b,\widetilde{X}^m}_\cdot:=\underset{n \to\infty}{\lim-L^p}\int_0^\cdot b_n(s,\widetilde X^m_s)\dif s\quad \text{and}\quad A^{b,\widetilde{X}}_\cdot:=\underset{n \to\infty}{\lim-L^p}\int_0^\cdot b_n (s,\widetilde{X}_s)\dif s
\end{align}
exist in $L^p(\Omega;C([0,T]))$ for any $p\geq 2$, and satisfy \eqref{S4:AA01}. Next, we are going to show that for any $t \in (0,T]$, 
\begin{align}\label{0805:00}
      A^{b_m,\widetilde{X}^m}_t :=\int_0^t b_m(s,\widetilde X^m_s)\dif s \xrightarrow[]{~ L^1(\widetilde \Omega,\widetilde \mP)~}  A^{b,\widetilde X}_t,\ \ \text{as $m\to \infty$}.
\end{align}
Indeed, for any $n\in\mN$,
\begin{align*}
    | A^{b_m,\widetilde X^m}_t-A^{b,\widetilde X}_t|&\le   | A^{b_m,\widetilde X^m}_t-A^{b_n,\widetilde X^m}_t|+ | A^{b_n,\widetilde X^m}_t -A^{b_n,\widetilde X}_t |+ | A^{b_n,\widetilde X}_t -A^{b,\widetilde X}_t |\\
&=:I_1^{(m,n)}+I_2^{(m,n)}+I_3^{(n)}.
\end{align*}
Thanks to Krylov's estimate \eqref{krylov}, we have
\begin{align*}
\| I_1^{(m,n)}\|_{L^2(\widetilde \Omega)}& =\| A^{b_m-b_n,\widetilde X^m}_t\|_{L^2(\widetilde \Omega)} \lesssim \| b_m - b_n \|_{L_T^\infty\bC^{-\beta}} \\
& \leq \| b_m - b \|_{L_T^\infty\bC^{-\beta}}+ \| b - b_n \|_{L_T^\infty\bC^{-\beta}}.
\end{align*}
Based on \eqref{0728:new06}, we obtain
\begin{align*}
    \|I_2^{(m,n)}\|_{L^1(\widetilde \Omega)} \leq \widetilde \mE \int_0^t \left| b_n(s,\widetilde X^m_s) - b_n (s,\widetilde X_s) \right| \dif s  \leq \|b_n\|_{L^\infty_T \sC^1} \mE\int_0^t \left[|\widetilde X^m_s-\widetilde X_s|\wedge 1\right]\dif s\to 0,
\end{align*}
as $m\to \infty$. Thus, we get that
\begin{align*}
    \lim_{m\to \infty}\| A^{b_m,\widetilde X^m}_t-A^{b,\widetilde X}_t\|_{L^1(\widetilde \Omega)}& \le  \lim_{m\to \infty} \|I_1^{(m,n)}\|_{L^2(\widetilde \Omega)} + \lim_{m\to \infty} \|I_2^{(m,n)}\|_{L^1(\widetilde \Omega)} + \|I_3^{(n)}\|_{L^2(\widetilde \Omega)}\\ 
    & \lesssim \|b_n-b\|_{L^\infty_T\bC^{-\beta}}+ \|I_3^{(n)}\|_{L^2(\widetilde \Omega)}\to 0,\quad \text{as $n\to \infty$},
\end{align*}
where provided by \eqref{eq:SD00}. This gives \eqref{0805:00} and the proof of the existence is completed.

\medskip\noindent
{\bf (Uniqueness)} Let $X^{(1)}$ and $X^{(2)}$ be two solutions to SDE \eqref{eq:SDE} in the sense of Definition \ref{Def1} with the same initial distribution. Let $u$ be the unique solution of the backward nonlocal PDE \eqref{0728:PDE} with $f \in C_c^\infty(\mR^d)$. Based on Corollary \ref{geito}, we have
$$
\mE \int_0^T f  (r, X^{(1)}_r) \dif r  =\mE u (0,X^{(1)}_0) =\mE u (0,X^{(2)}_0)= \mE \int_0^T f  (r, X^{(2)}_r) \dif r, 
$$
which implies that $X^{(1)}$ and $X^{(2)}$ have the same one-dimensional time marginal distributions. Hence, by \cite{EK86}, Corollary 4.4.3 (or \cite{SV06}, Theorem 6.2.3, p.147), we get $\mP\circ(X^{(1)})^{-1}=\mP\circ(X^{(2)})^{-1}$ through the standard induction approach. 

The proof is finished.
\end{proof}

\section{Stability of SDEs with $\bC^{-\beta}$ drifts} \label{sec:SDE-2}

In this section, we are going to prove Theorem \ref{thmSt}.

\subsection{Backward Cauchy problem with $f\equiv0$}

Consider the following backward Cauchy problem:
\begin{align*}
\p_tu+\sL^{(\alpha)}_{\sigma}  u+b\cdot\nabla u = 0,\quad u(T)=\varphi,
\end{align*}
where $t\in[0,T)$, and $\varphi\in C^\infty_b(\mR^d)$, and $\sL^{(\alpha)}_{\sigma}$ is defined by \eqref{eq:nolocal}. This subsection is devoted to giving some estimates of $u(t)$ that play an important role in proving stability results. Based on Theorems \ref{BB01} and \ref{BB01:00}, one sees that there is a unique solution $u$. Furthermore, by \cite{CHZ20}, Theorem 1.1, there is a semigroup $P_{t,s}$ such that
 \begin{align}\label{DH00}
u(t)=P_{t,T}\varphi+\int_t^TP_{t,s}(b\cdot\nabla u)(s)\dif s,
\end{align}
and for any $\delta\in[0,\alpha+1)$ and $\eta\in[-\delta,1)$, there is a constant $c>0$ such that for all $0\le t<s\le T$ and $\varphi\in C^\infty_b(\mR^d)$,
 \begin{align}\label{DH01}
 \|P_{t,s}\varphi\|_{\bC^{\delta}} \lesssim_c (s-t)^{-\frac{\delta+\eta}{\alpha}}\|\varphi\|_{\bC^{-\eta}}.
\end{align}

We state the following inequality for $u(t)$.

\bl
Let $T>0$, $\alpha \in (1,2)$, $\beta\in (0,\frac{\alpha-1}{2})$ (resp. $\beta\in [\frac{\alpha-1}{2}, \alpha-1)$), and  $ b\in L_T^\infty C_b^\infty \mR^d$.  Assume that $\sigma$ satisfies {\bf (H$^\sigma$)} with some constant $c_0$. For any $\gamma\in[0,\alpha-1-\beta)$ (resp. $\gamma\in[0,\alpha-1)$) and $\delta\in[0,\alpha-\beta]$, 
there is a constant $ c= c(\Theta,\gamma,\delta,\|b\|_{L_T^\infty\bC^{-\beta}})>0$ (resp. $ c= c(\Theta,\gamma,\delta,\|b\|_{L_T^\infty\bC^{-\beta}},\|\div b\|_{L_T^\infty\bC^{-\beta}})>0$) such that for any $\varphi\in C^\infty_b(\mR^d)$ and $t\in[0,T)$,
\begin{align}\label{DH02}
\|u(t)\|_{\bC^{\delta}} \lesssim_c (T-t)^{-\frac{\delta+\gamma}{\alpha}}\|\varphi\|_{\bC^{-\gamma}}.
\end{align}
\el

\begin{proof} 
We only consider the case $\beta\in [\frac{\alpha-1}{2},\alpha-1)$, since the argument for the case $\beta<\frac{\alpha-1}{2}$ is the same by rewriting $\eps$ as $\eps+\beta$. It follows from \eqref{DH00}, \eqref{DH01} and \eqref{eq:YY01} that for any $\delta\in[0,\alpha+1)$ and $\eps>0$,
\begin{align}
\|u(t)\|_{\bC^\delta}&\lesssim (T-t)^{-\frac{\delta+\gamma}{\alpha}}\|\varphi\|_{\bC^{-\gamma}}+\int_t^T(s-t)^{-\frac{\delta+\beta}{\alpha}}\|b\cdot\nabla u(s)\|_{\bC^{-\beta}}\dif s\nonumber\\ 
& \overset{\eqref{eq:YY01}}{\lesssim } (T-t)^{-\frac{\delta+\gamma}{\alpha}}\|\varphi\|_{\bC^{-\gamma}}+ \int_t^T(s-t)^{-\frac{\delta+\beta}{\alpha}}  \|u(s)\|_{\bC^{ 1+\eps}}\dif s.\label{DH04}
\end{align}
Hence, particularly, taking $0<\eps<\alpha$ and $\delta=1+\eps$ in \eqref{DH04}, we have
\begin{align*}
\|u(t)\|_{\bC^{1+\eps}}
&\lesssim (T-t)^{-\frac{1+\eps+\gamma}{\alpha}}\|\varphi\|_{\bC^{-\gamma}}+\int_t^T(s-t)^{-\frac{1+\beta+\eps}{\alpha}}\|u(s)\|_{\bC^{ 1+\eps}}\dif s,
\end{align*}
which yields
\begin{align*}
\|u(T-(T-t))\|_{\bC^{1+\eps}}
\lesssim & \int_0^{T-t}(T-t-s)^{-\frac{1+\beta+\eps}{\alpha}}\|u(T-s)\|_{\bC^{ 1+\eps}}\dif s\\
 & \quad+ (T-t)^{-\frac{1+\eps+\gamma}{\alpha}}\|\varphi\|_{\bC^{-\gamma}}.
\end{align*}
Then for any $\eps\in (0,\alpha-1-  \beta\vee\gamma)$, one sees that
\begin{align}\label{Hzm05}
\|u(t)\|_{\bC^{1+\eps}}\lesssim (T-t)^{-\frac{1+\eps+\gamma}{\alpha}}\|\varphi\|_{\bC^{-\gamma}}
\end{align}
by Gronwall's inequality of Volterra type (see \cite{CHZ20}, Lemma 2.6, or \cite{Zh10}, Lemma 2.2) with $\frac{1+\eps+\gamma}{\alpha},\frac{1+\beta+\eps}{\alpha}<1$. 

\medskip\noindent
{(i)~{\bf Case one:} $\delta \in [0,\alpha-\beta)$}. Substituting \eqref{Hzm05} to \eqref{DH04}, by the change of variable, we obtain that for any $0<\eps<\alpha-1- \beta\vee \gamma$,
\begin{align*}
\|u(t)\|_{\bC^\delta}&\lesssim (T-t)^{-\frac{\delta+\gamma}{\alpha}}\|\varphi\|_{\bC^{-\gamma}}+\|\varphi\|_{\bC^{-\gamma}} \int_t^T(s-t)^{-\frac{\delta+\beta}{\alpha}} (T-s)^{-\frac{1+\eps+\gamma}{\alpha}}\dif s\\
& \lesssim (T-t)^{-\frac{\delta+\gamma}{\alpha}}\|\varphi\|_{\bC^{-\gamma}} \left( 1+ (T-t)^{\frac{\alpha-1-\beta-\eps}{\alpha}}\int_0^{1} s^{-\frac{\delta+\beta}{\alpha}} (1 -s)^{-\frac{1+\eps+\gamma}{\alpha}}\dif s \right)\\
& \lesssim (T-t)^{-\frac{\delta+\gamma}{\alpha}}\|\varphi\|_{\bC^{-\gamma}} ,
\end{align*}
where we used the definitions of Beta functions with $0<\frac{\delta+\beta}{\alpha},\frac{1+\eps+\gamma}{\alpha}<1$.

\medskip\noindent
{(ii)~{\bf Case two:} $\delta = \alpha-\beta $}. Taking $\eta=\beta$ in \eqref{DH01}, we deduce that for any $\delta \in [0,\alpha+1)$ and $j \geq -1$,
$$
\|\cR_jP_{t,s}(b\cdot\nabla u(s))\|_\infty \lesssim 2^{-\delta j} (s-t)^{-\frac{  \delta+\beta}{\alpha}} \|b\cdot\nabla u(s)\|_{\bC^{-\beta}}.
$$
Thus, letting $\delta=0,  \alpha$ in turn, one sees that
\begin{align*}
\|\cR_jP_{t,s}(b\cdot\nabla u(s))\|_\infty &\lesssim \left[(2^{-\alpha j}(s-t)^{-\frac{  \alpha+\beta}{\alpha}})\wedge (s-t)^{-\frac{\beta}{\alpha}}\right] \|b\cdot\nabla u(s)\|_{\bC^{-\beta}},
\end{align*}
which, together with \eqref{eq:YY01} and \eqref{Hzm05}, derives that for any $0<\eps<\alpha-1-  \beta\vee \gamma$,
\begin{align*}
\|\cR_jP_{t,s}(b\cdot\nabla u(s))\|_\infty & \lesssim  \left[(2^{-\alpha j}(s-t)^{-\frac{  \alpha+\beta}{\alpha}})\wedge (s-t)^{-\frac{\beta}{\alpha}}\right]\|u(s)\|_{\bC^{1+\eps}}\\ 
& \lesssim  \left[(2^{-\alpha j}(s-t)^{-\frac{  \alpha+\beta}{\alpha}})\wedge (s-t)^{-\frac{\beta}{\alpha}}\right] (T-s)^{-\frac{1+\eps+\gamma}{\alpha}}\|\varphi\|_{\bC^{-\gamma}}.
\end{align*}
Consequently, by the change of variable and Lemma \ref{Ap:lem}, we obtain that 
\begin{align*} 
& \left\|\cR_j \int_t^T P_{t,s}(b\cdot\nabla u(s))\dif s\right\|_\infty\\
\lesssim&\,  \|\varphi\|_{\bC^{-\gamma}}
\int_0^{T-t} \left[(2^{- \alpha j}s^{-\frac{  \alpha+\beta}{\alpha}})\wedge s^{-\frac{\beta}{\alpha}}\right] (T-t-s)^{-\frac{1+\eps+\gamma}{\alpha}}\dif s\\
\lesssim &\, \|\varphi\|_{\bC^{-\gamma}} 2^{-j(\alpha-\beta)}(T-t)^{-\frac{1+\eps+\gamma}{\alpha}}
 \lesssim  \|\varphi\|_{\bC^{-\gamma}} 2^{-j(\alpha-\beta)}(T-t)^{-\frac{\alpha-\beta+\gamma}{\alpha}},
\end{align*}
which leads to the desired estimates trivially. 

Combining these two cases, we complete the proof.
\end{proof}

\subsection{Proof of Theorem \ref{thmSt}}
Based on the uniqueness, we can assume $b_1,b_2\in L^\infty_TC^\infty_b (\mR^d)$. Then,  for each $i=1,2$, it is well-known (see \cite{CZZ21} for example) that,  for any $x \in \mR^d$, there is a unique solutions $X^i$,  to the following classical SDE:
\begin{align*}
X^i_{t}=x+\int_0^tb_i(r,X^i_{r})\dif r+\int_0^t \sigma(r,X_{r-})\dif L_r^{(\alpha)}.
\end{align*}
It suffices to estimate 
$
\left|\mE\varphi(X^2_{t})-\mE\varphi(X^1_{t})\right|
$
for any $\varphi\in C^\infty_b(\mR^d)$ and $t\in(0,T]$.  Using It\^o-Tanaka trick, letting $\varphi\in C^\infty_b(\mR^d)$ be the terminal condition of the following backward PDE:
\begin{align}\label{eq:GH01}
\p_r u^t+\sL^{(\alpha)}_{\sigma}u^t+b_1\cdot\nabla u^t=0,\quad u^t(t)=\varphi,
\end{align}
where $r\in (0,t)$, $u^t$ is the shifted function $u^t (r,x): = u(t-r,x)$, and $\sL^{(\alpha)}_{\sigma}$ is defined by \eqref{eq:nolocal}. It follows \eqref{DH02} that for every $\widetilde \delta\in[0,\alpha-\beta ]$,
\begin{align}\label{DH07}
\|u^t(r)\|_{\bC^{\tilde\delta}}\lesssim (t-r)^{-{\tilde\delta}/{\alpha}}\|\varphi\|_\infty.
\end{align}
By It\^{o}'s formula (cf. \cite{IW89}), we have that for $i=1,2$,
 \begin{align*}
u^t(t,X^i_t) &- u^t(0,x) = \ \ 
 \int_0^t (\p_r u^t) (r,X^i_r) \dif r 
+ \int_0^t b_i(r, X^i_r) \cdot  \nabla  u^t (r,X^i_r) \dif r\\
 &+  \int_0^{t} \int_{\mR^d} \( u^t (r,X^i_{r-}+\sigma(r,X^i_{r-})z) - u^t (r,X^i_{r-}) \) \widetilde N(\dif r, \dif z)\\
&+ \int_0^t \int_{\mR^d}  \( u^t (r,X^i_r+\sigma(r,X^i_{r-})z) - u^t (r,X^i_r) \\
& \qquad\qquad \qquad\qquad\qquad\qquad- \sigma(r,X^i_{r-})z \cdot \nabla u^t(r,X^i_r) \) \nu^{(\alpha)}( \dif z) \dif r.
\end{align*}
Observe that the third term on the right-hand side of the above equality is a martingale. 
Thus, we obtain  
\begin{align*}
\mE u^t (t,X^i_t ) = & \int_0^t  (\p_r +\sL^{(\alpha)}_{\sigma}  + b_i\cdot  \nabla\) u^t (r,X^i_r ) \dif r + u^t(0,x)   \\
 \overset{\eqref{eq:GH01}}{=}& \int_0^t  \(( b_i -b_1)\cdot \nabla u^t \) (r,X^i_r ) \dif r+ u^t(0,x)  ,
\end{align*}
which implies 
\begin{align}\label{eq:WM01}
\mE\varphi(X^1_{t}) =\mE u^t(t,X^1_{t})=u^t(0,x),
\end{align}
and then
\begin{align}\label{eq:WM00}
\mE\varphi(X^2_{t})-\mE\varphi(X^1_{t})=\mE\int_0^t
\((b_2-b_1)\cdot\nabla u^t\)(r,X^2_r)\dif r.
\end{align}
Now the key point is to estimate
$$
\mE ((b_2-b_1)\cdot\nabla u^t)(r,\cdot )(X^2_r).
$$ 
Define
$
g_{r,t}(x):=((b_2-b_1)\cdot\nabla u^t)(r,\cdot)(x)
$
with $ r <t \leq T$. 
Using the It\^o-Tanaka trick again, we consider the following PDE:
\begin{align*}
\p_s w +\sL^{(\alpha)}_{\sigma} w +b_2\cdot\nabla w =0,\quad w (r)=g_{r,t},
\end{align*}
where $s \in (0,r)$ and $g$ belongs to $C_b^\infty(\mR^d)$. Adopting the same argument as \eqref{eq:WM01}, one sees that
$$
\mE g_{r,t}(X^2_{r})=w(0,x),
$$
which by \eqref{DH02} implies that for any $\theta\in [\beta,\alpha-1-\beta \1_{\{\beta<\frac{\alpha-1}{2}\}})$,
\begin{align*}
|\mE g_{r,t}(X^2_{r})|=|w(0,x)|\le \|w(0)\|_{\bC^\eps}\lesssim_c  r^{-\frac{\theta}{\alpha}}\|g\|_{\bC^{-\theta}},
\end{align*}
which derives that
\begin{align*}
\left|\mE\varphi(X^2_{t})-\mE\varphi(X^1_{t})\right| 
\overset{\eqref{eq:WM00}}{\lesssim}\int_0^t r^{-\frac{ \theta}{\alpha}}\| ((b_2-b_1)\cdot\nabla u^t)(r)\|_{\bC^{-\theta}}\dif r.
\end{align*}
From \eqref{eq:YY01} and \eqref{DH07}, for the case $\beta\in (0,\frac{\alpha-1}{2})$, we have that for any  $ 0<\eps< \alpha - 1 - \theta$,
\begin{align*}
\left|\mE\varphi(X^2_{t})-\mE\varphi(X^1_{t})\right|&\lesssim \|b_1-b_2\|_{\bC^{-\theta}}\int_0^t r^{-\frac{ \theta}{\alpha}}\| u^t(r)\|_{\bC^{\theta+1+\eps}}\dif r\\
& \overset{\eqref{DH07}}{\lesssim} \|b_1-b_2\|_{\bC^{-\theta}}\|\varphi\|_\infty\int_0^t r^{-\frac{ \theta}{\alpha}}(t-r)^{-\frac{\theta+1+\eps}{\alpha}}\dif r\\
& = \|b_1-b_2\|_{\bC^{-\theta}}\|\varphi\|_\infty   t^{\frac{\alpha-1-2\theta-\eps}{\alpha} } \int_0^1 r^{-\frac{\theta  }{\alpha}}(1-r)^{-\frac{\theta+1+\eps}{\alpha}}\dif r \\
& \lesssim \|b_1-b_2\|_{\bC^{-\theta}}\|\varphi\|_\infty t^{\frac{\alpha-1-2\theta-\eps}{\alpha} };
\end{align*}
and samely, for the case $\beta\in [ \frac{\alpha-1}{2}, \alpha-1 )$,
\begin{align*}
\left|\mE\varphi(X^2_{t})-\mE\varphi(X^1_{t})\right|&\lesssim \left(\|b_1-b_2\|_{\bC^{-\theta}}+\|\div b_1-\div b_2\|_{\bC^{-\theta}}\right)\int_0^t r^{-\frac{\theta}{\alpha}}\| u^t(r)\|_{\bC^{1+\eps}}\dif r\\
& \lesssim \left(\|b_1-b_2\|_{\bC^{-\theta}}+\|\div b_1-\div b_2\|_{\bC^{-\theta}}\right)\|\varphi\|_\infty t^{\frac{\alpha-1-\theta-\eps}{\alpha} }.
\end{align*}
This completes the proof.
 
\begin{appendix}

\section{Proof of Lemma \ref{Lem31}: Boundedness of operator $\sL_\sigma^{(\alpha)}$}


We first introduce some notation. For any $x,y,z\in\mR^d$, define  
\begin{align*}
\sigma_t^y (x):= \sigma(t,x+y),\ \ \hbox{and}\ \ \Lambda_{t,z}^y(x):=x+\sigma_t^y (x)  z.
\end{align*} 
It is easy to see that, under the condition {\bf (H$^\sigma$)} with some constant $c_0$, for any $t\ge0$, and $x_1,x_2\in\mR^d$, and $|z|\leq\frac{1}{2}c_0^{-1}$,
\begin{align*} 
\frac{1}{2}|x_1-x_2|\le |\Lambda_{t,z}^y(x_1)-\Lambda_{t,z}^y(x_2)|\le 2|x_1-x_2|.
\end{align*}
Define 
\begin{align*}
\sD_{\sigma_t^y,z} f(x):=f(x+\sigma_t^y(x)z)-f(x+\sigma_t^y(0)z)-(\sigma_t^y(x)-\sigma_t^y(0))z\cdot\nabla f(x).
\end{align*}
Then, by \cite{CHZ20}, Lemma 2.2, for any $f\in \bC^1 (\mR^d)$, $g \in L^1(\mR^d)$ with $\nabla g,\nabla^2 g  \in L^1(\mR^d)$, and $\theta\in[0,1]$, we have
\begin{align}\label{EE02}
|\<\sD_{\sigma_t^y,z} f,g\>|\leq {c}_{d,\theta}|z|^{1+\theta}\|f\|_{\bC^\theta}\(\sum_{j=0}^1\mu_j(|\nabla^j g|)+\mu_{1+\theta}(|\nabla^2g|)^\theta\mu_{1+\theta}(|\nabla g|)^{1-\theta}\),
\end{align}
where the constant $c_{d,\theta}>0$ is independent of the variables $t,y,z$, and
\begin{align*}
\mu_\theta(\dif x):=(|x|\wedge1)^\theta\dif x,\ \ \hbox{and} \ \  \mu_\theta(f):=\int_{\mR^d}f(x)\mu_\theta(\dif x).
\end{align*}

Now we give the

\begin{proof}[Proof of Lemma \ref{Lem31}]
For simplicity, we drop the time variable $t$ and the superscript $\alpha$ in $\nu^{(\alpha)}$ in the following proof. Observe that
$$
\cR_j\sL_{\sigma}^{(\alpha)} u(x) = [\cR_j, \sL_{\sigma}^{(\alpha)}]  u(x)+\sL_{\sigma}^{(\alpha)} \cR_j u(x) , 
$$
where $[\sA_1, \sA_2] := \sA_1 \sA_2  - \sA_2 \sA_1$ is called commutator operator. Next, we estimate these two terms in turn.

\medskip
\noindent
{\bf (Step 1)} Define
$$
u^x(y):=u(y+x),\quad \sigma^x(y):=\sigma(y+x).
$$
Recall the definition \eqref{eq:Block} of block operators $\cR_j$. By the change of variable, we have
\begin{align*}
I&:= \int_{\mR^d}\check\psi_j(x-y)\(u(y+\sigma(y)z)-u(y+\sigma(x)z)-(\sigma(y)-\sigma(x))z\cdot\nabla u(y)\)\dif y\\
&=  \int_{\mR^d} \check\psi_j(-y)\(u^x(y+\sigma^x(y)z)-u^x(y+\sigma^x(0)z)-(\sigma^x(y)-\sigma^x(0))z\cdot\nabla u^x(y)\)\dif y\\
&=\int_{\mR^d} \check\psi_j(-y)\sD_{\sigma^x,z} u^x(y) \dif y,
\end{align*}
which derives that
$$
 [\cR_j, \sL_{\sigma}^{(\alpha)}]  u(x) =   \int_{\mR^d}\<\sD_{\sigma^x,z} u^x,\check\psi_j(-\cdot)\> \nu  (\dif z).
 $$
For $|z|\le \delta\leq\frac{1}{2c_0}$, based on the scaling property 
\begin{align}
\mu_{\delta}\label{eq:RQ00}(|\nabla^k\check\psi_j|)\lesssim 2^{(k-\delta)j}, \ \ \delta\ge0,~~k\in\mN_0,
\end{align}
and \eqref{EE02} with $\theta\in(\alpha-1,  (\alpha-\beta)\wedge 1)$, one sees that
\begin{align*}
\sup_x|\<\sD_{\sigma^x,z} u^x,\check\psi_j(-\cdot)\> | &\lesssim |z|^{1+\theta}\|u\|_{\bC^{\theta}}\(1+2^{\theta(2-1-\theta)j}2^{(1-\theta)(1-1-\theta)j}\)\\
&
\lesssim |z|^{1+\theta}\|u\|_{\bC^{\alpha-\beta}}.
\end{align*}
For $|z|>\delta$, by the definition and the integral by parts, under the condition {\bf(H$^\sigma$)}, one sees that
\begin{align*}
|\langle\sD_{\sigma^x,z} u^x ,\check\psi_j(-\cdot)\rangle |&\le 2\|u^x\|_\infty \|\check{\psi_j}\|_1+\left|\int_{\mR^d}\check{\psi_j}(-y)(\sigma^x(y)-\sigma^x(0))z\cdot\nabla u^x(y)\dif y\right|\\
&\le 2\|u^x\|_\infty \|\check{\psi_j}\|_1+\left|\int_{\mR^d}\nabla\check{\psi_j}(-y)\cdot (\sigma^x(y)-\sigma^x(0))z u^x(y)\dif y\right|\\
&\qquad \qquad \qquad\qquad +\left|\int_{\mR^d}\check{\psi_j}(-y){\rm div} \sigma^x(y)\cdot z u^x(y)\dif y\right|\\
&\lesssim \|u^x\|_\infty \|\check{\psi_j}\|_1+|z|\|u^x\|_\infty\|\nabla\sigma\|_\infty
\left( \mu_1(|\nabla\check{\psi_j}|)+\|\check{\psi_j}\|_1 \right)\\
&\overset{\eqref{eq:RQ00}}{\lesssim}  (1+|z|) \|u\|_\infty \lesssim ( 1+|z|) \|u\|_{\bC^{\alpha-\beta}},
\end{align*}
where we used the fact $\alpha-\beta >0$ in the last inequality. Hence, by \eqref{eq:AE02}, we obtain that
\begin{align}\label{EE00}
\|[\cR_j, \sL_{\sigma}^{(\alpha)}]  u\|_\infty 
\leq \int_{\mR^d}\sup_x|\<\sD_{\sigma^x,z} u^x,\check\psi_j(-\cdot)\>|\nu (\dif z)  
\lesssim \|u\|_{\bC^{\alpha-\beta}}.
\end{align}

\medskip
\noindent
{\bf (Step 2)} On the other hand, by Bernstein's inequality (see Lemma \ref{Bern}),  we have
\begin{align*}
\sup_x|\cR_ju(x+\sigma(x)z)-\cR_ju(x)-\sigma(x)z\cdot\nabla \cR_j u(x)|\lesssim (|2^jz|\wedge |2^j z|^2 )\|\cR_ju\|_{\infty},
\end{align*}
which implies that
\begin{align*}
\|\sL_{\sigma}^{(\alpha)} \cR_j u\|_\infty \lesssim \|\cR_ju\|_{\infty} \int_{\mR^d}(|2^jz|\wedge|2^jz|^2)\nu (\dif z)\lesssim 2^{\alpha j}\|\cR_ju\|_{\infty},
\end{align*}
provided by the scaling property \eqref{eq:scal} and \eqref{eq:AE02}.
Therefore, combining this with \eqref{EE00}, we have
\begin{align*} 
2^{-\beta j} \|\cR_j\sL_\sigma^{(\alpha)} u\|_{\infty}
\lesssim 2^{-\beta j} \|u\|_{\bC^{\alpha-\beta}}+2^{(\alpha-\beta) j}\|\cR_ju\|_{\infty}
\lesssim \|u\|_{\bC^{\alpha-\beta}},
\end{align*}
which derives \eqref{EE01} by taking the supremum of $j$. The proof is completed.
\end{proof}

\section{Auxiliary inequalities}

\bl\label{Ap:lem}
Let $T>0$. For any $0<\beta,\gamma<1<\alpha$, there is a constant $c>0$ such that for all $ \lambda>0$ and $t\in(0,T]$,
\begin{align*}
\int_0^t \left[(\lambda s^{-\alpha})\wedge s^{-\beta}\right](t-s)^{-\gamma}\dif s \lesssim_c  \lambda^{\frac{1-\beta}{\alpha-\beta}}t^{-\gamma}.
\end{align*}
\el

\begin{proof}
First of all, by the change of variable, we have
\begin{align*}
\int_0^t \left[(\lambda s^{-\alpha})\wedge s^{-\beta}\right](t-s)^{-\gamma}\dif s=t^{1-\gamma-\beta}\int_0^1 \left[(\lambda t^{-(\alpha-\beta)} s^{-\alpha})\wedge s^{-\beta}\right](1-s)^{-\gamma}\dif s.
\end{align*}
Therefore, it is sufficient to show 
\begin{align*}
\sI:=\int_0^1 \left[(\lambda s^{-\alpha})\wedge s^{-\beta}\right](1-s)^{-\gamma}\dif s\lesssim \lambda^{\frac{1-\beta}{\alpha-\beta}}.
\end{align*}
Indeed, for any $0<\beta,\gamma<1<\alpha$,
\begin{align*}
\sI& \lesssim \int_0^{\tfrac12} \left[(\lambda s^{-\alpha})\wedge s^{-\beta}\right]\dif s+ \left[(\lambda (\tfrac12)^{-\alpha})\wedge (\tfrac12)^{-\beta}\right]\int_{\tfrac12}^1 (1-s)^{-\gamma}\dif s\\
& \lesssim \lambda^{\frac{1-\beta}{\alpha-\beta}}\int_0^{\tfrac{1}{2\lambda^{1/(\alpha-\beta)}}} r^{-\beta}\left[r^{-\alpha+\beta}\wedge 1\right]\dif r+ \lambda \wedge 1\\
& \lesssim \lambda^{\frac{1-\beta}{\alpha-\beta}}\int_0^{\infty} r^{-\beta}\left[r^{-\alpha+\beta}\wedge 1\right]\dif r+ \lambda\wedge 1 \lesssim \lambda^{\frac{1-\beta}{\alpha-\beta}}, 
\end{align*}
where we used the change of variable $s=\lambda^{1/(\alpha-\beta)}r$ in the second inequality. This completes the proof.
\end{proof}

\end{appendix}

\begin{acks}[Acknowledgments]
 
The first version of this paper was completed while the second author, Mingyan Wu, was a postdoctoral fellow at Huazhong University of Science and Technology under the supervision of Prof. Fuke Wu. The authors are deeply grateful to Prof. Fuke Wu for his valuable suggestions and for correcting some errors. 
\end{acks}

\begin{funding}
Zimo Hao is grateful for the financial support of NNSFC grants of China (Nos. 12131019, 11731009), and the DFG through the CRC 1283/2 2021 - 317210226 "Taming uncertainty and profiting from randomness and low regularity in analysis, stochastics and their applications". Mingyan Wu is partially supported by the National Natural Science Foundation of China (Grant No. 12201227 and No. 61873320) and the
Fundamental Research Funds for the Central Universities (Grant No. 20720240128).
\end{funding}

\newcommand{\eprint}[1]{\href{https://arxiv.org/abs/#1}{\it{#1}}}

\bibliographystyle{imsart-nameyear} 
\bibliography{refs_Hao-Wu.bib}       

\begin{thebibliography}{39}

\bibitem[\protect\citeauthoryear{Athreya, Butkovsky and Mytnik}{2020}]{ABM18}
\begin{barticle}[author]
\bauthor{\bsnm{Athreya},~\bfnm{Siva}\binits{S.}},
  \bauthor{\bsnm{Butkovsky},~\bfnm{Oleg}\binits{O.}} \AND
  \bauthor{\bsnm{Mytnik},~\bfnm{Leonid}\binits{L.}}
(\byear{2020}).
\btitle{Strong existence and uniqueness for stable stochastic differential
  equations with distributional drift}.
\bjournal{Ann. Probab.}
\bvolume{48}
\bpages{178--210}.
\bdoi{10.1214/19-AOP1358}
\bmrnumber{4079434}
\end{barticle}
\endbibitem

\bibitem[\protect\citeauthoryear{Bahouri, Chemin and Danchin}{2011}]{BCD11}
\begin{bbook}[author]
\bauthor{\bsnm{Bahouri},~\bfnm{Hajer}\binits{H.}},
  \bauthor{\bsnm{Chemin},~\bfnm{Jean-Yves}\binits{J.-Y.}} \AND
  \bauthor{\bsnm{Danchin},~\bfnm{Rapha\"{e}l}\binits{R.}}
(\byear{2011}).
\btitle{Fourier analysis and nonlinear partial differential equations}.
\bseries{Grundlehren der mathematischen Wissenschaften [Fundamental Principles
  of Mathematical Sciences]}
\bvolume{343}.
\bpublisher{Springer, Heidelberg}.
\bdoi{10.1007/978-3-642-16830-7}
\bmrnumber{2768550}
\end{bbook}
\endbibitem

\bibitem[\protect\citeauthoryear{Bass and Chen}{2001}]{BC01}
\begin{barticle}[author]
\bauthor{\bsnm{Bass},~\bfnm{Richard~F.}\binits{R.~F.}} \AND
  \bauthor{\bsnm{Chen},~\bfnm{Zhen-Qing}\binits{Z.-Q.}}
(\byear{2001}).
\btitle{Stochastic differential equations for {D}irichlet processes}.
\bjournal{Probab. Theory Related Fields}
\bvolume{121}
\bpages{422--446}.
\bdoi{10.1007/s004400100151}
\bmrnumber{1867429}
\end{barticle}
\endbibitem

\bibitem[\protect\citeauthoryear{Bony}{1981}]{Bo81}
\begin{barticle}[author]
\bauthor{\bsnm{Bony},~\bfnm{Jean-Michel}\binits{J.-M.}}
(\byear{1981}).
\btitle{Calcul symbolique et propagation des singularit\'{e}s pour les
  \'{e}quations aux d\'{e}riv\'{e}es partielles non lin\'{e}aires}.
\bjournal{Ann. Sci. \'{E}cole Norm. Sup. (4)}
\bvolume{14}
\bpages{209--246}.
\bmrnumber{631751}
\end{barticle}
\endbibitem

\bibitem[\protect\citeauthoryear{Cannizzaro and Chouk}{2018}]{CC18}
\begin{barticle}[author]
\bauthor{\bsnm{Cannizzaro},~\bfnm{Giuseppe}\binits{G.}} \AND
  \bauthor{\bsnm{Chouk},~\bfnm{Khalil}\binits{K.}}
(\byear{2018}).
\btitle{Multidimensional {SDE}s with singular drift and universal construction
  of the polymer measure with white noise potential}.
\bjournal{Ann. Probab.}
\bvolume{46}
\bpages{1710--1763}.
\bdoi{10.1214/17-AOP1213}
\bmrnumber{3785598}
\end{barticle}
\endbibitem

\bibitem[\protect\citeauthoryear{Cavallazzi}{}]{Ca22}
\begin{barticle}[author]
\bauthor{\bsnm{Cavallazzi},~\bfnm{Thomas}\binits{T.}}
\btitle{Quantitative weak propagation of chaos for stable-driven McKean-Vlasov
  SDEs}.
\bjournal{Available at arXiv:\eprint{2212.01079}}.
\end{barticle}
\endbibitem

\bibitem[\protect\citeauthoryear{Chaudru~de Raynal and Menozzi}{2022}]{CM19}
\begin{barticle}[author]
\bauthor{\bparticle{Chaudru~de}
  \bsnm{Raynal},~\bfnm{Paul-\'{E}ric}\binits{P.-E.}} \AND
  \bauthor{\bsnm{Menozzi},~\bfnm{St\'{e}phane}\binits{S.}}
(\byear{2022}).
\btitle{On multidimensional stable-driven stochastic differential equations
  with {B}esov drift}.
\bjournal{Electron. J. Probab.}
\bvolume{27}
\bpages{Paper No. 163, 52}.
\bdoi{10.1214/22-ejp864}
\bmrnumber{4525442}
\end{barticle}
\endbibitem

\bibitem[\protect\citeauthoryear{Chen, Hao and Zhang}{2020}]{CHZ20}
\begin{barticle}[author]
\bauthor{\bsnm{Chen},~\bfnm{Zhen-Qing}\binits{Z.-Q.}},
  \bauthor{\bsnm{Hao},~\bfnm{Zimo}\binits{Z.}} \AND
  \bauthor{\bsnm{Zhang},~\bfnm{Xicheng}\binits{X.}}
(\byear{2020}).
\btitle{H\"{o}lder regularity and gradient estimates for {SDE}s driven by
  cylindrical {$\alpha$}-stable processes}.
\bjournal{Electron. J. Probab.}
\bvolume{25}
\bpages{Paper No. 137, 23}.
\bdoi{10.1214/20-ejp542}
\bmrnumber{4179301}
\end{barticle}
\endbibitem

\bibitem[\protect\citeauthoryear{Chen and Zhang}{2018}]{CZ18b}
\begin{barticle}[author]
\bauthor{\bsnm{Chen},~\bfnm{Zhen-Qing}\binits{Z.-Q.}} \AND
  \bauthor{\bsnm{Zhang},~\bfnm{Xicheng}\binits{X.}}
(\byear{2018}).
\btitle{{$L^p$}-maximal hypoelliptic regularity of nonlocal kinetic
  {F}okker-{P}lanck operators}.
\bjournal{J. Math. Pures Appl. (9)}
\bvolume{116}
\bpages{52--87}.
\bdoi{10.1016/j.matpur.2017.10.003}
\bmrnumber{3826548}
\end{barticle}
\endbibitem

\bibitem[\protect\citeauthoryear{Chen, Zhang and Zhao}{2021}]{CZZ21}
\begin{barticle}[author]
\bauthor{\bsnm{Chen},~\bfnm{Zhen-Qing}\binits{Z.-Q.}},
  \bauthor{\bsnm{Zhang},~\bfnm{Xicheng}\binits{X.}} \AND
  \bauthor{\bsnm{Zhao},~\bfnm{Guohuan}\binits{G.}}
(\byear{2021}).
\btitle{Supercritical {SDE}s driven by multiplicative stable-like {L}\'{e}vy
  processes}.
\bjournal{Trans. Amer. Math. Soc.}
\bvolume{374}
\bpages{7621--7655}.
\bdoi{10.1090/tran/8343}
\bmrnumber{4328678}
\end{barticle}
\endbibitem

\bibitem[\protect\citeauthoryear{Delarue and Diel}{2016}]{DD16}
\begin{barticle}[author]
\bauthor{\bsnm{Delarue},~\bfnm{Fran\c{c}ois}\binits{F.}} \AND
  \bauthor{\bsnm{Diel},~\bfnm{Roland}\binits{R.}}
(\byear{2016}).
\btitle{Rough paths and 1d {SDE} with a time dependent distributional drift:
  application to polymers}.
\bjournal{Probab. Theory Related Fields}
\bvolume{165}
\bpages{1--63}.
\bdoi{10.1007/s00440-015-0626-8}
\bmrnumber{3500267}
\end{barticle}
\endbibitem

\bibitem[\protect\citeauthoryear{Ethier and Kurtz}{1986}]{EK86}
\begin{bbook}[author]
\bauthor{\bsnm{Ethier},~\bfnm{Stewart~N.}\binits{S.~N.}} \AND
  \bauthor{\bsnm{Kurtz},~\bfnm{Thomas~G.}\binits{T.~G.}}
(\byear{1986}).
\btitle{Markov processes}.
\bseries{Wiley Series in Probability and Mathematical Statistics: Probability
  and Mathematical Statistics}.
\bpublisher{John Wiley \& Sons, Inc., New York}
\bnote{Characterization and convergence}.
\bdoi{10.1002/9780470316658}
\bmrnumber{838085}
\end{bbook}
\endbibitem

\bibitem[\protect\citeauthoryear{Friz and Zhang}{2018}]{FZ18}
\begin{barticle}[author]
\bauthor{\bsnm{Friz},~\bfnm{Peter~K.}\binits{P.~K.}} \AND
  \bauthor{\bsnm{Zhang},~\bfnm{Huilin}\binits{H.}}
(\byear{2018}).
\btitle{Differential equations driven by rough paths with jumps}.
\bjournal{J. Differential Equations}
\bvolume{264}
\bpages{6226--6301}.
\bdoi{10.1016/j.jde.2018.01.031}
\bmrnumber{3770049}
\end{barticle}
\endbibitem

\bibitem[\protect\citeauthoryear{Gr\"afner and Perkowski}{}]{GP23}
\begin{barticle}[author]
\bauthor{\bsnm{Gr\"afner},~\bfnm{L.}\binits{L.}} \AND
  \bauthor{\bsnm{Perkowski},~\bfnm{N.}\binits{N.}}
\btitle{Energy solutions and generators of singular SPDEs}.
\bjournal{Lecture note,
  \url{https://www.mis.mpg.de/fileadmin/pdf/slides_randompde_5824.pdf}}.
\bdoi{https://www.mis.mpg.de/fileadmin/pdf/slides_randompde_5824.pdf}
\end{barticle}
\endbibitem

\bibitem[\protect\citeauthoryear{Gubinelli, Imkeller and
  Perkowski}{2015}]{GIP15}
\begin{barticle}[author]
\bauthor{\bsnm{Gubinelli},~\bfnm{Massimiliano}\binits{M.}},
  \bauthor{\bsnm{Imkeller},~\bfnm{Peter}\binits{P.}} \AND
  \bauthor{\bsnm{Perkowski},~\bfnm{Nicolas}\binits{N.}}
(\byear{2015}).
\btitle{Paracontrolled distributions and singular {PDE}s}.
\bjournal{Forum Math. Pi}
\bvolume{3}
\bpages{e6, 75}.
\bdoi{10.1017/fmp.2015.2}
\bmrnumber{3406823}
\end{barticle}
\endbibitem

\bibitem[\protect\citeauthoryear{Hao, Wang and Wu}{2024}]{HWW20}
\begin{barticle}[author]
\bauthor{\bsnm{Hao},~\bfnm{Zimo}\binits{Z.}},
  \bauthor{\bsnm{Wang},~\bfnm{Zhen}\binits{Z.}} \AND
  \bauthor{\bsnm{Wu},~\bfnm{Mingyan}\binits{M.}}
(\byear{2024}).
\btitle{Schauder {E}stimates for {N}onlocal {E}quations with {S}ingular
  {L}\'{e}vy {M}easures}.
\bjournal{Potential Anal.}
\bvolume{61}
\bpages{13--33}.
\bdoi{10.1007/s11118-023-10101-9}
\bmrnumber{4758470}
\end{barticle}
\endbibitem

\bibitem[\protect\citeauthoryear{Hao and Zhang}{}]{HZ23}
\begin{barticle}[author]
\bauthor{\bsnm{Hao},~\bfnm{Zimo}\binits{Z.}} \AND
  \bauthor{\bsnm{Zhang},~\bfnm{Xicheng}\binits{X.}}
\btitle{SDEs with supercritical distributional drifts}.
\bjournal{To appear in {\it } Available at arXiv:\eprint{2312.11145}}.
\bdoi{https://arxiv.org/abs/2312.11145}
\end{barticle}
\endbibitem

\bibitem[\protect\citeauthoryear{Ikeda and Watanabe}{1989}]{IW89}
\begin{bbook}[author]
\bauthor{\bsnm{Ikeda},~\bfnm{Nobuyuki}\binits{N.}} \AND
  \bauthor{\bsnm{Watanabe},~\bfnm{Shinzo}\binits{S.}}
(\byear{1989}).
\btitle{Stochastic differential equations and diffusion processes},
\bedition{second} ed.
\bseries{North-Holland Mathematical Library}
\bvolume{24}.
\bpublisher{North-Holland Publishing Co., Amsterdam; Kodansha, Ltd., Tokyo}.
\bmrnumber{1011252}
\end{bbook}
\endbibitem

\bibitem[\protect\citeauthoryear{Kallenberg}{2021}]{Ka3rd}
\begin{bbook}[author]
\bauthor{\bsnm{Kallenberg},~\bfnm{Olav}\binits{O.}}
(\byear{2021}).
\btitle{Foundations of modern probability (Third edition)}.
\bseries{Probability Theory and Stochastic Modelling}
\bvolume{99}.
\bpublisher{Springer, Cham}.
\bdoi{10.1007/978-3-030-61871-1}
\bmrnumber{4226142}
\end{bbook}
\endbibitem

\bibitem[\protect\citeauthoryear{Kremp and Perkowski}{2022}]{KP20}
\begin{barticle}[author]
\bauthor{\bsnm{Kremp},~\bfnm{Helena}\binits{H.}} \AND
  \bauthor{\bsnm{Perkowski},~\bfnm{Nicolas}\binits{N.}}
(\byear{2022}).
\btitle{Multidimensional {SDE} with distributional drift and {L}\'{e}vy noise}.
\bjournal{Bernoulli}
\bvolume{28}
\bpages{1757--1783}.
\bdoi{10.3150/21-bej1394}
\bmrnumber{4411510}
\end{barticle}
\endbibitem

\bibitem[\protect\citeauthoryear{Kremp and Perkowski}{2025}]{KP25}
\begin{barticle}[author]
\bauthor{\bsnm{Kremp},~\bfnm{Helena}\binits{H.}} \AND
  \bauthor{\bsnm{Perkowski},~\bfnm{Nicolas}\binits{N.}}
(\byear{2025}).
\btitle{Rough weak solutions for singular L\'evy SDEs}.
\bjournal{Probab. Theory Related Fields}.
\bdoi{10.1007/s00440-025-01371-y}
\end{barticle}
\endbibitem

\bibitem[\protect\citeauthoryear{Kunita}{2004}]{Ku04}
\begin{bincollection}[author]
\bauthor{\bsnm{Kunita},~\bfnm{Hiroshi}\binits{H.}}
(\byear{2004}).
\btitle{Stochastic differential equations based on {L}\'{e}vy processes and
  stochastic flows of diffeomorphisms}.
In \bbooktitle{Real and stochastic analysis}.
\bseries{Trends Math.}
\bpages{305--373}.
\bpublisher{Birkh\"{a}user Boston, Boston, MA}.
\bmrnumber{2090755}
\end{bincollection}
\endbibitem

\bibitem[\protect\citeauthoryear{Kurtz and Protter}{1991}]{KP91}
\begin{barticle}[author]
\bauthor{\bsnm{Kurtz},~\bfnm{Thomas~G.}\binits{T.~G.}} \AND
  \bauthor{\bsnm{Protter},~\bfnm{Philip}\binits{P.}}
(\byear{1991}).
\btitle{Weak limit theorems for stochastic integrals and stochastic
  differential equations}.
\bjournal{Ann. Probab.}
\bvolume{19}
\bpages{1035--1070}.
\bmrnumber{1112406}
\end{barticle}
\endbibitem

\bibitem[\protect\citeauthoryear{L\'epingle}{1976}]{Le76}
\begin{barticle}[author]
\bauthor{\bsnm{L\'epingle},~\bfnm{D.}\binits{D.}}
(\byear{1976}).
\btitle{La variation d'ordre {$p$} des semi-martingales}.
\bjournal{Z. Wahrscheinlichkeitstheorie und Verw. Gebiete}
\bvolume{36}
\bpages{295--316}.
\bdoi{10.1007/BF00532696}
\bmrnumber{420837}
\end{barticle}
\endbibitem

\bibitem[\protect\citeauthoryear{Ling and Zhao}{2022}]{LZ22}
\begin{barticle}[author]
\bauthor{\bsnm{Ling},~\bfnm{Chengcheng}\binits{C.}} \AND
  \bauthor{\bsnm{Zhao},~\bfnm{Guohuan}\binits{G.}}
(\byear{2022}).
\btitle{Nonlocal elliptic equation in {H}\"{o}lder space and the martingale
  problem}.
\bjournal{J. Differential Equations}
\bvolume{314}
\bpages{653--699}.
\bdoi{10.1016/j.jde.2022.01.025}
\bmrnumber{4369182}
\end{barticle}
\endbibitem

\bibitem[\protect\citeauthoryear{Novikov}{1975}]{No75}
\begin{barticle}[author]
\bauthor{\bsnm{Novikov},~\bfnm{A.~A.}\binits{A.~A.}}
(\byear{1975}).
\btitle{On discontinuous martingales}.
\bjournal{Theory Probab. Appl.}
\bvolume{20}
\bpages{11--26}.
\bdoi{https://doi.org/10.1137/1120002}
\end{barticle}
\endbibitem

\bibitem[\protect\citeauthoryear{Perkowski}{}]{Per23}
\begin{barticle}[author]
\bauthor{\bsnm{Perkowski},~\bfnm{N.}\binits{N.}}
\btitle{Energy solutions of singular SPDEs}.
\bjournal{Lecture notes}.
\bdoi{https://www.mi.fu-berlin.de/math/groups/stoch/members/publ_Perkowski/LN202408_MarylandSummerSchool.pdf}
\end{barticle}
\endbibitem

\bibitem[\protect\citeauthoryear{Priola}{2012}]{Pr12}
\begin{barticle}[author]
\bauthor{\bsnm{Priola},~\bfnm{Enrico}\binits{E.}}
(\byear{2012}).
\btitle{Pathwise uniqueness for singular {SDE}s driven by stable processes}.
\bjournal{Osaka J. Math.}
\bvolume{49}
\bpages{421--447}.
\bmrnumber{2945756}
\end{barticle}
\endbibitem

\bibitem[\protect\citeauthoryear{Sato}{1999}]{Sa99}
\begin{bbook}[author]
\bauthor{\bsnm{Sato},~\bfnm{Ken-iti}\binits{K.-i.}}
(\byear{1999}).
\btitle{L\'{e}vy processes and infinitely divisible distributions}.
\bseries{Cambridge Studies in Advanced Mathematics}
\bvolume{68}.
\bpublisher{Cambridge University Press, Cambridge}
\bnote{Translated from the 1990 Japanese original, Revised by the author}.
\bmrnumber{1739520}
\end{bbook}
\endbibitem

\bibitem[\protect\citeauthoryear{Song and Zhang}{2015}]{SZ15}
\begin{barticle}[author]
\bauthor{\bsnm{Song},~\bfnm{Yulin}\binits{Y.}} \AND
  \bauthor{\bsnm{Zhang},~\bfnm{Xicheng}\binits{X.}}
(\byear{2015}).
\btitle{Regularity of density for {SDE}s driven by degenerate {L}\'{e}vy
  noises}.
\bjournal{Electron. J. Probab.}
\bvolume{20}
\bpages{no. 21, 27}.
\bdoi{10.1214/EJP.v20-3287}
\bmrnumber{3325091}
\end{barticle}
\endbibitem

\bibitem[\protect\citeauthoryear{Stroock and Varadhan}{2006}]{SV06}
\begin{bbook}[author]
\bauthor{\bsnm{Stroock},~\bfnm{Daniel~W.}\binits{D.~W.}} \AND
  \bauthor{\bsnm{Varadhan},~\bfnm{S.~R.~Srinivasa}\binits{S.~R.~S.}}
(\byear{2006}).
\btitle{Multidimensional diffusion processes}.
\bseries{Classics in Mathematics}.
\bpublisher{Springer-Verlag, Berlin}
\bnote{Reprint of the 1997 edition}.
\bmrnumber{2190038}
\end{bbook}
\endbibitem

\bibitem[\protect\citeauthoryear{Tanaka, Tsuchiya and Watanabe}{1974}]{TTW74}
\begin{barticle}[author]
\bauthor{\bsnm{Tanaka},~\bfnm{Hiroshi}\binits{H.}},
  \bauthor{\bsnm{Tsuchiya},~\bfnm{Masaaki}\binits{M.}} \AND
  \bauthor{\bsnm{Watanabe},~\bfnm{Shinzo}\binits{S.}}
(\byear{1974}).
\btitle{Perturbation of drift-type for {L}\'evy processes}.
\bjournal{J. Math. Kyoto Univ.}
\bvolume{14}
\bpages{73--92}.
\bdoi{10.1215/kjm/1250523280}
\bmrnumber{368146}
\end{barticle}
\endbibitem

\bibitem[\protect\citeauthoryear{Triebel}{1992}]{Tr92}
\begin{bbook}[author]
\bauthor{\bsnm{Triebel},~\bfnm{Hans}\binits{H.}}
(\byear{1992}).
\btitle{Theory of function spaces. {II}}.
\bseries{Monographs in Mathematics}
\bvolume{84}.
\bpublisher{Birkh\"{a}user Verlag, Basel}.
\bdoi{10.1007/978-3-0346-0419-2}
\bmrnumber{1163193}
\end{bbook}
\endbibitem

\bibitem[\protect\citeauthoryear{Xia et~al.}{2020}]{XXZZ20}
\begin{barticle}[author]
\bauthor{\bsnm{Xia},~\bfnm{Pengcheng}\binits{P.}},
  \bauthor{\bsnm{Xie},~\bfnm{Longjie}\binits{L.}},
  \bauthor{\bsnm{Zhang},~\bfnm{Xicheng}\binits{X.}} \AND
  \bauthor{\bsnm{Zhao},~\bfnm{Guohuan}\binits{G.}}
(\byear{2020}).
\btitle{{$L^q(L^p)$}-theory of stochastic differential equations}.
\bjournal{Stochastic Process. Appl.}
\bvolume{130}
\bpages{5188--5211}.
\bdoi{10.1016/j.spa.2020.03.004}
\bmrnumber{4108486}
\end{barticle}
\endbibitem

\bibitem[\protect\citeauthoryear{Zhang}{2010}]{Zh10}
\begin{barticle}[author]
\bauthor{\bsnm{Zhang},~\bfnm{Xicheng}\binits{X.}}
(\byear{2010}).
\btitle{Stochastic {V}olterra equations in {B}anach spaces and stochastic
  partial differential equation}.
\bjournal{J. Funct. Anal.}
\bvolume{258}
\bpages{1361--1425}.
\bdoi{10.1016/j.jfa.2009.11.006}
\bmrnumber{2565842}
\end{barticle}
\endbibitem

\bibitem[\protect\citeauthoryear{Zhang}{2013a}]{Zh13b}
\begin{barticle}[author]
\bauthor{\bsnm{Zhang},~\bfnm{Xicheng}\binits{X.}}
(\byear{2013}a).
\btitle{Stochastic differential equations with {S}obolev drifts and driven by
  {$\alpha$}-stable processes}.
\bjournal{Ann. Inst. Henri Poincar\'{e} Probab. Stat.}
\bvolume{49}
\bpages{1057--1079}.
\bdoi{10.1214/12-AIHP476}
\bmrnumber{3127913}
\end{barticle}
\endbibitem

\bibitem[\protect\citeauthoryear{Zhang}{2013b}]{Zh13}
\begin{barticle}[author]
\bauthor{\bsnm{Zhang},~\bfnm{Xicheng}\binits{X.}}
(\byear{2013}b).
\btitle{Degenerate irregular {SDE}s with jumps and application to
  integro-differential equations of {F}okker-{P}lanck type}.
\bjournal{Electron. J. Probab.}
\bvolume{18}
\bpages{no. 55, 25}.
\bdoi{10.1214/EJP.v18-2820}
\bmrnumber{3065865}
\end{barticle}
\endbibitem

\bibitem[\protect\citeauthoryear{Zhang and Zhao}{}]{ZZ17}
\begin{barticle}[author]
\bauthor{\bsnm{Zhang},~\bfnm{Xicheng}\binits{X.}} \AND
  \bauthor{\bsnm{Zhao},~\bfnm{Guohuan}\binits{G.}}
\btitle{Heat kernel and ergodicity of SDEs with distributional drifts}.
\bjournal{Available at arXiv:\eprint{1710.10537}}.
\bdoi{https://arxiv.org/abs/1710.10537}
\end{barticle}
\endbibitem

\bibitem[\protect\citeauthoryear{Zhang and Zhao}{2018}]{ZZ18}
\begin{barticle}[author]
\bauthor{\bsnm{Zhang},~\bfnm{Xicheng}\binits{X.}} \AND
  \bauthor{\bsnm{Zhao},~\bfnm{Guohuan}\binits{G.}}
(\byear{2018}).
\btitle{Singular {B}rownian diffusion processes}.
\bjournal{Commun. Math. Stat.}
\bvolume{6}
\bpages{533--581}.
\bdoi{10.1007/s40304-018-0164-7}
\bmrnumber{3877717}
\end{barticle}
\endbibitem

\end{thebibliography}

\end{document}